\documentclass[a4paper,12pt]{amsart}
\usepackage{amssymb,amsmath,a4wide,graphicx,stmaryrd,fullpage,setspace,microtype,textcomp,tikz,enumitem,accents,mathtools,abraces,cite}
\usepackage[pagebackref=true]{hyperref}
\usepackage{float}
\usepackage[all,cmtip]{xy}
\title{Homotopy theory of Moore flows (II)}

\author[P. Gaucher]{Philippe Gaucher}

\address{Universit\'e de Paris, CNRS, IRIF, F-75006, Paris, France}

\urladdr{http://www.irif.fr/{\~{}}gaucher} 

\makeatletter
\@namedef{subjclassname@2020}{%
	\textup{2020} Mathematics Subject Classification}
\makeatother


\swapnumbers

\newcommand{\D}{\mathcal{D}}
\newcommand{\K}{\mathcal{K}}

\newcommand{\de}{\partial}
\newcommand{\p}{\times}
\renewcommand{\vec}{\overrightarrow}
\renewcommand{\P}{\mathbb{P}}

\newtheorem{thm}{Theorem}[section]
\newtheorem{prop}[thm]{Proposition}
\newtheorem{lem}[thm]{Lemma}
\newtheorem{exa}[thm]{Example}
\newtheorem{cor}[thm]{Corollary}

\newtheorem{defn}[thm]{Definition}
\newtheorem{nota}[thm]{Notation}

\newtheorem{defnot}[thm]{Definition and notation}
\newcommand{\bd}{\begin{defn}}
\newcommand{\ed}{\end{defn}}
\newcommand{\bdn}{\begin{defnot}}
\newcommand{\edn}{\end{defnot}}
\newcommand{\bp}{\begin{prop}}
\newcommand{\ep}{\end{prop}}
\newcommand{\bth}{\begin{thm}}
\renewcommand{\eth}{\end{thm}}
\newcommand{\bpf}{\begin{proof}}
\newcommand{\epf}{\end{proof}}
\newcommand{\bc}{\begin{cor}}
\newcommand{\ec}{\end{cor}}

\newcommand{\fL}[1]{\ar@{->}[ll]_-{#1}}
\newcommand{\fR}[1]{\ar@{->}[rr]^-{#1}}
\newcommand{\fRr}[1]{\ar@{->}[rrr]^-{#1}}
\newcommand{\fD}[1]{\ar@{->}[dd]_-{#1}}
\newcommand{\fU}[1]{\ar@{->}[uu]^-{#1}}
\newcommand{\f}[2]{\ar@{->}[#1]|{#2}}
\newcommand{\ff}[2]{\ar@2{->}[#1]|{#2}}
\newcommand{\frr}[1]{\ar@{->}[rrrr]^-{#1}}

\newcommand{\fl}[1]{\ar@{->}[l]_-{#1}}
\newcommand{\fr}[1]{\ar@{->}[r]^-{#1}}
\newcommand{\fd}[1]{\ar@{->}[d]_-{#1}}
\newcommand{\fu}[1]{\ar@{->}[u]^-{#1}}

\renewcommand{\top}{{\mathbf{Top}}}
\newcommand{\ho}{{\mathbf{Ho}}}
\newcommand{\iso}{\cong}

\newcommand{\vI}{\vec{I}}
\renewcommand{\leq}{\leqslant}
\renewcommand{\geq}{\geqslant}

\newcommand{\mdtop}{{\brm{MdTop}}}

\newcommand{\moore}{{\mathbb{M}}}
\newcommand{\lmoore}{\mathbb{M}_!}

\newcommand{\topdgr}{[\mathcal{G}^{op},\top]}

\def\cartesien{%
  \ar@{-}[]+R+<6pt,-2pt>;[]+RD+<6pt,-6pt>%
  \ar@{-}[]+D+<2pt,-6pt>;[]+RD+<6pt,-6pt>%
}
\def\cocartesien{%
  \ar@{-}[]+L+<-6pt,+2pt>;[]+LU+<-6pt,+6pt>%
  \ar@{-}[]+U+<-2pt,+6pt>;[]+LU+<-6pt,+6pt>%
}
\def\hocartesien{%
  \ar@{-}[]+R+<6pt,-2pt>;[]+RD+<6pt,-6pt>_{h}%
  \ar@{-}[]+D+<2pt,-6pt>;[]+RD+<6pt,-6pt>%
}
\def\hococartesien{%
  \ar@{-}[]+L+<-6pt,+2pt>;[]+LU+<-6pt,+6pt>_{h}%
  \ar@{-}[]+U+<-2pt,+6pt>;[]+LU+<-6pt,+6pt>%
}

\newcommand{\brm}[1]{\rm{\mathbf{#1}}}

\newcommand{\dtop}{{\brm{Flow}}}

\newcommand{\dtopG}{{\mathcal{G}\brm{Flow}}}

\newcommand{\set}{{\brm{Set}}}

\newcommand{\ttop}{{\brm{TOP}}}
\newcommand{\mtop}{{\brm{MTop}}}

\newcommand{\glob}{{\rm{Glob}}}

\newcommand{\globP}{{\rm{Glob}}}
\newcommand{\globG}{{\rm{Glob}}^{\mathcal{G}}}

\DeclareMathOperator{\id}{Id}

\DeclareMathOperator{\Obj}{Obj}
\DeclareMathOperator{\Mor}{Mor}
\DeclareMathOperator{\pr}{pr}

\newcommand{\liminj}{\varinjlim}
\newcommand{\limproj}{\varprojlim}

\newcommand{\rest}{\!\upharpoonright\!}

\DeclareMathOperator{\carrier}{Carrier}

\newcommand{\ptop}[1]{{\brm{{#1}dTop}}}

\makeatletter
\def\varholim@#1#2{%
  \vtop{\m@th\ialign{##\cr
    \hfil$#1\operator@font holim$\hfil\cr
    \noalign{\nointerlineskip\kern1.5\ex@}#2\cr
    \noalign{\nointerlineskip\kern-\ex@}\cr}}%
}
\def\holimproj{%
  \mathop{\mathpalette\varholim@{\leftarrowfill@\textstyle}}\nmlimits@
}
\def\holiminj{%
  \mathop{\mathpalette\varholim@{\rightarrowfill@\textstyle}}\nmlimits@
}
\makeatother

\makeatletter
\newcommand*{\@opargbegintheorem}[3]{\trivlist
	\item[\hskip \labelsep{\bfseries #1\ #2}] \textbf{(#3)}\ \itshape}
\makeatother

\DeclareMathOperator{\cocyl}{{Path}}
\newcommand{\adj}[4]{\xymatrix@1{{#1}\ar@/^0.8em/[r]^-{#2} \ar@{}[r]|-{\perp} & \ar@/^0.8em/[l]^-{#3} {#4}}}

\newcommand{\ot}{\otimes}
\newcommand{\mins}[1]{\texorpdfstring{$#1$}{Lg}}

\makeatletter
\newcommand*{\addFileDependency}[1]{
	\typeout{(#1)}
	\@addtofilelist{#1}
	\IfFileExists{#1}{}{\typeout{No file #1.}}
}
\makeatother

\begin{document}

\begin{abstract} 
	This paper proves that the q-model structures of Moore flows and of multipointed $d$-spaces are Quillen equivalent. The main step is the proof that the counit and unit maps of the Quillen adjunction are isomorphisms on the q-cofibrant objects (all objects are q-fibrant). As an application, we provide a new proof of the fact that the categorization functor from multipointed $d$-spaces to flows has a total left derived functor which induces a category equivalence between the homotopy categories. The new proof sheds light on the internal structure of the categorization functor which is neither a left adjoint nor a right adjoint. It is even possible to write an inverse up to homotopy of this functor using Moore flows.
\end{abstract}

\maketitle

\tableofcontents

\section{Introduction}
\subsection*{Presentation}

This paper is the companion paper of \cite{Moore1}. The purpose of these two papers is to exhibit, by means of the q-model
category of \textit{Moore flows} (cf. Definition~\ref{def-Moore-flow}), a zig-zag of Quillen equivalences between the q-model structure of \textit{multipointed $d$-spaces} introduced in \cite{mdtop} and the q-model structure of \textit{flows} introduced in \cite{model3}. The only known functor which was a good candidate for a Quillen equivalence from multipointed $d$-spaces to flows (Definition~\ref{cat-func}) has indeed a total left derived functor in the sense of \cite{HomotopicalCategory} which induces an equivalence of categories between the homotopy categories (\cite[Theorem~7.5]{mdtop}). However, this functor is neither a left adjoint nor a right adjoint by Theorem~\ref{noadjoint}.

Multipointed $d$-spaces and flows can be used to model concurrent processes. For example, the paper \cite{ccsprecub} shows how to model all process algebras for any synchronization algebra using flows. There are many geometric models of concurrency available in the literature \cite{mg,equilogical,distinguishedcube2,distinguishedcube1,MR2545830} (the list does not pretend to be exhaustive). Most of them are used to study the fundamental category of a concurrent process or any derived concept. It is something which can be also carried out with the formalisms of flows and multipointed $d$-spaces. The fundamental category functor is easily calculable indeed, at least for flows since it is a left adjoint, and for cellular multipointed $d$-spaces by using Corollary~\ref{calculation-pathspace3}, and it interacts very well with the underlying simplicial structures.

The Quillen equivalence between flows and Moore flows is proved in \cite[Theorem~10.9]{Moore1}. The Quillen equivalence between multipointed $d$-spaces and Moore flows is proved in Theorem~\ref{adj-zigzag}. The latter theorem is a consequence of the structural properties of the adjunction between multipointed $d$-spaces and Moore flows which can be summarized as follows:

\begin{thm} \label{synth} (Theorem~\ref{c2}, Corollary~\ref{c1} and Theorem~\ref{adj-zigzag})
	The adjunction \[\lmoore^{\mathcal{G}}\dashv \moore^{\mathcal{G}}:\dtopG \leftrightarrows \ptop{\mathcal{G}}\] between Moore flows and multipointed $d$-spaces is a Quillen equivalence. The counit map and the unit map of this Quillen adjunction are isomorphisms on q-cofibrant objects (recall that all objects are q-fibrant).
\end{thm}

Another standard example of the situation of Theorem~\ref{synth} is the Quillen equivalence between the q-model structures of $\Delta$-generated spaces and of $k$-spaces (cf. Appendix~\ref{kspaces}). 

This paper is the first use in a real practical situation of the closed semimonoidal category of $\mathcal{G}$-spaces (Definition~\ref{gspace} and Theorem~\ref{closedsemimonoidal}). It illustrates the interest of this structure already for calculating spaces of execution paths of cellular multipointed $d$-spaces. The interest of this structure is beyond directed homotopy theory, as remarked in the introduction of \cite{Moore1} where possible connections with type theory are briefly discussed. The potential of this semimonoidal structure is visible in the proofs of Proposition~\ref{comp-gl}, Theorem~\ref{pre-calculation-pathspace}, Theorem~\ref{pre-calculation-pathspace2} and Corollary~\ref{calculation-pathspace}. 

The Moore flows enable us to write explicitly an ``inverse up to homotopy'' of the categorization functor of Definition~\ref{cat-func} in Definition~\ref{inv-up-to-homotopy}. Two applications of the existence of this inverse up to homotopy are given. The first one is a new proof of \cite[Theorem~7.5]{mdtop} provided in Theorem~\ref{lepb4} which is totally independent from \cite{model2,mdtop}. The second one is a concise and very natural definition of the underlying homotopy type of a flow in Proposition~\ref{underlyinghomotopytype}. 

As a curiosity, it is also proved in passing a kind of second Dini theorem for spaces of execution paths of finite cellular multipointed $d$-spaces without loops in Corollary~\ref{Dini-cellular}.

This paper is written in the setting of $\Delta$-Hausdorff or not $\Delta$-generated spaces. The setting of weakly Hausdorff or not $k$-spaces is of very little interest for the study of multipointed $d$-spaces and flows not only because all concrete examples coming from computer science are cellular objects of the q-model structures, and also because it is not known how to left Bousfield localize because of the cofibration identifying two states. The locally presentable setting has many other advantages like the existence of adjoints \cite[Theorem~5.10]{leftproperflow}, the smallness conditions of \cite[Definition~2.1.3]{MR99h:55031} always satisfied and the existence of left-determined model categories in the tractable cases \cite{henry2020minimal}. It might be interesting anyway to make some comments about $k$-spaces to emphasize some topological arguments of this paper. These comments are postponed to Appendix~\ref{kspaces}.

\subsection*{Outline of the paper}

\begin{itemize}[leftmargin=*]
	\item Section~\ref{about-multi-space} is a reminder about multipointed $d$-spaces and about their q-model structure. It also contains new results about the topology of the space of execution paths. The section starts with a short reminder about $\Delta$-generated spaces. The notion of \textit{$\Delta$-inclusion} is introduced to clarify some topological arguments: they are for $\Delta$-generated spaces what $k$-inclusions are for $k$-spaces.
	\item The functor $\Omega$ which forgets the set of execution paths of a multipointed $d$-space is topological. Section~\ref{final-by-Moore} gives an explicit description of the $\Omega$-final structure in term of Moore composition. It culminates with Theorem~\ref{final-structure-revisited}. The calculations are a bit laborious but some of them are used further in the paper.
	\item Section~\ref{adj}, after a reminder about Moore flows and their q-model structure, describes the adjunctions between multipointed $d$-spaces and Moore flows. The right adjoint from multipointed $d$-spaces to Moore flows is quite easy to define. The existence of the left adjoint is straightforward. Appendix~\ref{lmoore-explicit} provides an explicit construction of this left adjoint $\lmoore^{\mathcal{G}}:\dtopG\to \ptop{\mathcal{G}}$. It uses results dating back to \cite{model3} obtained for flows, i.e. small semicategories enriched over topological spaces, and adapted to Moore flows, i.e. small semicategories enriched over $\mathcal{G}$-spaces. This explicit construction is not necessary to establish the results of the main part of the paper. It is the reason why it is postponed to an appendix.
	\item Section~\ref{cellular-obj} gathers some geometric properties of cellular multipointed $d$-spaces concerning their underlying topologies, the topologies of their spaces of execution paths and some of their structural properties like Theorem~\ref{bounded0} which has important consequences. The main tools are the notion of \textit{carrier} of an execution path (Definition~\ref{def-carrier}) and the notion of \textit{achronal slice} of a globular cell (Definition~\ref{achronal}) studied in Proposition~\ref{middle0} and Proposition~\ref{good-achronal}. It also contains Theorem~\ref{calcul_final_structure} which provides a kind of normal form for the execution paths of a cellular multipointed $d$-space obtained as a pushout along a generating q-cofibration. 
	\item Section~\ref{chain} studies chains of globes. It is an important geometric object for the proofs of this paper. It enables us to understand what happens \textit{locally} in the space of execution paths of a cellular multipointed $d$-space. The main theorem is Theorem~\ref{img-closed} which can be viewed as a workaround of the fact that the space $\mathcal{G}(1,1)$ of nondecreasing homeomorphisms from $[0,1]$ to itself equipped with the compact-open topology is not sequentially compact. As a byproduct, it is also proved in Corollary~\ref{Dini-cellular} a second Dini theorem for finite cellular multipointed $d$-spaces without loops.
	\item Section~\ref{unit} is the core of the paper. It proves that the unit and the counit of the adjunction are isomorphisms on q-cofibrant objects in Theorem~\ref{c2} and in Corollary~\ref{c1}. The main tool of this part is Corollary~\ref{calculation-pathspace} which proves that the right adjoint constructed in Section~\ref{adj} preserves pushouts of cellular multipointed $d$-spaces along q-cofibrations. It relies on Theorem~\ref{pre-calculation-pathspace} whose proof performs an analysis of the execution paths in a pushout along a q-cofibration and on Theorem~\ref{pre-calculation-pathspace2} whose proof carries out a careful analysis of the underlying topology of the spaces of execution paths involved. 
	\item Section~\ref{conclusion} is the concluding section. It establishes that the adjunction between multipointed $d$-spaces and Moore flows yields a Quillen equivalence between the q-model structures. It provides, as an application, a more conceptual proof of the fact that the categorization functor $cat$ from multipointed $d$-spaces to flows has a total left derived functor which induces a category equivalence between the homotopy categories of the q-model structures of multipointed $d$-spaces and flows. And finally, it is shown how to recover the underlying homotopy type of a flow in a very intuitive way.
\end{itemize}

\subsection*{Prerequisites}

We refer to \cite{TheBook} for locally presentable categories, to \cite{MR2506258} for combinatorial model categories.  We refer to \cite{MR99h:55031,ref_model2} for more general model categories. We refer to \cite{KellyEnriched} and to \cite[Chapter~6]{Borceux2} for enriched categories. All enriched categories are topologically enriched categories: \textit{the word topologically is therefore omitted}. What follows is some notations and conventions.

\begin{itemize}[leftmargin=*]
	\item $A:=B$ means that $B$ is the definition of $A$.
	\item $\iso$ denotes an isomorphism, $\simeq$ denotes a weak equivalence.
	\item $f\rest_A$ denotes the restriction of $f$ to $A$.
	\item $\set$ is the category of sets. 
	\item $\mathcal{T\!O\!P}$ is the category of general topological spaces together with the continuous maps.
	\item $\K^{op}$ denotes the opposite category of $\K$.
	\item $\Obj(\K)$ is the class of objects of $\K$.
	\item $\Mor(\K)$ is the category of morphisms of $\K$ with the commutative squares for the morphisms.
	\item $\K^I$ is the category of functors and natural transformations from a small category $I$ to $\K$.
	\item $\Delta_I(Z)$ is the constant diagram over the small category $I$ with unique value $Z$.
	\item $\varnothing$ is the initial object, $\mathbf{1}$ is the final object, $\id_X$ is the identity of $X$.
	\item $\K(X,Y)$ is the set of maps in a set-enriched, i.e. locally small, category $\K$.
	\item $\K(X,Y)$ is the space of maps in an enriched category $\K$. The underlying set of maps may be denoted by $\K_0(X,Y)$ if it is necessary to specify that we are considering the underlying set. 
	\item The composition of two maps $f:A\to B$ and $g:B\to C$ is denoted by $gf$ or, if it is helpful for the reader, by $g.f$; the composition of two functors is denoted in the same way.
	\item The notations $\ell,\ell',\ell_i,L,\dots$ mean a strictly positive real number unless specified something else.
	\item \textbf{$[\ell,\ell']$ denotes a segment. Unless specified, it is always understood that $\ell<\ell'$.}
	\item A \textit{cellular object} of a combinatorial model category is an object $X$ such that the canonical map $\varnothing\to X$ is a transfinite composition of pushouts of generating cofibrations.
	\item The notation $(-)^{cof}$ denotes a cofibrant replacement functor of a combinatorial model structure; note that all model categories of this paper contain only fibrant objects.
	\item \textbf{A compact space is a quasicompact Hausdorff space.}
	\item A sequentially compact space is a space such that each sequence has a limit point.
	\item The set of rational numbers is denoted by $\mathbb{Q}$, the set of real numbers by $\mathbb{R}$.
	\item The complement of $A\subset B$ is denoted by $A^c$ if there is no ambiguity.
	\item Let $n\geq 1$. Denote by $\mathbf{D}^n = \{b\in \mathbb{R}^n, |b| \leq 1\}$ the $n$-dimensional disk, and by $\mathbf{S}^{n-1}= \{b\in \mathbb{R}^n, |b| = 1\}$ the $(n-1)$-dimensional sphere. By convention, let $\mathbf{D}^{0}=\{0\}$ and $\mathbf{S}^{-1}=\varnothing$. 
\end{itemize}

\subsection*{Acknowledgment}

I am indebted to Tyrone Cutler for drawing my attention to the paper \cite{MR3270173}. I thank the anonymous referee for reading this very technical paper.

\section{Multipointed \mins{d}-spaces}
\label{about-multi-space}

Throughout the paper, we work with the category, denoted by $\top$, either of \textit{$\Delta$-generated spaces} or of \textit{$\Delta$-Hausdorff $\Delta$-generated spaces} (cf. \cite[Section~2 and Appendix~B]{leftproperflow}) equipped with its q-model structure (we use the terminology of \cite{ParamHomTtheory}). We summarize some basic properties of $\top$ for the convenience of the reader: 
\begin{itemize}[leftmargin=*]
	\item $\top$ is locally presentable.
	\item All objects of $\top$ are sequential topological spaces.
	\item A closed subset of a $\Delta$-generated space equipped with the relative topology is not necessarily $\Delta$-generated (e.g. the Cantor set), but it is always sequential.
	\item All locally path-connected first-countable topological spaces are $\Delta$-generated by \cite[Proposition~3.11]{MR3270173}, in particular all locally path-connected metrizable topological spaces are $\Delta$-generated.
	\item The inclusion functor from the full subcategory of $\Delta$-generated spaces to the category of general topological spaces together with the continuous maps has a right adjoint called the $\Delta$-kelleyfication functor. The latter functor does not change the underlying set.
	\item Let $A\subset B$ be a subset of a space $B$ of $\top$. Then $A$ equipped with the $\Delta$-kelleyfication of the relative topology belongs to $\top$. 
	\item The colimit in $\top$ is given by the final topology in the following situations: 
	\begin{itemize}
		\item A transfinite compositions of one-to-one maps.
		\item A pushout along a closed inclusion.
		\item A quotient by a closed subset or by an equivalence relation having a closed graph.
	\end{itemize}
	In these cases, the underlying set of the colimit is therefore the colimit of the underlying sets. In particular, the CW-complexes, and more generally all cellular spaces are equipped with the final topology. 
	\item Cellular spaces are weakly Hausdorff. It implies that the image by any continuous application of any compact is closed and compact, i.e. closed, quasicompact and Hausdorff. Cellular spaces are also $\Delta$-Hausdorff and therefore has unique sequential limits by \cite[Proposition~B.17]{leftproperflow}.
	\item $\top$ is cartesian closed. The internal hom $\ttop(X,Y)$ is given by taking the $\Delta$-kelleyfication of the compact-open topology on the set $\mathcal{T\!O\!P}(X,Y)$ of all continuous maps from $X$ to $Y$. 
\end{itemize}

\bd
A one-to-one map of $\Delta$-generated spaces $i:A\to B$ is a {\rm $\Delta$-inclusion} if for all $\Delta$-generated spaces $Z$, the set map $Z\to A$ is continuous if and only if the composite set map $Z\to A\to B$ is continuous.
\ed

\bp \label{DeltaIncl} Let $i:A\to B$ be a one-to-one continuous map. The following assertions are equivalent:
\begin{enumerate}
	\item $i$ is a $\Delta$-inclusion.
	\item $A$ is homeomorphic to $i(A)$ equipped with the $\Delta$-kelleyfication of the relative topology.
	\item A set map $[0,1]\to A$ is continuous if and only if the composite set map $[0,1]\to A\to B$ is continuous.
\end{enumerate}
\ep

\bpf The proof is similar to the same statement for $k$-inclusions of $k$-spaces.
\epf

\begin{cor} \label{DeltaHomeo2}
A continuous bijection $f:U\to V$ of $\top$ is a homeomorphism if and only if it is a $\Delta$-inclusion.
\end{cor}

\begin{nota}
	The notation $[0,\ell_1]\iso^+ [0,\ell_2]$ for two real numbers $\ell_1,\ell_2>0$ means a nondecreasing homeomorphism from $[0,\ell_1]$ to $[0,\ell_2]$. It takes $0$ to $0$ and $\ell_1$ to $\ell_2$. 
\end{nota}

The enriched small category $\mathcal{G}$ is defined as follows: 
\begin{itemize}[leftmargin=*]
	\item The set of objects is the open interval $]0,\infty[$.
	\item The space $\mathcal{G}(\ell_1,\ell_2)$ is the set $\{[0,\ell_1]\iso^+ [0,\ell_2]\}$ for all $\ell_1,\ell_2>0$ equipped with the $\Delta$-kelleyfication of the relative topology induced by the set inclusion $\mathcal{G}(\ell_1,\ell_2) \subset \ttop([0,\ell_1],[0,\ell_2])$. In other terms, a set map $[0,1]\to \mathcal{G}(\ell_1,\ell_2)$ is continuous if and only if the composite set map $[0,1]\to \mathcal{G}(\ell_1,\ell_2)\subset \ttop([0,\ell_1],[0,\ell_2])$ is continuous. 
	\item For every $\ell_1,\ell_2,\ell_3>0$, the composition map \[\mathcal{G}(\ell_1,\ell_2)\p \mathcal{G}(\ell_2,\ell_3) \to \mathcal{G}(\ell_1,\ell_3)\] is induced by the composition of continuous maps. It induces a continuous map since the composite set map \[\mathcal{G}(\ell_1,\ell_2)\p \mathcal{G}(\ell_2,\ell_3) \to \mathcal{G}(\ell_1,\ell_3) \subset \ttop([0,\ell_1],[0,\ell_3])\] corresponds by the adjunction to the continuous map \[[0,\ell_1]\p \mathcal{G}(\ell_1,\ell_2)\p \mathcal{G}(\ell_2,\ell_3) \to [0,\ell_3]\] which takes $(t,x,y)$ to $y(x(t))$. 
\end{itemize} 

The enriched category $\mathcal{G}$ is an example of a reparametrization category in the sense of \cite[Definition~4.3]{Moore1} which is different from the terminal category. It is introduced in \cite[Proposition~4.9]{Moore1}. Another example is given in \cite[Proposition~4.11]{Moore1}.

\bp \label{morphG-metrizable}
The topology of $\mathcal{G}(\ell_1,\ell_2)$ is the compact-open topology. In particular, it is metrizable. A sequence $(\phi_n)_{n\geq 0}$ of $\mathcal{G}(\ell_1,\ell_2)$ converges to $\phi\in \mathcal{G}(\ell_1,\ell_2)$ if and only if it converges pointwise. 
\ep

It means that the topology of the pointwise convergence of $\mathcal{G}(\ell_1,\ell_2)$ is $\Delta$-generated. Proposition~\ref{morphG-metrizable} has an interesting generalization in Corollary~\ref{Dini-cellular}.

\bpf
The compact-open topology on $\mathcal{G}(\ell_1,\ell_2)$ is metrizable by \cite[Proposition~A.13]{MR1867354}. The metric is given by the distance of the uniform convergence. Consider a ball $B(\phi_,\epsilon)$ for this metric. Let $\psi\in B(\phi_,\epsilon)$. Then for all $h\in [0,1]$,
\[
|\big(h\psi(t) + (1-h)\phi(t)\big) -\phi(t)| = |h(\psi(t)-\phi(t))| <h\epsilon\leq\epsilon.
\]
Thus, the compact-open topology is locally path-connected. The compact-open topology is therefore equal to its $\Delta$-kelleyfication. The last assertion is then a consequence of the second Dini theorem. 
\epf

A \textit{multipointed space} is a pair $(|X|,X^0)$ where
\begin{itemize}[leftmargin=*]
	\item $|X|$ is a topological space called the \textit{underlying space} of $X$.
	\item $X^0$ is a subset of $|X|$ called the \textit{set of states} of $X$.
\end{itemize}
A morphism of multipointed spaces $f:X=(|X|,X^0) \rightarrow Y=(|Y|,Y^0)$ is a commutative square
\[
\xymatrix@C=3em@R=3em{
	X^0 \fr{f^0}\fd{} & Y^0 \fd{} \\ 
	|X| \fr{|f|} & |Y|.}
\] 
The corresponding category is denoted by $\mtop$. 

\begin{nota}
	The maps $f^0$ and $|f|$ will be often denoted by $f$ if there is no possible confusion. 
\end{nota}

We have the well-known proposition: 

\bp \label{Moore-comp} (The Moore composition) Let $U$ be a topological space. Let \[\gamma_i:[0,\ell_i]\to U\] $n$ continuous maps with $1\leq i \leq n$ with $n\geq 1$. Suppose that $\gamma_i(\ell_{i})=\gamma_{i+1}(0)$ for $1\leq i < n$. Then there exists a unique continuous map \[\gamma_1*\dots *\gamma_n:[0,\sum_i\ell_i]\to U\] such that 
\[
(\gamma_1*\dots *\gamma_n)(t) = \gamma_i\big(t-\sum_{j<i}\ell_i\big) \hbox{ for }\sum_{j<i} \ell_i\leq t \leq \sum_{j\leq i} \ell_i.
\]
In particular, there is the equality $(\gamma_1* \gamma_2)* \gamma_3=\gamma_1*(\gamma_2* \gamma_3)$.
\ep

\begin{nota}
	Let $\ell>0$. Let $\mu_{\ell}:[0,\ell]\to [0,1]$ be the homeomorphism defined by $\mu_\ell(t) = t/\ell$. 
\end{nota}

\bd \label{composition_map} The map $\gamma_1*\gamma_2$ is called the {\rm Moore composition} of $\gamma_1$ and $\gamma_2$.  The composite 
\[\xymatrix@C=4em{\gamma_1 *_N \gamma_2: [0,1] \fr{(\mu_2)^{-1}}& [0,2]
	\fr{\gamma_1*\gamma_2}& U}\] 
is called the {\rm normalized composition}. One has 
\[
(\gamma_1 *_N \gamma_2)(t) = \begin{cases}
\gamma_1(2t)& \hbox{ if }0\leq t\leq \frac{1}{2},\\
\gamma_2(2t-1)& \hbox{ if }\frac{1}{2}\leq t\leq 1.
\end{cases}
\]
The normalized composition being not associative, a notation like $\gamma_1 *_N \dots *_N \gamma_n$ will mean, by convention, that $*_N$ is applied from the left to the right. 
\ed

A \textit{multipointed $d$-space $X$} is a triple $(|X|,X^0,\P^{\mathcal{G}}X)$ where
\begin{itemize}[leftmargin=*]
	\item The pair $(|X|,X^0)$ is a multipointed space. The space $|X|$ is called the \textit{underlying space} of $X$ and the set $X^0$ the \textit{set of states} of $X$.
	\item The set $\P^{\mathcal{G}}X$ is a set of continous maps from $[0,1]$ to $|X|$ called the \textit{execution paths}, satisfying the following axioms:
	\begin{itemize}
		\item For any execution path $\gamma$, one has $\gamma(0),\gamma(1)\in X^0$.
		\item Let $\gamma$ be an execution path of $X$. Then any composite $\gamma\phi$ with $\phi\in \mathcal{G}(1,1)$ is an execution path of $X$.
		\item Let $\gamma_1$ and $\gamma_2$ be two composable execution paths of $X$; then the normalized composition $\gamma_1 *_N \gamma_2$ is an execution path of $X$.
	\end{itemize}
\end{itemize}
A map $f:X\to Y$ of multipointed $d$-spaces is a map of multipointed spaces from $(|X|,X^0)$ to $(|Y|,Y^0)$ such that for any execution path $\gamma$ of $X$, the map \[\P^{\mathcal{G}}f:\gamma\mapsto f. \gamma\] is an execution path of $Y$. 

\begin{nota}
	The mapping $\P^{\mathcal{G}}f$ will be often denoted by $f$ if there is no ambiguity.
\end{nota} 

The following examples play an important role in the sequel. 
\begin{enumerate}[leftmargin=*]
	\item Any set $E$ will be identified with the multipointed $d$-space $(E,E,\varnothing)$.
	\item The \textit{topological globe of $Z$ of length $\ell>0$}, which is denoted by $\glob^{\mathcal{G}}_\ell(Z)$, is the multipointed $d$-space defined as follows
	\begin{itemize}
		\item the underlying topological space is the quotient space~\footnote{It is the suspension of $Z$.} \[\frac{\{{0},{1}\}\sqcup (Z\p[0,\ell])}{(z,0)=(z',0)={0},(z,1)=(z',1)={1}}\]
		\item the set of states is $\{{0},{1}\}$
		\item the set of execution paths is the set of continuous maps \[\{\delta_z\phi\mid \phi\in \mathcal{G}(1,\ell),z\in  Z\}\]
		with $\delta_z(t) = (z,t)$.	It is equal to the underlying set of $\mathcal{G}(1,\ell)\p Z$.
	\end{itemize}
	In particular, $\glob^{\mathcal{G}}_\ell(\varnothing)$ is the multipointed $d$-space $\{{0},{1}\} = (\{{0},{1}\},\{{0},{1}\},\varnothing)$. For $\ell=1$, we set \[\glob^{\mathcal{G}}(Z)=\glob^{\mathcal{G}}_1(Z).\]
	\item The \textit{directed segment} is the multipointed $d$-space $\vI^{\mathcal{G}}=\glob^{\mathcal{G}}(\{0\})$. 
\end{enumerate}

The category of multipointed $d$-spaces is denoted by $\ptop{\mathcal{G}}$. The subset of execution paths from $\alpha$ to $\beta$ is the set of $\gamma\in\P^{\mathcal{G}} X$  such that $\gamma(0)=\alpha$ and $\gamma(1)=\beta$; it is denoted by $\P^{\mathcal{G}}_{\alpha,\beta} X$: $\alpha$ is called the \textit{initial state} and $\beta$ the \textit{final state} of such a $\gamma$. An execution path having the same initial and final state is called a \textit{loop}. The set $\P^{\mathcal{G}}_{\alpha,\beta} X$ is equipped with the $\Delta$-kelleyfication of the relative topology induced by the inclusion $\P^{\mathcal{G}}_{\alpha,\beta} X\subset \ttop([0,1],|X|)$. In other terms, a set map $U\to \P^{\mathcal{G}}_{\alpha,\beta} X$ is continuous if and only if the composite set map $U\to \P^{\mathcal{G}}_{\alpha,\beta} X\subset \ttop([0,1],|X|)$ is continuous. The category $\ptop{\mathcal{G}}$ is locally presentable by \cite[Theorem~3.5]{mdtop}.

\bp (\cite[Proposition~6.5]{QHMmodel}) \label{functtop}
The mapping $\Omega:X\mapsto (|X|,X^0)$ induces a functor from $\ptop{\mathcal{G}}$ to $\mtop$ which is topological and fibre-small. 
\ep

The $\Omega$-final structure is generated by the finite normalized composition of execution paths. We will come back on this point in Theorem~\ref{final-structure-revisited}. Note that Proposition~\ref{functtop} holds both by working with $\Delta$-generated spaces and with $\Delta$-Hausdorff $\Delta$-generated spaces.

The following proposition is implicitly assumed (for $\ell=1$) in all the previous papers about multipointed $d$-spaces: 

\bp \label{calcul-topology-glob}
Let $Z$ be a topological space. Then there is the homeomorphism \[\P_{0,1}^{\mathcal{G}}\glob_\ell^{\mathcal{G}}(Z)\iso \mathcal{G}(1,\ell)\p Z.\]
\ep

\bpf
The set map 
\[
\begin{cases}
&\Psi:\mathcal{G}(1,\ell)\p Z \longrightarrow \P_{0,1}^{\mathcal{G}}\glob_\ell^{\mathcal{G}}(Z)\\
&(\phi,z) \mapsto \delta_{z}\phi
\end{cases}
\]
is continuous because the mapping $(t,\phi,z)\mapsto (z,\phi(t))$ from $[0,1]\p \mathcal{G}(1,\ell)\p Z$ to $|\glob_\ell^{\mathcal{G}}(Z)|$ is continuous. It is a bijection since, by definition of $\glob_\ell^{\mathcal{G}}(Z)$, the underlying set of $\P_{0,1}^{\mathcal{G}}\glob_\ell^{\mathcal{G}}(Z)$ is equal to the underlying set of $\mathcal{G}(1,\ell)\p Z$. The composite set map 
\[
\begin{cases}
& \P_{0,1}^{\mathcal{G}}\glob_\ell^{\mathcal{G}}(Z)\longrightarrow (\P_{0,1}^{\mathcal{G}}\glob_\ell^{\mathcal{G}}(Z))_{co} \longrightarrow Z\p ]0,1[ \stackrel{\pr_1}\longrightarrow Z\\
& \gamma \mapsto \pr_1(\gamma(\frac{1}{2}))
\end{cases}
\]
where $(\P_{0,1}^{\mathcal{G}}\glob_\ell^{\mathcal{G}}(Z))_{co}$ is the set $\P_{0,1}^{\mathcal{G}}\glob_\ell^{\mathcal{G}}(Z)$ equipped with the compact-open topology is continuous. The continuous map $Z\to \{0\}$ induces a continuous map 
\[
\begin{cases}
& \P_{0,1}^{\mathcal{G}}\glob_\ell^{\mathcal{G}}(Z)\longrightarrow \P_{0,1}^{\mathcal{G}}\glob_\ell^{\mathcal{G}}(\{0\}) \iso \mathcal{G}(1,\ell) \\
& \gamma \mapsto p.\gamma,
\end{cases}
\]
where $p:|\glob_\ell^{\mathcal{G}}(Z)|\to [0,1]$ is the projection map. Therefore the set map 
\[
\begin{cases}
& \Psi^{-1}:\P_{0,1}^{\mathcal{G}}\glob_\ell^{\mathcal{G}}(Z)\longrightarrow \mathcal{G}(1,\ell) \p Z \\
& \gamma \mapsto (p.\gamma,\pr_1(\gamma(\frac{1}{2})))
\end{cases}
\]
is continuous and $\Psi$ is a homeomorphism. 
\epf

\bd
Let $X$ be a multipointed $d$-space $X$. Denote again by $\P^{\mathcal{G}} X$ the topological space \[\P^{\mathcal{G}} X = \bigsqcup_{(\alpha,\beta)\in X^0\p X^0} \P_{\alpha,\beta}^{\mathcal{G}}X.\]
\ed

A straightforward consequence of the definition of the topology of $\P^{\mathcal{G}} X$ is: 

\bp\label{init-final-cst}
Let $X$ be a multipointed $d$-space. Let $f:[0,1] \to \P^{\mathcal{G}}X$ be a continuous map. Then $f$ factors as composite of continuous maps  $f:[0,1] \to \P_{\alpha,\beta}^{\mathcal{G}}X\to \P^{\mathcal{G}}X$ for some $(\alpha,\beta)\in X^0\p X^0$.
\ep

\bpf
It is due to the fact that $[0,1]$ is connected. 
\epf

\bp \label{init-final-cst-cellular} Let $X$ be a multipointed $d$-space such that $X^0$ is a totally disconnected subset of $|X|$. Then the topology of $\P^{\mathcal{G}}X$ is the $\Delta$-kelleyfication of the relative topology induced by the inclusion $\P^{\mathcal{G}}X \subset \ttop([0,1],|X|)$.
\ep

\bpf
Call for this proof $(\P^{\mathcal{G}}X)_+$ the set $\P^{\mathcal{G}}X$ equipped with the $\Delta$-kelleyfication of the relative topology induced by the inclusion $\P^{\mathcal{G}}X_\lambda \subset \ttop([0,1],|X|)$. There is a continuous bijection $\P^{\mathcal{G}}X \to (\P^{\mathcal{G}}X)_+$. Using Corollary~\ref{DeltaHomeo2}, the proof is complete since $X^0$ a totally disconnected subset of $|X|$ and since $[0,1]$ is connected. 
\epf

\bth \label{pre-right-adj}
The functor $\P^{\mathcal{G}}:\mdtop \to \top$ is a right adjoint. In particular, it is limit preserving and accessible.
\eth

\bpf
The left adjoint is constructed in \cite[Proposition~4.9]{mdtop} in the case of $\Delta$-generated spaces. The proof still holds for $\Delta$-Hausdorff $\Delta$-generated spaces. It relies on the fact that $\top$ is cartesian closed and that every $\Delta$-generated space is homeomorphic to the disjoint sum of its path-connected components which are also its connected components. The construction is similar to the construction of the left adjoint of the path $\mathcal{P}$-space functor for $\mathcal{P}$-flows \cite[Theorem~6.13]{Moore1} and to the construction of the left adjoint of the path functor for flows \cite[Theorem~5.9]{leftproperflow}.
\epf

The \textit{q-model structure} of multipointed $d$-spaces is the unique combinatorial model structure such that 
\[\{\glob^{\mathcal{G}}(\mathbf{S}^{n-1})\subset \glob^{\mathcal{G}}(\mathbf{D}^{n}) \mid n\geq 0\} \cup \{C:\varnothing \to \{0\},R:\{0,1\} \to \{0\}\}\]
is the set of generating cofibrations, the maps between globes being induced by the closed inclusion $\mathbf{S}^{n-1}\subset \mathbf{D}^{n}$, and such that 
\[
\{\glob^{\mathcal{G}}(\mathbf{D}^{n}\p\{0\})\subset \glob^{\mathcal{G}}(\mathbf{D}^{n+1}) \mid n\geq 0\}
\]
is the set of generating trivial cofibrations, the maps between globes being induced by the closed inclusion $(x_1,\dots,x_n)\mapsto (x_1,\dots,x_n,0)$ (e.g. \cite[Theorem~6.16]{QHMmodel}). The weak equivalences are the maps of multipointed $d$-spaces $f:X\to Y$  inducing a bijection $f^0:X^0\iso Y^0$ and a weak homotopy equivalence $\P^{{\mathcal{G}}} f:\P^{{\mathcal{G}}} X \to \P^{{\mathcal{G}}} Y$ and the fibrations are the maps of multipointed $d$-spaces $f:X\to Y$  inducing a q-fibration $\P^{{\mathcal{G}}} f:\P^{{\mathcal{G}}} X \to \P^{{\mathcal{G}}} Y$ of topological spaces.

\section{Moore composition and \mins{\Omega}-final structure}
\label{final-by-Moore}

\begin{nota}
	Let $\phi_i:[0,\ell_i] \iso^+ [0,\ell'_i]$ for $n\geq 1$ and $1\leq i \leq n$. Then the map
	\[
	\phi_1 \ot \dots \ot \phi_n : [0,\sum_i \ell_i] \iso^+ [0,\sum_i \ell'_i]
	\]
	denotes the homeomorphism defined by 
	\[
	(\phi_1 \ot \dots \ot \phi_n)(t) = \begin{cases}
	\phi_1(t) & \hbox{if } 0\leq t\leq \ell_1\\
	\phi_2(t-\ell_1)+\ell'_1 & \hbox{if } \ell_1\leq t\leq \ell_1+\ell_2\\
	\dots \\
	\phi_i(t-\sum_{j<i}\ell_j) + \sum_{j<i}\ell'_j& \hbox{if } \sum_{j<i}\ell_j\leq t \leq \sum_{j\leq i}\ell_j\\
	\dots\\
	\phi_n(t-\sum_{j<n}\ell_j) + \sum_{j<n}\ell'_j & \hbox{if } \sum_{j<n}\ell_j\leq t \leq \sum_{j\leq n}\ell_j.
	\end{cases}
	\] 
\end{nota}

\bp \label{decomposition-tenseur}
Let $\phi:[0,\ell] \iso^+ [0,\ell']$. Let $n\geq 1$. Consider $\ell_1,\dots,\ell_n>0$ with $n\geq 1$ such that $\sum_{i=1}^{i=n}\ell_i =\ell$. Then there exists a unique decomposition of $\phi$ of the form \[\phi=\phi_1\ot \dots \ot \phi_n\] such that $\phi_i:[0,\ell_i] \iso^+ [0,\ell'_i]$ for $1\leq i \leq n$.  
\ep

\bpf 
By definition of $\phi_1\ot \dots \ot \phi_n$, we necessarily have 
\[
\phi\big(\sum_{j\leq i}\ell_j\big) = \phi_i(\sum_{j\leq i}\ell_j-\sum_{j<i}\ell_j) + \sum_{j<i}\ell'_j = \phi_i(\ell_i) + \sum_{j<i}\ell'_j = \sum_{j\leq i}\ell'_j
\]
The real numbers $\ell'_i$ are therefore defined by induction on $i\geq 1$ by the formula 
\[
\phi\big(\sum_{j\leq i}\ell_j\big) - \sum_{j< i}\ell'_j = \ell'_i.
\]
In other terms, we have 
\[
\forall 1\leq i\leq n,\phi\big(\sum_{j\leq i}\ell_j\big) = \sum_{j\leq i}\ell'_j.
\]
Let \[\phi_i(t) = \phi(t + \sum_{j<i}\ell_j) - \sum_{j<i}\ell'_j\] for all $t\in [0,\ell_i]$. Then, by definition of $\phi_1 \ot \dots \ot \phi_n$, we obtain 
\[\begin{aligned}
(\phi_1 \ot \dots \ot \phi_n)(t) &=\phi_i\big(t-\sum_{j<i}\ell_j\big) + \sum_{j<i}\ell'_j & \hbox{if } \sum_{j<i}\ell_j\leq t \leq \sum_{j\leq i}\ell_j\\
&= \phi\big(t  -\sum_{j<i}\ell_j + \sum_{j<i}\ell_j\big) - \sum_{j<i}\ell'_j+ \sum_{j<i}\ell'_j& \hbox{if } \sum_{j<i}\ell_j\leq t \leq \sum_{j\leq i}\ell_j\\
&=\phi(t) & \hbox{if } \sum_{j<i}\ell_j\leq t \leq \sum_{j\leq i}\ell_j,
\end{aligned}\]
the first equality by definition of $\phi_1 \ot \dots \ot \phi_n$, the second equality by definition of $\phi_i$ and the third equality by algebraic simplification. Consider a second decomposition \[\phi=\phi'_1\ot \dots \ot \phi'_n\] such that $\phi'_i:[0,\ell_i] \iso^+ [0,\ell'_i]$ for $1\leq i \leq n$. Then for $\sum_{j<i}\ell_j\leq t \leq \sum_{j\leq i}\ell_j$, we have 
\[
\phi_i(t-\sum_{j<i}\ell_j) + \sum_{j<i}\ell'_j = \phi(t) = \phi'_i(t-\sum_{j<i}\ell_j) + \sum_{j<i}\ell'_j
\]
by definition of $\phi_1 \ot \dots \ot \phi_n$ and of $\phi'_1 \ot \dots \ot \phi'_n$. We deduce that $\phi_i=\phi'_i$ for all $1\leq i \leq n$. 
\epf

\begin{cor} \label{decomposition-tenseur2}
	Let $\phi\in \mathcal{G}(1,1)$. Let $n\geq 1$. Assume that \[\sum_{i=1}^{i=n}\ell_i = \sum_{i=1}^{i=n}\ell'_i = 1\] and that 
	\[
	\forall 1\leq i\leq n,\phi\big(\sum_{j\leq i}\ell_j\big) = \sum_{j\leq i}\ell'_j.
	\]
	Then there exist (unique) $\phi_i:[0,\ell_i]\iso^+ [0,\ell'_i]$ for $1\leq i \leq n$ such that $\phi = \phi_1 \ot \dots \ot \phi_n$. 
\end{cor}

\bp \label{lem-1}
Let $U$ be a topological space. Let $\gamma_i:[0,1]\to U$ be $n$ continuous maps with $1\leq i \leq n$ and $n\geq 1$. Let $\phi_i:[0,\ell_i] \iso^+ [0,\ell'_i]$ for $1\leq i \leq n$. Then we have 
\[
\left((\gamma_1\mu_{\ell'_1})*\dots *(\gamma_n\mu_{\ell'_n})\right)(\phi_1\ot\dots\ot\phi_n) = (\gamma_1\mu_{\ell'_1}\phi_1)*\dots *(\gamma_n\mu_{\ell'_n}\phi_n).
\]
\ep

\bpf 
For $\sum_{j<i}\ell_j\leq t \leq \sum_{j\leq i}\ell_j$, we have 
\[\begin{aligned}
\big((\gamma_1\mu_{\ell'_1})*\dots *(\gamma_n\mu_{\ell'_n})\big)(\phi_1 \ot \dots \phi_n)(t) &= (\gamma_i\mu_{\ell'_i})\big((\phi_1 \ot \dots \ot \phi_n)(t)- \sum_{j<i}\ell'_j\big)\\
&= (\gamma_i\mu_{\ell'_i})\big(\big(\phi_i\bigl(t-\sum_{j<i}\ell_j\big) + \sum_{j<i}\ell'_j\big)- \sum_{j<i}\ell'_j\big)\\
&= (\gamma_i\mu_{\ell'_i})\big(\phi_i\bigl(t-\sum_{j<i}\ell_j\big)\big)\\
& = \big((\gamma_1\mu_{\ell'_1}\phi_1)*\dots *(\gamma_n\mu_{\ell'_n}\phi_n)\big)(t),
\end{aligned}\]
the first and the fourth equality by definition of the Moore composition, the second equality by definition of $\phi_1 \ot \dots \phi_n$, and the third equality by algebraic simplification. 
\epf

\bp \label{lem0}
Let $U$ be a topological space. Let $\gamma_i:[0,1]\to U$ be $n$ continuous maps with $1\leq i \leq n$ and $n\geq 1$. Let $\ell_i>0$ with $1\leq i \leq n$ nonzero real numbers with $\sum_i \ell_i = 1$. Then for all $\ell >0$, we have 
\[
\big((\gamma_1\mu_{\ell_1})*\dots *(\gamma_n\mu_{\ell_n})\big)\mu_{\ell} = (\gamma_1\mu_{\ell_1\ell})*\dots *(\gamma_n\mu_{\ell_n\ell}).
\]
\ep

\bpf
For all $1\leq j\leq n$, we have by definition of the Moore composition
\[
\big((\gamma_1\mu_{\ell_1})*\dots *(\gamma_n\mu_{\ell_n})\big)\mu_{\ell}(t) = \gamma_j\bigg(\frac{1}{\ell_j}\bigg(\frac{t}{\ell}-\sum_{i<j} \ell_i\bigg)\bigg)
\]
if $\sum_{i<j} \ell_i\leq t/\ell\leq\sum_{i\leq j} \ell_i$ and still by definition of the Moore composition, we have 
\[
\big((\gamma_1\mu_{\ell_1\ell})*\dots *(\gamma_n\mu_{\ell_n\ell})\big)(t) = \gamma_j\bigg(\frac{t-\sum_{i<j} \ell_i\ell}{\ell_j\ell}\bigg)
\]
if $\sum_{i<j} \ell_i\ell\leq t\leq\sum_{i\leq j} \ell_i\ell$. 
\epf

\bp \label{starN}
Let $U$ be a topological space. Let $\gamma_i:[0,1]\to U$ be $n$ continuous maps with $n\geq 2$ and $1\leq i \leq n$ such that $\gamma_1 *_N \dots *_N \gamma_n$ exists. Then there is the equality 
\[
\gamma_1 *_N \dots *_N \gamma_n = \big(\gamma_1\mu_{\frac{1}{2^{n-1}}}\big) * \big(\gamma_2\mu_{\frac{1}{2^{n-1}}}\big) * \big(\gamma_3\mu_{\frac{1}{2^{n-2}}}\big) * \dots * \big(\gamma_n\mu_{\frac{1}{2^{}}}\big).
\]
In particular, for $n=2$, we have $\gamma_1 *_N \gamma_2 = (\gamma_1\mu_{\frac{1}{2}})  * (\gamma_2 \mu_{\frac{1}{2}})$. 
\ep

\bpf
The proof is by induction on $n\geq 2$. The map $\mu_{\frac{1}{2}}:[0,\frac{1}{2}] \iso^+ [0,1]$ which takes $t$ to $2t$ gives rise to a homeomorphism $\mu_{\frac{1}{2}} \ot \mu_{\frac{1}{2}}:[0,1] \iso^+ [0,2]$ which is equal to $\mu^{-1}_2:[0,1] \iso^+[0,2]$. We then write 
\[\begin{aligned}
\gamma_1 *_N \gamma_2 & = (\gamma_1 * \gamma_2)\mu^{-1}_2 &\hbox{ by definition of $*_N$}\\
& = (\gamma_1 * \gamma_2)(\mu_{\frac{1}{2}} \ot \mu_{\frac{1}{2}})&\hbox{ because $\mu^{-1}_2=\mu_{\frac{1}{2}} \ot \mu_{\frac{1}{2}}$}\\
& = (\gamma_1\mu_{\frac{1}{2}})  * (\gamma_2 \mu_{\frac{1}{2}}) & \hbox{ by Proposition~\ref{lem-1}}.
\end{aligned}\]
The statement is therefore proved for $n=2$. Assume that the statement is proved for some $n\geq 2$ and for $n=2$. Then we obtain
\[\begin{aligned}
\gamma_1 &*_N \dots *_N\gamma_{n+1} \\
&= \big(\big(\gamma_1\mu_{\frac{1}{2^{n-1}}}\big) * \big(\gamma_2\mu_{\frac{1}{2^{n-1}}}\big) * \big(\gamma_3\mu_{\frac{1}{2^{n-2}}}\big) * \dots * \big(\gamma_n\mu_{\frac{1}{2^{}}}\big)\big)  *_N \gamma_{n+1} \\
&= \big(\big(\big(\gamma_1\mu_{\frac{1}{2^{n-1}}}\big) * \big(\gamma_2\mu_{\frac{1}{2^{n-1}}}\big) * \big(\gamma_3\mu_{\frac{1}{2^{n-2}}}\big) * \dots * \big(\gamma_n\mu_{\frac{1}{2^{}}}\big)\big)\mu_{\frac{1}{2^{}}}\big)  * \big(\gamma_{n+1}\mu_{\frac{1}{2^{}}}\big)\\
&= \big(\big(\gamma_1\mu_{\frac{1}{2^{n}}}\big) * \big(\gamma_2\mu_{\frac{1}{2^{n}}}\big) * \big(\gamma_3\mu_{\frac{1}{2^{n-1}}}\big) * \dots * \big(\gamma_n\mu_{\frac{1}{2^{2}}}\big)\big)* \big(\gamma_{n+1}\mu_{\frac{1}{2^{}}}\big)\\
&= \big(\gamma_1\mu_{\frac{1}{2^{n}}}\big) * \big(\gamma_2\mu_{\frac{1}{2^{n}}}\big) * \big(\gamma_3\mu_{\frac{1}{2^{n-1}}}\big) * \dots * \big(\gamma_n\mu_{\frac{1}{2^{2}}}\big)* \big(\gamma_{n+1}\mu_{\frac{1}{2^{}}}\big),
\end{aligned}\]
the first equality by induction hypothesis, the second equality by the case $n=2$, the third equality by Proposition~\ref{lem0}, and the last equality by associativity of the Moore composition. We have proved the statement for $n+1$. 
\epf

\bp \label{startN/2}
Let $U$ be a topological space. Let $\gamma_i:[0,1]\to U$ be $n$ continuous maps with $n\geq 1$ and $1\leq i \leq n$ such that $\gamma_1 *_N \dots *_N \gamma_n$ exists. Let $\phi\in \mathcal{G}(1,1)$. Then there exist $\phi_1:[0,\ell_1]\iso^+ [0,\frac{1}{2^{n-1}}]$, $\phi_2:[0,\ell_2]\iso^+ [0,\frac{1}{2^{n-1}}]$, $\phi_3:[0,\ell_3]\iso^+ [0,\frac{1}{2^{n-2}}]$, etc... until $\phi_n:[0,\ell_n]\iso^+ [0,\frac{1}{2^{}}]$ such that $\phi=\phi_1 \ot \dots \ot \phi_n$ (which implies $\sum_i \ell_i = 1$) and there is the equality 
\[
\big(\gamma_1 *_N \dots *_N \gamma_n\big)\phi = \big(\gamma_1\mu_{\frac{1}{2^{n-1}}}\phi_1\big) * \big(\gamma_2\mu_{\frac{1}{2^{n-1}}}\phi_2\big) * \big(\gamma_3\mu_{\frac{1}{2^{n-2}}}\phi_3\big) * \dots * \big(\gamma_n\mu_{\frac{1}{2^{}}}\phi_n\big).
\]
\ep

\bpf Let $\ell_1,\dots,\ell_n >0$ such that $\sum_i \ell_i = 1$ and such that 
\[
\begin{cases}
\phi(\ell_1)=\frac{1}{2^{n-1}} \\
\phi(\ell_1+\ell_2) = \frac{1}{2^{n-1}} + \frac{1}{2^{n-1}} \\
\phi(\ell_1+\ell_2 + \ell_3) = \frac{1}{2^{n-1}} + \frac{1}{2^{n-1}} + \frac{1}{2^{n-2}} \\
\dots \\
\phi(\ell_1+\ell_2 + \ell_3 + \dots + \ell_n) = \frac{1}{2^{n-1}} + \frac{1}{2^{n-1}} + \frac{1}{2^{n-2}} + \dots + \frac{1}{2^{}} =1.
\end{cases}
\]
By Proposition~\ref{decomposition-tenseur}, there exist $\phi_1:[0,\ell_1]\iso^+ [0,\frac{1}{2^{n-1}}]$, $\phi_2:[0,\ell_2]\iso^+ [0,\frac{1}{2^{n-1}}]$, $\phi_3:[0,\ell_3]\iso^+ [0,\frac{1}{2^{n-2}}]$, etc... until $\phi_n:[0,\ell_n]\iso^+ [0,\frac{1}{2^{}}]$ such that $\phi=\phi_1 \ot \dots \ot \phi_n$. We obtain 
\[\begin{aligned}
\big(\gamma_1 &*_N \dots *_N\gamma_{n+1}\big)\phi \\
&= \big(\big(\gamma_1\mu_{\frac{1}{2^{n}}}\big) * \big(\gamma_2\mu_{\frac{1}{2^{n}}}\big) * \big(\gamma_3\mu_{\frac{1}{2^{n-1}}}\big) * \dots * \big(\gamma_n\mu_{\frac{1}{2^{2}}}\big)* \big(\gamma_{n+1}\mu_{\frac{1}{2^{}}}\big)\big)\phi\\
&= \big(\gamma_1\mu_{\frac{1}{2^{n-1}}}\phi_1\big) * \big(\gamma_2\mu_{\frac{1}{2^{n-1}}}\phi_2\big) * \big(\gamma_3\mu_{\frac{1}{2^{n-2}}}\phi_3\big) * \dots * \big(\gamma_n\mu_{\frac{1}{2^{}}}\phi_n\big),
\end{aligned}\]
the first equality by Proposition~\ref{starN} and the second equality by Proposition~\ref{lem-1}. 
\epf

\bp \label{startN2}
Let $U$ be a topological space. Let $\gamma_i:[0,1]\to U$ be $n$ continuous maps with $n\geq 2$ and $1\leq i \leq n$ such that $\gamma_1 *_N \dots *_N \gamma_n$ exists. Let $\ell_1,\dots,\ell_n >0 $ be nonzero real numbers such that $\sum_i \ell_i = 1$. Let $\phi_1:[0,\frac{1}{2^{n-1}}] \iso^+ [0,\ell_1]$, $\phi_2:[0,\frac{1}{2^{n-1}}] \iso^+ [0,\ell_2]$, $\phi_3:[0,\frac{1}{2^{n-2}}] \iso^+ [0,\ell_3]$, etc... until $\phi_n:[0,\frac{1}{2^{}}] \iso^+ [0,\ell_n]$ and let $\phi=\phi_1 \ot \dots \ot \phi_n$. Then $\phi\in \mathcal{G}(1,1)$ and there is the equality 
\begin{multline*}
\big(\big(\gamma_1\mu_{\ell_1}\big) * \big(\gamma_2\mu_{\ell_2}\big) * \big(\gamma_3\mu_{\ell_3}\big) * \dots * \big(\gamma_n\mu_{\ell_n}\big)\big)\phi = \\
\big(\gamma_1\mu_{\ell_1}\phi_1\mu^{-1}_{\frac{1}{2^{n-1}}}\big) *_N \big(\gamma_2\mu_{\ell_2}\phi_2\mu^{-1}_{\frac{1}{2^{n-1}}}\big) *_N \big(\gamma_3\mu_{\ell_3}\phi_3\mu^{-1}_{\frac{1}{2^{n-2}}}\big) *_N \dots *_N \big(\gamma_n\mu_{\ell_n}\phi_n\mu^{-1}_{\frac{1}{2^{}}}\big).
\end{multline*}
\ep

\bpf  
We have 
\[\begin{aligned}
&\big(\big(\gamma_1\mu_{\ell_1}\big) * \big(\gamma_2\mu_{\ell_2}\big) * \big(\gamma_3\mu_{\ell_3}\big) * \dots * \big(\gamma_n\mu_{\ell_n}\big)\big)\phi \\
&=\big(\gamma_1\mu_{\ell_1}\phi_1\big) * \big(\gamma_2\mu_{\ell_2}\phi_2\big) * \big(\gamma_3\mu_{\ell_3}\phi_3\big) * \dots * \big(\gamma_n\mu_{\ell_n}\phi_n\big)\\
&=\big(\gamma_1\mu_{\ell_1}\phi_1\mu^{-1}_{\frac{1}{2^{n-1}}}\mu_{\frac{1}{2^{n-1}}}\big) * \big(\gamma_2\mu_{\ell_2}\phi_2\mu^{-1}_{\frac{1}{2^{n-1}}}\mu_{\frac{1}{2^{n-1}}}\big) \\&\hspace{5cm}* \big(\gamma_3\mu_{\ell_3}\phi_3\mu^{-1}_{\frac{1}{2^{n-2}}}\mu_{\frac{1}{2^{n-2}}}\big) * \dots * \big(\gamma_n\mu_{\ell_n}\phi_n\mu^{-1}_{\frac{1}{2^{}}}\mu_{\frac{1}{2^{}}}\big)\\
&=\big(\gamma_1\mu_{\ell_1}\phi_1\mu^{-1}_{\frac{1}{2^{n-1}}}\big) *_N \big(\gamma_2\mu_{\ell_2}\phi_2\mu^{-1}_{\frac{1}{2^{n-1}}}\big) *_N \big(\gamma_3\mu_{\ell_3}\phi_3\mu^{-1}_{\frac{1}{2^{n-2}}}\big) *_N \dots *_N \big(\gamma_n\mu_{\ell_n}\phi_n\mu^{-1}_{\frac{1}{2^{}}}\big),
\end{aligned}\]
where the first equality is due to Proposition~\ref{lem-1}, the second equality is due to the fact that $\mu^{-1}_\ell\mu_\ell$ is the identity of $[0,\ell]$ for all nonzero real numbers $\ell >0$, and the last equality is a consequence of Proposition~\ref{starN}. 
\epf

\bth \label{final-structure-revisited}
Consider a cocone $(\Omega(X_i))\stackrel{\bullet}\to (|X|,X^0)$ of $\mtop$. Let $X$ be the $\Omega$-final lift. Let $f_i:X_i\to X$ be the canonical maps. Then the set of execution paths of $X$ is the set of finite Moore compositions of the form $(f_1\gamma_1\mu_{\ell_1}) * \dots * (f_n\gamma_n\mu_{\ell_n})$ such that $\gamma_i$ is an execution path of $X_i$ for all $1\leq i \leq n$ with $\sum_i \ell_i = 1$.
\eth

\bpf Let $\mathcal{P}(X)$ be the set of execution paths of $X$ of the form $(f_1\gamma_1\mu_{\ell_1}) * \dots * (f_n\gamma_n\mu_{\ell_n})$ such that $\gamma_i$ is an execution path of $X_i$ for all $1\leq i \leq n$ with $\sum_i \ell_i = 1$. The final structure is generated by the finite normalized composition of execution paths $(f_1\gamma_1) *_N \dots *_N (f_n\gamma_n)$ (with the convention that the $*_N$ are calculated from the left to the right) and all reparametrizations by $\phi$ running over $\mathcal{G}(1,1)$. By Proposition~\ref{startN/2}, there exist $\phi_1:[0,\ell_1]\iso^+ [0,\frac{1}{2^{n-1}}]$, $\phi_2:[0,\ell_2]\iso^+ [0,\frac{1}{2^{n-1}}]$, $\phi_3:[0,\ell_3]\iso^+ [0,\frac{1}{2^{n-2}}]$, etc... until $\phi_n:[0,\ell_n]\iso^+ [0,\frac{1}{2^{}}]$ such that $\phi=\phi_1 \ot \dots \ot \phi_n$ and we have 
\[\begin{aligned}
\big((f_1\gamma_1)  & *_N \dots  *_N (f_n\gamma_n)\big)\phi \\& = \big(f_1\gamma_1\mu_{\frac{1}{2^{n-1}}}\phi_1\big) * \big(f_2\gamma_2\mu_{\frac{1}{2^{n-1}}}\phi_2\big) * \big(f_3\gamma_3\mu_{\frac{1}{2^{n-2}}}\phi_3\big) * \dots * \big(f_n\gamma_n\mu_{\frac{1}{2^{}}}\phi_n\big)\\
& = \big(f_1\gamma_1\mu_{\frac{1}{2^{n-1}}}\phi_1\mu^{-1}_{\ell_1}\mu_{\ell_1}\big) * \big(f_2\gamma_2\mu_{\frac{1}{2^{n-1}}}\phi_2\mu^{-1}_{\ell_2}\mu_{\ell_2}\big) \\ &\hspace{3cm}* \big(f_3\gamma_3\mu_{\frac{1}{2^{n-2}}}\phi_3\mu^{-1}_{\ell_3}\mu_{\ell_3}\big) * \dots * \big(f_n\gamma_n\mu_{\frac{1}{2^{}}}\phi_n\mu^{-1}_{\ell_n}\mu_{\ell_n}\big)\\
&=\big(f_1\gamma'_1\mu_{\ell_1}\big) * \big(f_2\gamma'_2\mu_{\ell_2}\big) * \big(f_3\gamma'_3\mu_{\ell_3}\big) * \dots * \big(f_n\gamma'_n\mu_{\ell_n}\big),
\end{aligned}\]
the first equality by Proposition~\ref{startN/2}, the second equality because $\mu^{-1}_\ell\mu_\ell$ is the identity of $[0,\ell]$ for all $\ell >0$ and the third equality because of the following notations: 
\[
\begin{cases}
\gamma'_1 = \gamma_1\big(\mu_{\frac{1}{2^{n-1}}}\phi_1\mu^{-1}_{\ell_1}\big)\\
\gamma'_2 = \gamma_2\big(\mu_{\frac{1}{2^{n-1}}}\phi_1\mu^{-1}_{\ell_2}\big)\\
\gamma'_3 = \gamma_3\big(\mu_{\frac{1}{2^{n-2}}}\phi_3\mu^{-1}_{\ell_3}\big)\\
\dots\\
\gamma'_n = \gamma_n\big(\mu_{\frac{1}{2^{}}}\phi_n\mu^{-1}_{\ell_n}\big).
\end{cases}
\]
It implies that the set $\mathcal{P}(X)$ contains the final structure. Conversely, let $(f_1\gamma_1\mu_{\ell_1}) * \dots * (f_n\gamma_n\mu_{\ell_n})$ be an element of $\mathcal{P}(X)$. Choose $\phi_1:[0,\frac{1}{2^{n-1}}] \iso^+ [0,\ell_1]$, $\phi_2:[0,\frac{1}{2^{n-1}}] \iso^+ [0,\ell_2]$, $\phi_3:[0,\frac{1}{2^{n-2}}] \iso^+ [0,\ell_3]$, etc... until $\phi_n:[0,\frac{1}{2^{}}] \iso^+ [0,\ell_n]$ and let $\phi=\phi_1 \ot \dots \ot \phi_n$. Using Proposition~\ref{startN2}, we obtain  
\begin{multline*}
\big(\big(f_1\gamma_1\mu_{\ell_1}\big) * \big(f_2\gamma_2\mu_{\ell_2}\big) * \big(f_3\gamma_3\mu_{\ell_3}\big) * \dots * \big(f_n\gamma_n\mu_{\ell_n}\big)\big)\phi = \\
\big(f_1\gamma_1\mu_{\ell_1}\phi_1\mu^{-1}_{\frac{1}{2^{n-1}}}\big) *_N \big(f_2\gamma_2\mu_{\ell_2}\phi_2\mu^{-1}_{\frac{1}{2^{n-1}}}\big) *_N \big(f_3\gamma_3\mu_{\ell_3}\phi_3\mu^{-1}_{\frac{1}{2^{n-2}}}\big) \\{*_N} \dots *_N \big(f_n\gamma_n\mu_{\ell_n}\phi_n\mu^{-1}_{\frac{1}{2^{}}}\big).
\end{multline*}
The continuous maps $\mu_{\ell_1}\phi_1\mu^{-1}_{\frac{1}{2^{n-1}}}, \mu_{\ell_2}\phi_2\mu^{-1}_{\frac{1}{2^{n-1}}},\mu_{\ell_3}\phi_3\mu^{-1}_{\frac{1}{2^{n-2}}}, \dots,\mu_{\ell_n}\phi_n\mu^{-1}_{\frac{1}{2^{}}}$ from $[0,1]$ to itself belong to $\mathcal{G}(1,1)$. Thus $\gamma'_1,\dots,\gamma'_n$ defined by the equalities 
\[
\begin{cases}
\gamma'_1 = \gamma_1\big(\mu_{\ell_1}\phi_1\mu^{-1}_{\frac{1}{2^{n-1}}}\big)\\
\gamma'_2 = \gamma_2\big(\mu_{\ell_2}\phi_1\mu^{-1}_{\frac{1}{2^{n-1}}}\big)\\
\gamma'_3 = \gamma_3\big(\mu_{\ell_3}\phi_3\mu^{-1}_{\frac{1}{2^{n-2}}}\big)\\
\dots\\
\gamma'_n = \gamma_n\big(\mu_{\ell_n}\phi_n\mu^{-1}_{\frac{1}{2^{}}}\big)
\end{cases}
\]
are execution paths of $X_1,\dots,X_n$ respectively. We obtain 
\begin{multline*}
\big(\big(f_1\gamma_1\mu_{\ell_1}\big) * \big(f_2\gamma_2\mu_{\ell_2}\big) * \big(f_3\gamma_3\mu_{\ell_3}\big) * \dots * \big(f_n\gamma_n\mu_{\ell_n}\big)\big)\phi = \\
\big(f_1\gamma'_1\big) *_N \big(f_2\gamma'_2\big) *_N \big(f_3\gamma'_3\big) *_N \dots *_N \big(f_n\gamma'_n\big).
\end{multline*}
We deduce that the set of paths $\mathcal{P}(X)$ is included in the $\Omega$-final structure. 
\epf

\section{From multipointed \mins{d}-spaces to Moore flows}
\label{adj}

\bdn \label{gspace}
	The enriched category of enriched presheaves from $\mathcal{G}$ to $\top$ is denoted by $\topdgr$. The underlying set-enriched category of enriched maps of enriched presheaves is denoted by $\topdgr_0$. The objects of $\topdgr_0$ are called the {\rm $\mathcal{G}$-spaces}. Let \[\mathbb{F}^{\mathcal{G}^{op}}_{\ell}U=\mathcal{G}(-,\ell)\p U \in \topdgr_0\] where $U$ is a topological space and where $\ell>0$.
\edn

\bp\label{ev-adj} \cite[Proposition~5.3 and Proposition~5.5]{dgrtop}
The category $\topdgr_0$ is a full reflective and coreflective subcategory of $\top^{\mathcal{G}^{op}_0}$. 
For every $\mathcal{G}$-space $F:\mathcal{G}^{op}\to \top$, every $\ell>0$ and every topological space $X$, we have the natural bijection of sets \[\topdgr_0(\mathbb{F}^{\mathcal{G}^{op}}_{\ell}X,F) \iso \top(X,F(\ell)).\] 
\ep

\bth (\cite[Theorem~5.14]{Moore1}) \label{closedsemimonoidal}
Let $D$ and $E$ be two $\mathcal{G}$-spaces. Let 
\[
D \ot E := \int^{(\ell_1,\ell_2)} \mathcal{G}(-,\ell_1+\ell_2) \p D(\ell_1) \p E(\ell_2).
\]
The pair $(\topdgr_0,\ot)$ has the structure of a closed symmetric semimonoidal category, i.e. a closed symmetric nonunital monoidal category.
\eth

\begin{nota}
	Let $D$ be a $\mathcal{G}$-space. Let $\phi:\ell\to \ell'$ be a map of $\mathcal{G}$. Let $x\in D(\ell')$. We will use the notation 
	\[
	x.\phi := D(\phi)(x).
	\]
	Intuitively, $x$ is a path of length $\ell'$ and $x.\phi$ is a path of length $\ell$ which is the reparametrization by $\phi$ of $x$.
\end{nota} 

Proposition~\ref{underlying-set-tensor} sheds light on the meaning of the tensor product of $\mathcal{G}$-spaces. It is used in the proof of Theorem~\ref{pre-calculation-pathspace}. It is not in \cite{Moore1}. The proof is given in this section and not in Section~\ref{unit} to recall \cite[Corollary~5.13]{Moore1} which also helps to understand the geometric contents of the tensor product of $\mathcal{G}$-spaces.

\bp \label{underlying-set-tensor}
Let $D_1,\dots,D_n$ be $n$ $\mathcal{G}$-spaces with $n\geq 1$. Then the mapping \[(x_1,\dots,x_n) \mapsto (\id,x_1,\dots,x_n)\] yields a surjective continuous map 
\[
\Phi_{D_1,\dots,D_n}:\displaystyle\bigsqcup\limits_{\substack{(\ell_1,\dots,\ell_n)\\\ell_1+\dots+\ell_n=L}}  D_1(\ell_1)\p \dots \p D_n(\ell_n)\longrightarrow (D_1\ot \dots \ot D_n)(L).
\]
\ep

\bpf
By \cite[Corollary~5.13]{Moore1}, the space $(D_1 \ot \dots \ot D_n)(L)$ is the quotient of the space
\[
\bigsqcup_{(\ell_1,\dots,\ell_n)} \mathcal{G}(L,\ell_1+\dots + \ell_n) \p D_1(\ell_1) \p \dots D_n(\ell_n).
\]
by the identifications 
\[
(\psi,x_1\phi_1,\dots,x_n\phi_n) = ((\phi_1\ot \dots \ot\phi_n)\psi,x_1,\dots,x_n)
\]
for all $\ell_1,\ell'_1,\dots,\ell_n,\ell'_n>0$, all $\psi\in \mathcal{G}(L,\ell_1+\dots+\ell_n)$, all $x_i\in D_i(\ell'_i)$ and all $\phi_i\in \mathcal{G}(\ell_i,\ell'_i)$. Let $\ell''_1,\dots,\ell''_n>0$ defined by induction on $i$ by the equation
\[
\forall 1\leq i \leq n, \ell''_i = \psi^{-1}\bigg(\sum_{1\leq j\leq i}\ell_j\bigg) - \sum_{1\leq j< i}\ell''_j.
\]
Note that $L=\ell''_1+\dots+\ell''_n$. We obtain 
\[
\forall 1\leq i \leq n, \psi\bigg(\sum_{1\leq j\leq i}\ell''_j\bigg) = \sum_{1\leq j\leq i}\ell_j.
\]
By Proposition~\ref{decomposition-tenseur}, there is a (unique) decomposition $\psi = \psi_1\ot \dots \ot \psi_n$ with $\psi_i\in \mathcal{G}(\ell''_i,\ell_i)$ for $1\leq i \leq n$. Then \[(\psi,x_1.\phi_1,\dots,x_n.\phi_n) = (\id,x_1.\phi_1.\psi_1,\dots,x_n.\phi_n.\psi_n)\] in $(D_1 \ot \dots \ot D_n)(L)$. Therefore the continuous map 
\[
\displaystyle\bigsqcup\limits_{\substack{(\ell_1,\dots,\ell_n)\\\ell_1+\dots+\ell_n=L}}  D_1(\ell_1)\p \dots \p D_n(\ell_n)\longrightarrow (D_1\ot \dots \ot D_n)(L).
\]
induced by the mapping $(x_1,\dots,x_n) \mapsto (\id,x_1,\dots,x_n)$ is surjective.
\epf

A \textit{semicategory}, also called \textit{nonunital category} in the literature, is a category without identity maps in the structure. It is \textit{enriched} over a closed symmetric semimonoidal category $(\mathcal{V},\ot)$ if it satisfied all axioms of enriched category except the one involving the identity maps, i.e. the enriched composition is associative and not necessarily unital. 

\bd \cite[Definition~6.2]{Moore1} \label{def-Moore-flow}
A {\rm Moore flow} is a small semicategory enriched over the closed semimonoidal category $(\topdgr_0,\ot)$ of Theorem~\ref{closedsemimonoidal}.
\ed

A Moore flow $X$ has therefore a set of objects denoted by $X^0$, and called \textit{states} in this context, and for each $(\alpha,\beta)\in X^0\p X^0$ a $\mathcal{G}$-space $\P_{\alpha,\beta} X$: the elements of \[\P_{\alpha,\beta}^{\ell} X=\P_{\alpha,\beta} X(\ell)\] for $\ell>0$ are called \textit{the execution paths of length $\ell$}. 

The category of Moore flows, denoted by $\dtopG$, is locally presentable by \cite[Theorem~6.11]{Moore1}. A map of Moore flows $f:X\to Y$ induces a set map $f^0:X^0\to Y^0$  and a map of $\mathcal{G}$-spaces $\P_{\alpha,\beta} f:\P_{\alpha,\beta} X\to \P_{f(\alpha),f(\beta)} Y$ for each $(\alpha,\beta)\in X^0\p X^0$. Let 
\[
\begin{aligned}
&\P X = \bigsqcup_{(\alpha,\beta) \in X^0\p X^0} \P_{\alpha,\beta} X\\
&\P Y = \bigsqcup_{(\alpha,\beta) \in Y^0\p Y^0} \P_{\alpha,\beta} Y\\
&\P f= \bigsqcup_{(\alpha,\beta) \in X^0\p X^0} \P_{\alpha,\beta} f.
\end{aligned}
\]

\begin{nota}
	The map $\P f:\P X\longrightarrow \P Y$ can be denoted by $f:\P X\to \P Y$ is there is no ambiguity. The set map $f^0:X^0\longrightarrow Y^0$ can be denoted by $f:X^0\longrightarrow Y^0$ is there is no ambiguity.
\end{nota}

 Every set $S$ can be viewed as a Moore flow with an empty $\mathcal{G}$-space of execution paths denoted in the same way. Let $D:\mathcal{G}^{op}\to \top$ be a $\mathcal{G}$-space. We denote by $\globP(D)$ the Moore flow defined as follows: 
	\[
	\begin{aligned}
	&\globP(D)^0 = \{0,1\}\\
	&\P_{0,0}\globP(D)=\P_{1,1}\globP(D)=\P_{1,0}\globP(D)=\Delta_{\mathcal{G}_0}\varnothing\\
	&\P_{0,1}\globP(D)=D.
	\end{aligned}	
	\]
There is no composition law. This construction yields a functor \[\globP:\topdgr_0\to \dtopG.\]
There exists a unique model structure on $\dtopG$ such that 
\[
\{\globP(\mathbb{F}^{\mathcal{G}^{op}}_{\ell}\mathbf{S}^{n-1}) \subset \globP(\mathbb{F}^{\mathcal{G}^{op}}_{\ell}\mathbf{D}^{n})\mid n\geq 0,\ell>0\}\cup \{C:\varnothing \to \{0\},R:\{0,1\} \to \{0\}\}
\]
is the set of generating cofibrations and such that all objects are fibrant. The set of generating trivial cofibrations is 
\[
\{\globP(\mathbb{F}^{\mathcal{G}^{op}}_{\ell}\mathbf{D}^{n}) \subset \globP(\mathbb{F}^{\mathcal{G}^{op}}_{\ell}\mathbf{D}^{n+1})\mid n\geq 0,\ell>0\}
\]
where the maps $\mathbf{D}^{n}\subset \mathbf{D}^{n+1}$ are induced by the mappings $(x_1,\dots,x_n) \mapsto (x_1,\dots,x_n,0)$. The weak equivalences are the map of Moore flows $f:X\to Y$ inducing a bijection $X^0\iso Y^0$ and such that for all $(\alpha,\beta)\in X^0\p X^0$, the map of $\mathcal{G}$-spaces $\P_{\alpha,\beta}X\to \P_{f(\alpha),f(\beta)}Y$ is an objectwise weak homotopy equivalence. The fibrations are the map of Moore flows $f:X\to Y$ such that for all $(\alpha,\beta)\in X^0\p X^0$, the map of $\mathcal{G}$-spaces $\P_{\alpha,\beta}X\to \P_{f(\alpha),f(\beta)}Y$ is an objectwise q-fibration of spaces. It is called \textit{the q-model structure} and we use the terminology of \textit{q-cofibration} and \textit{q-fibration} for naming the cofibrations and the fibrations respectively.

\bd \label{path-with-length}
Let $X$ be a multipointed $d$-space. Let $\P^\ell_{\alpha,\beta}X$ be the subspace of continuous maps from $[0,\ell]$ to $|X|$ defined by \[\P^\ell_{\alpha,\beta}X = \{t\mapsto \gamma\mu_\ell\mid \gamma\in \P^{\mathcal{G}}_{\alpha,\beta}X\}.\] Its elements are called {\rm the execution paths of length $\ell$} from $\alpha$ to $\beta$. Let \[\P^\ell X = \bigsqcup_{(\alpha,\beta) \in X^0\p X^0} \P^\ell_{\alpha,\beta}X.\]
A map of multipointed $d$-spaces $f:X\to Y$ induces for each $\ell>0$ a continuous map $\P^\ell f:\P^\ell X \to \P^\ell Y$ by composition by $f$ (in fact by $|f|$).
\ed

Note that $\P^1_{\alpha,\beta}X = \P^{\mathcal{G}}_{\alpha,\beta}X$, that there is a homeomorphism $\P^{\ell}_{\alpha,\beta}X \iso \P^{\mathcal{G}}_{\alpha,\beta}X$ for all $\ell>0$, and that for any topological space $Z$, we have the homeomorphism \[\P^\ell_{{0},{1}}(\globG(Z)) \iso \mathcal{G}(\ell,1) \p Z \] for any $\ell > 0$ by Proposition~\ref{calcul-topology-glob}. 

The definition above of an execution path of length $\ell>0$ is not restrictive. Indeed, we have: 

\bp Let $X$ be a multipointed $d$-space. Let $\phi:[0,\ell]\iso^+ [0,\ell]$. Let $\gamma\in \P^\ell X$. Then $\gamma\phi\in \P^\ell X$. \ep

\bpf By definition of $\P^\ell X$, there exists $\overline{\gamma}\in \P^{\mathcal{G}} X$ such that $\gamma = \overline{\gamma}\mu_{\ell}$. We obtain $\gamma\phi = \overline{\gamma}\mu_{\ell}\phi \mu^{-1}_{\ell}\mu_{\ell}$. Since $\mu_{\ell}\phi \mu^{-1}_{\ell}\in \mathcal{G}(1,1)$, we deduce that $\overline{\gamma}\mu_{\ell}\phi \mu^{-1}_{\ell}\in \P^{\mathcal{G}} X$ and that $\gamma\phi\in \P^\ell X$. 
\epf

\bp \label{variable-length}
Let $X$ be a multipointed $d$-space. Let $\gamma_1$ and $\gamma_2$ be two execution paths of $X$ with $\gamma_1(1) = \gamma_2(0)$. Let $\ell_1,\ell_2>0$. Then \[\big(\gamma_1\mu_{\ell_1} *  \gamma_2\mu_{\ell_1}\big)\mu^{-1}_{\ell_1+\ell_2}\] is an execution path of $X$. 
\ep

\bpf 
Let $\phi_1:[0,\frac{1}{2}]\iso^+ [0,\ell_1]$ and $\phi_2:[0,\frac{1}{2}]\iso^+ [0,\ell_2]$. Then we have \[\phi_1\ot \phi_2:[0,1]\iso^+ [0,\ell_1+\ell_2].\] We obtain the sequence of equalities
\[\begin{aligned}
\big((\gamma_1\mu_{\ell_1}) *  (\gamma_2\mu_{\ell_2})\big)\mu^{-1}_{\ell_1+\ell_2} &= \big((\gamma_1\mu_{\ell_1}) *  (\gamma_2\mu_{\ell_2})\big)\big(\phi_1\ot \phi_2\big)\big(\phi_1\ot \phi_2\big)^{-1}\mu^{-1}_{\ell_1+\ell_2}\\
&= \big((\gamma_1\mu_{\ell_1}\phi_1) *  (\gamma_2\mu_{\ell_1}\phi_2)\big)\big(\phi_1\ot \phi_2\big)^{-1}\mu^{-1}_{\ell_1+\ell_2}\\
&= \big((\gamma_1\mu_{\ell_1}\phi_1\mu^{-1}_{\frac{1}{2}}\mu_{\frac{1}{2}}) *  (\gamma_2\mu_{\ell_2}\phi_2\mu^{-1}_{\frac{1}{2}}\mu_{\frac{1}{2}})\big)\big(\phi_1\ot \phi_2\big)^{-1}\mu^{-1}_{\ell_1+\ell_2}\\
&= \big((\gamma_1\underbrace{\mu_{\ell_1}\phi_1\mu^{-1}_{\frac{1}{2}}}_{\in \mathcal{G}(1,1)}) *_N (\gamma_2\underbrace{\mu_{\ell_2}\phi_2\mu^{-1}_{\frac{1}{2}}}_{\in \mathcal{G}(1,1)}) \big)\underbrace{\big(\phi_1\ot \phi_2\big)^{-1}\mu^{-1}_{\ell_1+\ell_2}}_{\in \mathcal{G}(1,1)},
\end{aligned}\]
the first equality because $\phi_1\ot \phi_2$ is invertible, the second equality by Proposition~\ref{lem-1}, the third equality because $\mu_{\frac{1}{2}}$ is invertible, and finally the last equality by Proposition~\ref{startN/2}. The proof is complete because the set of execution paths of $X$ is invariant by the action of $\mathcal{G}(1,1)$. 
\epf

\bp \label{addlength} Let $X$ be a multipointed $d$-space. Let $\ell_1,\ell_2>0$. 
The Moore composition of continuous maps yields a continuous maps $\P^{\ell_1}X \p \P^{\ell_2}X \to \P^{\ell_1+\ell_2}X$. 
\ep

\bpf It is a consequence of Definition~\ref{path-with-length} and Proposition~\ref{variable-length}
\epf

\bth \label{adj-multi-moore} Let $X$ be a multipointed $d$-space. Then the following data
\begin{itemize}[leftmargin=*]
	\item The set of states $X^0$ of $X$
	\item For all $\alpha,\beta\in X^0$ and all real numbers $\ell>0$, let
	\[
	\P_{\alpha,\beta}^{\ell}\moore^{\mathcal{G}}(X) := \P_{\alpha,\beta}^{\ell}X.
	\]
	\item For all maps $[0,\ell]\iso^+[0,\ell']$, a map $f:[0,\ell']\to |X|$ of $\P_{\alpha,\beta}^{\ell'}\moore^{\mathcal{G}}(X)$ is mapped to the map $[0,\ell]\iso^+[0,\ell']\stackrel{f}\to |X|$ of $\P_{\alpha,\beta}^{\ell}\moore^{\mathcal{G}}(X)$ 
	\item For all $\alpha,\beta,\gamma\in X^0$ and all real numbers $\ell,\ell'>0$, the composition maps \[*:\P_{\alpha,\beta}^{\ell}\moore^{\mathcal{G}}(X) \p \P_{\beta,\gamma}^{\ell'}\moore^{\mathcal{G}}(X) \to \P_{\alpha,\gamma}^{\ell+\ell'}\moore^{\mathcal{G}}(X)\] of Proposition~\ref{addlength}.
\end{itemize}
assemble to a Moore flow $\moore^{\mathcal{G}}(X)$. This mapping induces a functor \[\moore^{\mathcal{G}}:\ptop{\mathcal{G}}\longrightarrow\dtopG\] which is a right adjoint.
\eth

Note that the left adjoint $\moore^{\mathcal{G}}_!:\dtopG \longrightarrow \ptop{\mathcal{G}}$ preserves the set of states as well as the functor $\moore^{\mathcal{G}}:\ptop{\mathcal{G}}\longrightarrow\dtopG$. 

\bpf These data give rise to a $\mathcal{G}$-space $\P_{\alpha,\beta}\moore^{\mathcal{G}}(X)$ for each pair $(\alpha,\beta)$ of states of $X^0$ and, thanks to Proposition~\ref{addlength}, to an associative composition law $*:\P^{\ell_1}_{\alpha,\beta}\moore^{\mathcal{G}}(X) \p \P^{\ell_2}_{\beta,\gamma}\moore^{\mathcal{G}}(X) \to \P^{\ell_1+\ell_2}_{\alpha,\gamma}\moore^{\mathcal{G}}(X)$ which is natural with respect to $(\ell_1,\ell_2)$. By \cite[Section~6]{Moore1}, these data assemble to a Moore flow. Since limits and colimits of $\mathcal{G}$-spaces are calculated objectwise, the functor $\moore^{\mathcal{G}}:\ptop{\mathcal{G}}\longrightarrow\dtopG$ is limit-preserving and accessible by Theorem~\ref{pre-right-adj}. Therefore it is a right adjoint by \cite[Theorem~1.66]{TheBook}. 
\epf

\bp \label{precalcul} Let $X$ be a multipointed $d$-space. Let $\ell>0$ be a real number. Let $Z$ be a topological space. Then there is a bijection of sets \[\ptop{\mathcal{G}}(\glob^{\mathcal{G}}_\ell(Z),X) \iso \bigsqcup_{(\alpha,\beta)\in X^0\p X^0}\top(Z,\P_{\alpha,\beta}^{\ell}X)\]
which is natural with respect to $Z$ and $X$.
\ep

\bpf
A map $f$ of multipointed $d$-spaces from $\glob^{\mathcal{G}}_\ell(Z)$ to $X$ is determined by 
\begin{itemize}[leftmargin=*]
	\item The image by $f$ of ${0}$ and ${1}$ which will be denoted by $\alpha$ and $\beta$ respectively
	\item A continuous map (still denoted by $f$) from $|\glob^{\mathcal{G}}_\ell(Z)|$ to $|X|$ such that for all $x\in Z$ and all $\phi:[0,1]\iso^+[0,\ell]$, the map $t\mapsto f(x,\phi(t))$ from $[0,1]$ to $|X|$ belongs to $\P^{\mathcal{G}}_{\alpha,\beta}X$.
\end{itemize}
By definition of $\P_{\alpha,\beta}^{\ell}X$, for every $x\in Z$, the continuous map $f(x,-)$ from $[0,\ell]$ to $|X|$ belongs to  $\P_{\alpha,\beta}^{\ell}X$ since $f(x,-) = f(x,\phi(-)).\phi^{-1}$ for any $\phi:[0,1]\iso^+[0,\ell]$. Since $f$ is continuous and since $\top$ is cartesian closed, the mapping $x\mapsto f(x,-)$ actually yields a continuous map from $Z$ to $\P_{\alpha,\beta}^{\ell}X$. Conversely, starting from a continuous map $g:Z\to \P_{\alpha,\beta}^{\ell}X$, one can define a map of multipointed $d$-spaces from $\glob^{\mathcal{G}}_\ell(Z)$ to $X$ by taking ${0}$ and ${1}$ to $\alpha$ and $\beta$ respectively and by taking $(x,t)\in |\glob^{\mathcal{G}}_\ell(Z)|$ to $g(x)(t)$. 
\epf

We want to recall for the convenience of the reader:

\bp \label{map-from-glob} \cite[Proposition~6.10]{Moore1} Let $D:\mathcal{G}^{op}\to \top$ be a $\mathcal{G}$-space. Let $X$ be a Moore flow. Then there is the natural bijection \[\dtopG(\globP(D),X)\iso \bigsqcup_{(\alpha,\beta)\in X^0\p X^0} \topdgr_0(D,\P_{\alpha,\beta}X).\]
\ep

\bp \label{calculM}
For all topological spaces $Z$ and all $\ell>0$, there are the natural isomorphisms 
\[\begin{aligned}
&\moore^{\mathcal{G}}(\globG_\ell(Z)) \iso  \globP(\mathbb{F}^{\mathcal{G}^{op}}_\ell(Z)), \\
&\mathbb{M}_!^{\mathcal{G}}(\globP(\mathbb{F}^{\mathcal{G}^{op}}_\ell(Z))) \iso  \globG_\ell(Z).
\end{aligned}\]
\ep

\bpf By definition of $\moore^{\mathcal{G}}$ and by Proposition~\ref{calcul-topology-glob}, the only nonempty path $\mathcal{G}$-space of $\moore^{\mathcal{G}}(\globG_\ell(Z))$ is \[\P_{0,1}\moore^{\mathcal{G}}(\globG_\ell(Z)) = \mathcal{G}(-,\ell)\p Z\] and we obtain the first isomorphism. There is the sequence of natural bijections, for any multipointed $d$-space $X$,  
\[\begin{aligned}
\ptop{\mathcal{G}}\big(\mathbb{M}_!^{\mathcal{G}}(\globP(\mathbb{F}^{\mathcal{G}^{op}}_\ell(Z))),X \big)&\iso  \dtopG\big(\globP(\mathbb{F}^{\mathcal{G}^{op}}_\ell(Z)),\moore^{\mathcal{G}}X\big)\\
&\iso \bigsqcup_{(\alpha,\beta)\in X^0\p X^0} \topdgr_0\big(\mathbb{F}^{\mathcal{G}^{op}}_\ell(Z),\P_{\alpha,\beta}X\big)\\
&\iso \bigsqcup_{(\alpha,\beta)\in X^0\p X^0} \top(Z,\P_{\alpha,\beta}^{\ell}X)\\
&\iso \ptop{\mathcal{G}}(\glob^{\mathcal{G}}_\ell(Z),X),
\end{aligned}\]
the first bijection by adjunction, the second bijection by Proposition~\ref{map-from-glob}, the third bijection by Proposition~\ref{ev-adj} and the last bijection by Proposition~\ref{precalcul}. The proof of the second isomorphism is then complete thanks to the Yoneda lemma.
\epf

\section{Cellular multipointed \mins{d}-spaces}
\label{cellular-obj}

Let $\lambda$ be an ordinal. In this section, we work with a colimit-preserving functor \[X:\lambda \longrightarrow \ptop{\mathcal{G}}\] such that
\begin{itemize}[leftmargin=*]
	\item The multipointed $d$-space $X_0$ is a set, in other terms $X_0=(X^0,X^0,\varnothing)$ for some set $X^0$.
	\item For all $\nu<\lambda$, there is a pushout diagram of multipointed $d$-spaces 
	\[
	\xymatrix@C=3em@R=3em
	{
		\globG(\mathbf{S}^{n_\nu-1}) \fd{} \fr{g_\nu} & X_\nu \ar@{->}[d]^-{} \\
		\globG(\mathbf{D}^{n_\nu}) \fr{\widehat{g_\nu}} & \cocartesien X_{\nu+1}
	}
	\]
	with $n_\nu \geq 0$. 
\end{itemize}

Let $X_\lambda = \liminj_{\nu<\lambda} X_\nu$. Note that for all $\nu\leq \lambda$, there is the equality $X_\nu^0=X^0$. Denote by \[c_\nu = |\globG(\mathbf{D}^{n_\nu})|\backslash |\globG(\mathbf{S}^{n_\nu-1})|\] the $\nu$-th cell of $X_\lambda$. It is called a \textit{globular cell}. Like in the usual setting of CW-complexes, $\widehat{g_\nu}$ induces a homeomorphism from $c_\nu$ to $\widehat{g_\nu}(c_\nu)$ equipped with the relative topology which will be therefore denoted in the same way. It also means that $\widehat{g_\nu}(c_\nu)$ equipped with the relative topology is $\Delta$-generated. The closure of $c_\nu$ in $|X_\lambda|$ is denoted by \[\widehat{c_{\nu}} = \widehat{g_{\nu}}(|\globG(\mathbf{D}^{n_\nu})|).\] The boundary of $c_\nu$ in $|X_\lambda|$ is denoted by \[\de c_\nu = \widehat{g_{\nu}}(|\globG(\mathbf{S}^{n_\nu-1})|).\]
The state $\widehat{g_\nu}(0)\in X^0$ ($\widehat{g_\nu}(1)\in X^0$ resp.)  is called the \textit{initial (final resp.) state} of $c_\nu$. The integer $n_\nu+1$ is called the \textit{dimension} of the globular cell $c_\nu$. It is denoted by $\dim c_\nu$. The states of $X^0$ are also called the \textit{globular cells of dimension $0$}.

\bd \label{finitecellular}
The cellular multipointed $d$-space $X_\lambda$ is {\rm finite} if $\lambda$ is a finite ordinal and $X^0$ is finite.
\ed

\bp  \label{p1}
The space $|X_\lambda|$ is a cellular space. It contains $X^0$ as a discrete closed subspace. The space $|X_\lambda|$ is weakly Hausdorff. For every $0\leq \nu_1 \leq \nu_2 \leq \lambda$, the continuous map $|X_{\nu_1}| \to |X_{\nu_2}|$ is a q-cofibration of spaces, and in particular a closed $T_1$-inclusion.
\ep

\bpf
By \cite[Theorem~8.2]{4eme}, the continuous map \[|\globG(\mathbf{S}^{n_\nu-1})|\to |\globG(\mathbf{D}^{n_\nu})|\] is a q-cofibration of spaces for all $\nu \geq 0$ between cellular spaces. Since the functor $X\mapsto |X|$ is colimit-preserving, the space $|X_\lambda|$ is a cellular space. It is therefore weakly Hausdorff. For every $0\leq \nu_1 \leq \nu_2 \leq \lambda$, the continuous map $|X_{\nu_1}| \to |X_{\nu_2}|$ is a transfinite composition of q-cofibrations, and hence a q-cofibration. The map $X^0\to X_\lambda$ is a transfinite composition of q-cofibrations, and therefore a q-cofibration, and in particular a closed $T_1$-inclusion. Every subset of $X^0$ is closed since $X^0$ is equipped with the discrete topology. Consequently, $X^0$ is a discrete closed subspace of $|X_\lambda|$. 
\epf

\bp \label{restriction_path}
For all $0\leq \nu_1\leq \nu_2\leq\lambda$, there is the equality \[\P^{\mathcal{G}}X_{\nu_1} = \P^{\mathcal{G}} X_{\nu_2} \cap \ttop([0,1],|X_{\nu_1}|).\]
\ep

\bpf It is trivial for $\nu_1=\nu_2$. For $\nu_2=\nu_1+1$, there is a pushout diagram of multipointed $d$-spaces 
\[
\xymatrix@C=3em@R=3em
{
	\globG(\mathbf{S}^{n_{\nu_1}-1}) \fd{} \fr{g_{\nu_1}} & X_{\nu_1} \ar@{->}[d]^-{} \\
	\globG(\mathbf{D}^{n_{\nu_1}}) \fr{\widehat{g_{\nu_1}}} & \cocartesien X_{\nu_2}.
}
\]
The equality holds because the set of execution paths of $X_{\nu_2}$ is obtained as a $\Omega$-final structure. We conclude by a transfinite induction on $\nu_2$. 
\epf

\bp \label{pushout-qcof-plus}
For all $0\leq \nu_1\leq \nu_2\leq\lambda$, the continuous map $\P^{\mathcal{G}}X_{\nu_1}\to \P^{\mathcal{G}}X_{\nu_2}$ is a $\Delta$-inclusion.
\ep

\bpf
Consider a set map $[0,1]\to \P^{\mathcal{G}}X_{\nu_1}$ such that the composite set map 
\[
[0,1]\longrightarrow \P^{\mathcal{G}}X_{\nu_1} \longrightarrow \P^{\mathcal{G}}X_{\nu_2}
\]
is continuous. Then by adjunction, we obtain a continuous map 
\[
[0,1] \p [0,1] \longrightarrow |X_{\nu_2}|.
\]
By hypothesis, it factors as a composite of set maps 
\[
[0,1] \p [0,1] \longrightarrow |X_{\nu_1}| \longrightarrow |X_{\nu_2}|.
\]
By Proposition~\ref{p1}, the left-hand map is continuous since $[0,1] \p [0,1]$ is compact. The proof is complete thanks to Proposition~\ref{restriction_path} and Proposition~\ref{DeltaIncl}.
\epf

\bp \label{cpt-intersect-finite} 
Let $K$ be a compact subspace of $|X_\lambda|$. Then $K$ intersects finitely many $c_\nu$. 
\ep

\bpf We mimick the proof of \cite[Proposition~A.1]{MR1867354} for the transfinite case. Assume that there exists an infinite set $S=\{m_j\mid j\geq 0\}$ with $m_j\in K\cap c_{\nu_j}$, where $(\nu_j)_{j\geq 0}$ is a sequence of mutually distinct ordinals. By transfinite induction on $\nu\geq 0$, let us prove that $S\cap |X_\nu|$ is a closed subset of $|X_\nu|$. The assertion is trivial for $\nu=0$. There is the pushout diagram of spaces for all $\nu<\lambda$
\[
\xymatrix@C=3em@R=3em
{
	|\globG(\mathbf{S}^{n_{\nu}-1})| \fd{} \fr{g_{\nu}} & |X_{\nu}| \ar@{->}[d]^-{} \\
	|\globG(\mathbf{D}^{n_{\nu}})| \fr{\widehat{g_{\nu}}} & \cocartesien {|X_{\nu+1}|}.
}
\]
By induction hypothesis, $g_{\nu}^{-1}(S\cap |X_\nu|)$ is a closed subset of $|\globG(\mathbf{S}^{n_{\nu}-1})|$ and $\widehat{g_{\nu}}^{-1}(S\cap |X_{\nu+1}|)$ is equal to $g_{\nu}^{-1}(S\cap |X_\nu|)$ union at most one point. Therefore, $S\cap |X_{\nu+1}|$ is a closed subset of $|X_{\nu+1}|$ because the latter space is equipped with the final topology by Proposition~\ref{p1}. Suppose that we have proved that for all $\nu<\nu'$, $S\cap |X_\nu|$ is a closed subset of $|X_\nu|$ where $\nu'$ is a limit ordinal. Then, since the topology of $|X_{\nu'}|$ is the final topology (it is a tower of one-to-one maps), $S\cap |X_{\nu'}|$ is a closed subset of $|X_{\nu'}|$. Thus, by transfinite induction on $\nu\geq 0$, we prove that $S$ is closed in $|X_\nu|$ for all $0\leq \nu\leq \lambda$. The same argument proves that every subset of $S$ is closed in $|X_\lambda|$. Thus $S$ has the discrete topology. But it is compact, being a closed subset of the compact space $K$, and therefore finite. Contradiction.
\epf

Colimits of multipointed $d$-spaces are calculated by taking the colimit of the underlying spaces and of the sets of states and by taking the $\Omega$-final structure which is generated by the free finite compositions of execution paths. Consequently, the composite functor 
\[
\xymatrix@C=4em
{
	\ptop{\mathcal{G}}\fr{\P^{\mathcal{G}}} &  \top \fr{\subset} & \set
}
\]
is finitely accessible. It is unlikely that the functor $\P^{\mathcal{G}}:\ptop{\mathcal{G}}\to \top$, which is a right adjoint by Theorem~\ref{pre-right-adj}, is finitely accessible. However, we have: 

\bth  \label{topological-path-almost-accessible}
The composite functor 
\[
\lambda \stackrel{X} \longrightarrow \ptop{\mathcal{G}} \stackrel{\P^{\mathcal{G}}}\longrightarrow \top
\]
is colimit-preserving. In particular the continuous bijection \[\liminj (\P^{\mathcal{G}}.X) \longrightarrow \P^{\mathcal{G}} \liminj X\] is a homeomorphism. Moreover the topology of $\P^{\mathcal{G}} \liminj X$ is the final topology.
\eth

Note that Theorem~\ref{topological-path-almost-accessible} holds both for $\Delta$-generated spaces and $\Delta$-Hausdorff $\Delta$-generated spaces. 

\bpf
Consider the set of ordinals
\[
\bigg\{
\nu\leq \lambda \mid \nu \hbox{ limit ordinal and } \liminj_{\nu'<\nu}(\mathbb{P}^{\mathcal{G}}X_{\nu'}) \longrightarrow \mathbb{P}^{\mathcal{G}}X_\nu
\hbox{ not isomorphism}\bigg\}
\]
Assume this set nonempty. Let $\nu$ be its smallest element. The topology of $\liminj_{{\nu'}<\nu}\P^{\mathcal{G}}X_{\nu'}$ is the final topology because the continuous maps $\P^{\mathcal{G}}X_{\nu'}\to \P^{\mathcal{G}}X_{{\nu'}+1}$ are one-to-one. Let $f:[0,1]\to \P^{\mathcal{G}} X_\nu$ be a continuous map.  Therefore the composite map 
\[
[0,1]\stackrel{f}\longrightarrow \P^{\mathcal{G}} X_{\nu} \subset \ttop([0,1],|X_\nu|)
\]
is continuous. It gives rise by adjunction to a continuous map $[0,1] \p [0,1] \to |X_{\nu}|$. Since the functor $X:\lambda \to \ptop{\mathcal{G}}$ is colimit-preserving, there is the homeomorphism $|X_\nu| \iso \liminj_{{\nu'}<\nu} |X_{\nu'}|$. Since $[0,1] \p [0,1]$ is compact, the latter continuous map then factors as a composite $[0,1] \p [0,1] \to |X_{\nu'}| \to |X_{\nu}|$ for some ordinal ${\nu'}<\nu$ by Proposition~\ref{p1}. Since $\P^{\mathcal{G}}X_{\nu'} = \P^{\mathcal{G}} X_\nu \cap \ttop([0,1],|X_{\nu'}|)$ by Proposition~\ref{restriction_path}, $f$ factors as a composite $[0,1]\to \P^{\mathcal{G}} X_{\nu'}\to \P^{\mathcal{G}} X_\nu$. Using Corollary~\ref{DeltaHomeo2}. we obtain the homeomorphism $\liminj_{{\nu'}<\nu}\mathbb{P}^{\mathcal{G}}X_{\nu'} \longrightarrow \mathbb{P}^{\mathcal{G}}X_\nu$: contradiction.
\epf

\bth \label{cof-accessible}
The composite functor 
\[
\lambda \stackrel{X} \longrightarrow \ptop{\mathcal{G}} \stackrel{\moore^{\mathcal{G}}}\longrightarrow \dtopG
\]
is colimit-preserving. In particular the natural map \[\liminj_{\nu<\lambda} \moore^{\mathcal{G}}(X_\nu) \longrightarrow \moore^{\mathcal{G}} X_\lambda\] is an isomorphism.
\eth

\bpf
Theorem~\ref{topological-path-almost-accessible} states that there is the homeomorphism 
\[\liminj_{\nu<\lambda} \P^{\mathcal{G}}X_\nu \longrightarrow \P^{\mathcal{G}}X_\lambda.\]
We have, by definition of the functor $\moore^{\mathcal{G}}$, the equality of functors $\P^{\mathcal{G}}=\P^1.\moore^{\mathcal{G}}$. It means that there is the homeomorphism 
\[\liminj_{\nu<\lambda} \P^1\moore^{\mathcal{G}}(X_\nu) \longrightarrow \P^1\moore^{\mathcal{G}}(X_\lambda).\]
Since all maps the reparametrization category $\mathcal{G}$ are isomorphisms, we obtain for all $\ell>0$ the homeomorphism
\[\liminj_{\nu<\lambda} \P^\ell\moore^{\mathcal{G}}(X_\nu) \longrightarrow \P^\ell\moore^{\mathcal{G}}(X_\lambda).\]
Since colimits of $\mathcal{G}$-spaces are calculated objectwise, we obtain the isomorphism of $\mathcal{G}$-spaces 
\[\liminj_{\nu<\lambda} \P\moore^{\mathcal{G}}X_\nu \longrightarrow \P\moore^{\mathcal{G}}X_\lambda.\]
The proof is complete thanks to the universal property of the colimits.
\epf

\bd An execution path $\gamma$ of a multipointed $d$-space $X$ is {\rm  minimal} if \[\gamma(]0,1[) \cap X^0 = \varnothing.\]
\ed

For any (q-cofibrant or not) topological space $Z$, every execution path of the multipointed $d$-space $\globG(Z)$ is minimal. The following theorem proves that execution paths of cellular multipointed $d$-spaces have a normal form.

\bth \label{normal-form}
Let $\gamma$ be an execution path of $X_\lambda$. Then there exist minimal execution paths $\gamma_1,\dots,\gamma_n$ and $\ell_1,\dots,\ell_n>0$ with $\sum_i \ell_i=1$ such that 
\[
\gamma = (\gamma_1\mu_{\ell_1}) * \dots * (\gamma_n\mu_{\ell_n}).
\]
Moreover, if there is the equality 
\[
\gamma = (\gamma_1\mu_{\ell_1}) * \dots * (\gamma_n\mu_{\ell_n}) = (\gamma'_1\mu_{\ell'_1}) * \dots * (\gamma_{n'}\mu_{\ell'_{n'}})
\]
such that all $\gamma'_j$ are also minimal and with $\ell'_1,\dots,\ell'_{n'}>0$, then $n=n'$ and $\gamma_i=\gamma'_i$ and $\ell_i=\ell'_i$ for all $1\leq i \leq n$.
\eth

\bpf
The set of execution paths of $X_\lambda$ is obtained as a $\Omega$-final structure. Using Theorem~\ref{final-structure-revisited}, we obtain \[
\gamma = (\gamma_1\mu_{\ell_1}) * \dots * (\gamma_n\mu_{\ell_n}).
\]
for some $n\geq 1$ with $\ell_1+\dots+\ell_n=1$ such that for all $1\leq i\leq n$, there exists a globular cell $c_{\nu_i}$ such that 
\[
\begin{aligned}
&\gamma_i(]0,1[) \subset c_{\nu_i},\\
&\gamma_i(0)=\widehat{g_{\nu_i}}(0),\\
&\gamma_i(0)=\widehat{g_{\nu_i}}(1).
\end{aligned}
\]
Therefore there exists a finite set $\{t_0,\dots,t_n\}$ with $t_0=0<t_1<\dots <t_n=1$ and $n\geq 1$ such that $\gamma([0,1])\cap X^0 = \{\gamma(t_i)\mid 0\leq i\leq n\}$. We necessarily have $\ell_i = t_i-t_{i-1}$ for $1\leq i\leq n$. Let $\ell_0=0$. Then we deduce that $\sum_{j<i} \ell_j = t_{i-1}$ and $\sum_{j\leq i} \ell_j = t_{i}$. The equality $\gamma = (\gamma_1\mu_{\ell_1}) * \dots * (\gamma_n\mu_{\ell_n})$ therefore implies that $\gamma(t) = (\gamma_i\mu_{\ell_i})(t-t_{i-1})$ for $t_{i-1}\leq t \leq t_i$ for all $1\leq i\leq n$ by definition of the Moore composition of paths. We deduce that we necessarily have the equalities $\gamma_i(t) = \gamma(\ell_it+t_{i-1})$ for $t\in [0,1]$. 
\epf

Let $\gamma$ be an execution path of $X_\lambda$. Consider the normal form 
\[
\gamma = (\gamma_1\mu_{\ell_1}) * \dots * (\gamma_n\mu_{\ell_n}).
\]
of Theorem~\ref{normal-form}. There exists a unique sequence $[c_{\nu_1},\dots,c_{\nu_n}]$ of globular cells such that for all $1\leq i \leq n$, $\gamma_i(]0,1[)\subset c_{\nu_i}$, $\gamma_i(0)=\widehat{g_{\nu_i}}(0)$ and $\gamma_i(1)=\widehat{g_{\nu_i}}(1)$. This leads to the following notion: 

\bd \label{def-carrier} With the notations above. The sequence of globular cells \[\carrier(\gamma)=[c_{\nu_1},\dots,c_{\nu_n}]\] is called the {\rm carrier} of $\gamma$. The integer $n$ is called the {\rm length} of the carrier. 
\ed

\bp \label{carac-dec0} An execution path of $X_\lambda$ is minimal if and only if the length of its carrier is $1$.
\ep

\bpf It is a consequence of Theorem~\ref{normal-form}.
\epf

\bp \label{carac-dec}
An execution path $\gamma$ of $X_\lambda$ is non-minimal if and only if there exist two execution paths $\gamma_1$ and $\gamma_2$ such that $\gamma = \gamma_1 *_N \gamma_2$.
\ep

Proposition~\ref{carac-dec} does not hold for non q-cofibrant multipointed $d$-spaces. Consider e.g. the multipointed $d$-space $X$ obtained by starting from the directed segment $\vI^{\mathcal{G}}$ and by adding to the set of states $\{0,1\}$ the point $\frac{1}{2}$. Then all execution paths of $X$ are non-minimal and $\P_{0,\frac{1}{2}}X=\P_{\frac{1}{2},1}X=\varnothing$. Note that the q-cofibrant replacement of $X$ consists of the disjoint sum $\vI^{\mathcal{G}} \sqcup \{\frac{1}{2}\}$.

\bpf It is a consequence of Theorem~\ref{normal-form} and Proposition~\ref{carac-dec0}. 
\epf

\bp \label{all-regular}
All execution path of $X_\lambda$ are locally injective. 
\ep

In \cite{MR2369163}, the terminology of \textit{regular} paths is used.

\bpf All execution paths of globes $\globG(Z)$ are one-to-one for all topological spaces $Z$. Therefore all minimal execution paths are locally injective (it can be a loop). The proof is complete thanks to Theorem~\ref{normal-form}.
\epf

\bp \label{noescape} 
Consider a minimal execution path $\gamma$ of $X_\lambda$ with $\carrier(\gamma)=[c_{\nu_0}]$. Let $c_\nu$ be a globular cell of $X_\lambda$ with $\nu\neq \nu_0$. Then the following two assertions are equivalent: 
\begin{enumerate}
	\item $\gamma(]0,1[) \cap \widehat{c_{\nu}}\neq \varnothing$
	\item $\gamma([0,1]) \subset \de c_{\nu}$.
\end{enumerate}
Moreover, when the previous assertions are satisfied, there exists an execution path $\gamma'$ from the initial state of $c_{\nu}$ to its final state such that $\gamma'= (\gamma_1\mu_{\ell_1}) * \dots * (\gamma_n\mu_{\ell_n})$ with $\gamma=\gamma_i$ for at least one $i\in \{1,\dots,n\}$, $\gamma_1,\dots,\gamma_n$ minimal and $\sum_i\ell_i = 1$. 
\ep

\bpf 
Since $\de c_{\nu} \subset \widehat{c_{\nu}}$, we deduce $(2)\Rightarrow (1)$. Assume $(1)$. Since $\gamma(]0,1[)\subset c_{\nu_0}$ and $\nu\neq \nu_0$, one has $\nu>\nu_0$. It means that there exists a point $\widehat{g_\nu}(z,t)$ of $\de c_{\nu}$ which belongs to $c_{\nu_0}$ with $z\in \mathbf{S}^{n_\nu-1}$ and, since $c_{\nu_0}\cap X^0=\varnothing$, with $t\in ]0,1[$. Therefore the carrier of the execution path $\widehat{g_\nu}\delta_z$ contains the globular cell $c_{\nu_0}$. We deduce that there exists $\phi\in \mathcal{G}(1,1)$ such that 
\[
\widehat{g_\nu}\delta_z\phi = (\gamma_1\mu_{\ell_1}) * \dots * (\gamma_n\mu_{\ell_n})
\]
with $\gamma=\gamma_i$ for at least one $i\in \{1,\dots,n\}$, $\gamma_1,\dots,\gamma_n$ minimal and $\sum_i\ell_i = 1$. In particular, we deduce that $\gamma([0,1]) \subset \de c_{\nu}$: we have proved $(1)\Rightarrow(2)$.
\epf

\bd \label{achronal}
	Let $c_\nu$ be a globular cell of $X_\lambda$. Let $0<h<1$. Let
	\[
	\widehat{c_{\nu}}[h] = \bigg\{\widehat{g_{\nu}}(z,h)\mid (z,h)\in |\globG(\mathbf{D}^{n_{\nu}})|\bigg\}
	\]
	It is called an {\rm achronal slice} of the globular cell $c_\nu$.
\ed

\bp \label{middle0}
For any globular cell $c_\nu$ of $X_\lambda$ and any minimal execution path $\gamma$ and any $h\in ]0,1[$, the cardinal of the set \[\bigg\{t\in ]0,1[\mid \gamma(t) \in \widehat{c_{\nu}}[h]\bigg\}\] is at most one. 
\ep

In other terms, a minimal execution path of $X_\lambda$ intersects any achronal slice at most one time. Remember that execution paths of $X_\lambda$ are locally injective, i.e. they do not contain zero speed points. Proposition~\ref{middle0} does not hold in general for a non-minimal execution path because it could go back to the initial state of the globular cell after reaching its final state, which moreover could be equal to the initial state of the globular cell.

\bpf If the set $\gamma(]0,1[) \cap \widehat{c_{\nu}}[h]$ is nonempty, then the minimal execution path $\gamma$ has at least one point of $\gamma(]0,1[)$ belonging to $\widehat{c_{\nu}}$. If $[c_{\nu}]$ is the carrier of $\gamma$, then $\gamma=\delta_z \phi$ with $z\in \mathbf{D}^{n_{\nu}}\backslash\mathbf{S}^{n_{\nu}-1}$ and $\phi\in \mathcal{G}(1,1)$. We then have 
\[
\bigg\{t\in ]0,1[\mid \gamma(t) \in \widehat{c_{\nu}}[h]\bigg\}=\bigg\{\phi^{-1}(h)\bigg\}.
\]
Otherwise, by Proposition~\ref{noescape}, there is the inclusion $\gamma([0,1]) \subset \de c_{\nu}$ and there exists an execution path $\widehat{g_{\nu}}\delta_z\phi$ for some $z\in \mathbf{S}^{n_{\nu}-1}$ and $\phi\in \mathcal{G}(1,1)$ from the initial state of $c_{\nu}$ to its final state with \[\widehat{g_{\nu}}\delta_z\phi= (\gamma_1\mu_{\ell_1}) * \dots * (\gamma_n\mu_{\ell_n})\] with all $\gamma_i$ minimal and $\gamma\in \{\gamma_1,\dots,\gamma_n\}$. Since $\gamma(]0,1[) \cap \widehat{c_{\nu}}[h]$ is nonempty, we have 
\[
h\in \bigg]\phi(\sum_{j<i}\ell_j),\phi(\sum_{j\leq i}\ell_j)\bigg[
\]
for some $i\in \{1,\dots,n\}$ and $\gamma=\gamma_i$ ($h$ belongs to the interior of the interval because $\gamma(]0,1[)\cap X^0=\varnothing$). We obtain 
\[
\gamma(t) = \widehat{g_{\nu}}\bigg(z,\phi\big(\ell_it+\sum_{j<i}\ell_j\big) \bigg)
\]
for all $t\in [0,1]$. We deduce the equality 
\[\bigg\{t\in ]0,1[\mid \gamma(t) \in \widehat{c_{\nu}}[h]\bigg\}= \bigg\{\frac{\phi^{-1}(h)-\sum_{j<i}\ell_j}{\ell_i}\bigg\}.\]
\epf

\begin{figure}
	\def\n{5}
	\begin{tikzpicture}[black,scale=4,pn/.style={circle,inner sep=0pt,minimum width=4pt,fill=red}]
	\fill [color=gray!15] (0,0) -- (0,1) -- (1,1) -- (1,0) -- cycle;
	\draw[-] [very thick] (0,1) -- (0,0);
	\draw[-] [very thick] (0,0) -- (1,1);
	\foreach \n in {1,2,3,4,5,6,7,8,9,10}
	{\draw[red][->][thick](0,\n/10) -- (\n/10,\n/10);}
	\draw (0,0) node[pn] {} node[black,below left] {$(0\,,0)$};
	\draw (1,1) node[black,above right] {$(1\,,1)$};
	\draw (0,1) node[black,above left] {$(0\,,1)$};
	\draw (1,0) node[black,below right] {$(1\,,0)$};
	\end{tikzpicture}
	\caption{$|X|=[0,1]\p [0,1]$, $X^0=\{0\}\p [0,1] \cup \{(x,x)\mid x\in [0,1]\}$, $\P^{\mathcal{G}}_{(0,t),(t,t)}X= \mathcal{G}(1,1)$ for all $t\in ]0,1]$, $\P^{\mathcal{G}}_{(0,0),(0,0)}X= \{(0,0)\}$ and $\P^{\mathcal{G}}_{\alpha,\beta}X=\varnothing$ otherwise, there is no composable execution paths.}
	\label{contracting}
\end{figure}
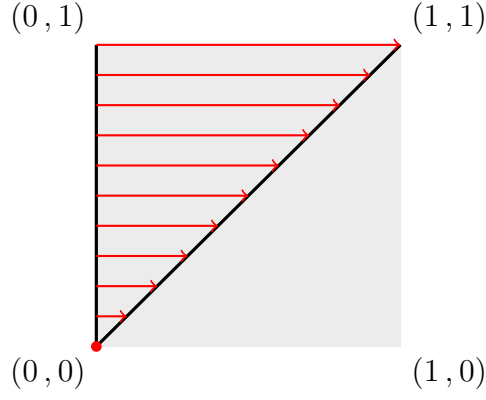

\bp \label{good-achronal}
Let $c_\nu$ be a globular cell of $X_\lambda$ for some $\nu<\lambda$. There exists $b\in ]0,1[$ such that 
\[
\forall h\in ]0,b], \widehat{c_{\nu}}[h]\cap X^0=\varnothing.
\]
\ep

Proposition~\ref{good-achronal} means that, close enough to the initial state of a globular cell, an achronal slice does not contain any state of $X^0$. Similarly, it is possible to prove that close enough to the final state of a globular cell, an achronal slice does not contain any state of $X^0$ either. It is due to the two following geometric facts. The first one is that close enough to the initial state of the globular cell, there is no other states of $X^0$ than the initial state because $X^0$ is discrete. The second one is that a non-constant execution path cannot be deformed in a continuous way to a point in the space of execution paths of a cellular multipointed $d$-space. It is possible in more general multipointed $d$-spaces as the one depicted in Figure~\ref{contracting}. Note that the q-cofibrant replacement of the latter is equal to the disjoint sum of the q-cofibrant replacement of the terminal multipointed $d$-space and of uncountably many directed segments.

\bpf One has $c_\nu\cap X^0=\varnothing$. Consequently, if $\widehat{g_{\nu}}(z,h)\in X^0$ for some $h\in ]0,1[$, then  $z\in\mathbf{S}^{n_\nu-1}$. Thus, if $n_\nu=0$, then $\mathbf{S}^{n_\nu-1}=\varnothing$ and for any $h\in ]0,1[$, one has $\widehat{c_{\nu}}[h]\cap X^0=\varnothing$. Assume now that $n_\nu\geq 1$. Consider the set \[J_1 = \{h\in ]0,1[ \mid \widehat{c_{\nu}}[h] \cap (X^0\backslash \{\widehat{g_{\nu}}(0)\}) \neq \varnothing\}.\] If $J_1$ is nonempty, then consider a sequence $(h^1_n)_{n\geq 0}$ of $J_1$ converging to the greatest lower bound $\inf J_1$ of $J_1$. For all $n\geq 0$, let $z^1_n\in \mathbf{S}^{n_\nu-1}$ such that \[\widehat{g_{\nu}}(z^1_n,h^1_n)\in X^0\backslash \{\widehat{g_{\nu}}(0)\}.\] By extracting a subsequence, we can suppose that the sequence $(z^1_n)_{n\geq 0}$ converges to $z^1_\infty\in \mathbf{S}^{n_\nu-1}$. Since the space $|\globG(\mathbf{D}^{n_{\nu}})|$ is compact, the subspace $\widehat{c_{\nu}}$ is a compact subspace of the weakly Hausdorff space $|X_\lambda|$. The set $\widehat{c_{\nu}} \cap X^0$ is therefore finite because $X^0$ is discrete in $|X_\lambda|$ by Proposition~\ref{p1}. Since \[(\widehat{c_{\nu}} \cap X^0)\backslash \{\widehat{g_{\nu}}(0)\} \subset \widehat{c_{\nu}} \cap X^0\] is discrete finite as well, the sequence $(\widehat{g_{\nu}}(z^1_n,h^1_n))_{n\geq 0}$, which converges to $\widehat{g_{\nu}}(z^1_\infty,\inf J_1)$ by continuity of $\widehat{g_{\nu}}$, eventually becomes constant. Thus, \[\widehat{g_{\nu}}(z^1_\infty,\inf J_1)\neq \widehat{g_{\nu}}(0).\] It implies that \[\inf J_1>0.\] It means that whether $J_1$ is empty or not, there exists $a\in ]0,1[$ such that for all $h\in ]0,a]$, one has $\widehat{c_{\nu}}[h] \cap X^0 \subset \{\widehat{g_{\nu}}(0)\}$. Consider the set \[J_2 = \{h\in ]0,a]\mid \widehat{c_{\nu}}[h] \cap X^0 = \{\widehat{g_{\nu}}(0)\}\}.\] If $J_2$ is nonempty, then consider a sequence $(h^2_n)_{n\geq 0}$ of $J_2$ converging to $\inf J_2$. For all $n\geq 0$, there exists $z^2_n\in \mathbf{S}^{n_\nu-1}$ such that \[\widehat{g_{\nu}}(z^2_n,h^2_n)=\widehat{g_{\nu}}(0).\] By extracting a subsequence, one can suppose that the sequence $(z^2_n)_{n\geq 0}$ of $\mathbf{S}^{n_\nu-1}$ converges to $z^2_\infty\in \mathbf{S}^{n_\nu-1}$. Consider the sequence of globular cells $(c_{\nu_n})_{n\geq 0}$ such that for all $n\geq 0$, the globular cell $c_{\nu_n}$ is the first globular cell appearing in $\carrier(\widehat{g_{\nu}}\delta_{z^2_n})$, i.e. 
\[
\carrier(\widehat{g_{\nu}}\delta_{z^2_n}) = [c_{\nu_n},\underline{c}_n]
\]
where $\underline{c}_n$ is a sequence of globular cells which is necessarily nonempty because $h^2_n<1$. Using Proposition~\ref{cpt-intersect-finite}, we have that the compact subspace $\widehat{c_{\nu}}$ intersects finitely many globular cells of $X_\lambda$. Consequently, by extracting a subsequence again, we can suppose that the sequence of globular cells $(c_{\nu_n})_{n\geq 0}$ is constant and equal to the globular cell $c_{\nu'}$ for some $\nu'<\nu$. Write \[\widehat{g_{\nu}}\delta_{z^2_n} = (\widehat{g_{\nu'}}\delta_{z'_n}\phi_n \mu_{t_n}) * (\gamma_{n}\mu_{1-t_n})\] with, for all $n\geq 0$, 
\[
\begin{aligned}
&0<t_n\leq h^2_n<1,\\
&z'_n\in \mathbf{D}^{n_{\nu'}}\backslash\mathbf{S}^{n_{\nu'}-1},\\
&\phi_n\in \mathcal{G}(1,1),\\
&\widehat{g_{\nu'}}(0) = \widehat{g_{\nu}}(0) \hbox{ (the globular cells $c_\nu$ and $c_{\nu'}$ have the same initial state)},\\
&\widehat{g_{\nu'}}(1) = (\widehat{g_{\nu}}\delta_{z^2_n})(t_n),\\
&\gamma_n\in \P^{\mathcal{G}}_{\widehat{g_{\nu'}}(1),\widehat{g_{\nu}}(1)}X\hbox{ with }\carrier(\gamma_n)=[\underline{c}_n].
\end{aligned}
\]
By extracting a subsequence, one can suppose that the sequence $(z'_n)_{n\geq 0}$ of $\mathbf{D}^{n_{\nu'}}$ converges to $z'_\infty$. Since $\carrier(\widehat{g_{\nu'}}\delta_{z'_\infty})$ exists (it is a sequence of globular cells intersecting $\widehat{c_{\nu'}}$), the execution path $\widehat{g_{\nu'}}\delta_{z'_\infty}$ is not constant. Thus, there exists $T\in ]0,1[$ such that \[\widehat{g_{\nu'}}(z'_\infty,T)\neq \widehat{g_{\nu'}}(0).\] By extracting again a subsequence, one can suppose that the sequence $(t_n\phi_n^{-1}(T))_{n\geq 0}$ of $[0,1]$ converges to $t_\infty$. We have \[\widehat{g_{\nu}}\delta_{z^2_n}(t_n\phi_n^{-1}(T))=\widehat{g_{\nu'}}\delta_{z'_n}\bigg(\phi_n \mu_{t_n}\bigg(t_n\phi_n^{-1}(T)\bigg)\bigg)=\widehat{g_{\nu'}}(z'_n,T)\] for all $n\geq 0$. We obtain by passing to the limit \[\widehat{g_{\nu}}\delta_{z^2_\infty}(t_\infty)=\widehat{g_{\nu'}}(z'_\infty,T).\] We deduce that $\widehat{g_{\nu}}\delta_{z^2_\infty}(t_\infty)\neq \widehat{g_{\nu}}(0)$ and therefore that \[0<t_\infty.\] From the inequalities \[t_n\phi_n^{-1}(T) \leq t_n\leq h^2_n\] for all $n\geq 0$, we obtain by passing to the limit the inequalities \[0<t_\infty\leq \inf J_2.\] It means that whether $J_2$ is empty or not, there exists $b\in ]0,1[$ such that for all $h\in ]0,b]$, one has $\widehat{c_{\nu}}[h] \cap X^0 =\varnothing$. 
\epf

\bth \label{bounded0}
Let $\gamma_\infty$ be an execution path of $X_\lambda$. Let $\nu_0<\lambda$. There exists an open neighborhood $\Omega$ of $\gamma_\infty$ in $\P^{\mathcal{G}}X_\lambda$ such that for all execution paths $\gamma\in \Omega$, the number of copies of $c_{\nu_0}$ in the carrier of $\gamma$ cannot exceed the length of the carrier of $\gamma_\infty$.
\eth

\bpf
Let $\carrier(\gamma_\infty)=[c_{\nu_1},\dots,c_{\nu_n}]$. Consider the decomposition of Theorem~\ref{normal-form}
\[
\gamma_\infty = (\gamma^1_\infty\mu_{\ell_1}) * \dots * (\gamma^n_\infty\mu_{\ell_n})
\]
with $\sum_i\ell_i=1$ and all execution paths $\gamma^i_\infty$ minimal for $i=1,\dots,n$. For $1\leq i \leq n$, let $\nu_i<\lambda$, $\phi_i\in \mathcal{G}(1,1)$ and $z_i\in \mathbf{D}^{n_{\nu_i}}\backslash \mathbf{S}^{n_{\nu_i}-1}$ such that
\[
\begin{aligned}
&\carrier(\gamma^i_\infty)=[c_{\nu_i}],\\
&\gamma^i_\infty(]0,1[)\subset c_{\nu_i},\\
&\gamma^i_\infty=\delta_{z_{i}}\phi_i.
\end{aligned}
\]
Using Proposition~\ref{good-achronal}, pick $h\in ]0,1[$ such that $\widehat{c_{\nu_0}}[h] \cap X^0=\varnothing$. For all $1\leq i\leq n$, the set \[\bigg\{t\in ]0,1[\mid \gamma^i_\infty(t)\in \widehat{c_{\nu_0}}[h]\bigg\}\] contains at most one point $t_{i}$ by Proposition~\ref{middle0}; if the set above is empty, let $t_{i}=\frac{1}{2}$. For all $1\leq i\leq n$, let $L_i$ and $L'_i$ be two real numbers such that \[0<L_i< t_i < L'_i<1.\] For $1\leq i\leq n$, consider the covering of the segment $[\sum_{j<i}\ell_j,\sum_{j\leq i}\ell_j]$ in three nonoverlapping segments of strictly positive length: 
\[\begin{aligned}
&K_i^- = \bigg[\sum_{j<i}\ell_j,\sum_{j<i}\ell_j+\mu_{\ell_i}^{-1}\phi_i^{-1}(L_i)\bigg], \\
&K_i^m = \bigg[\sum_{j< i}\ell_j+\mu_{\ell_i}^{-1}\phi_i^{-1}(L_i),\sum_{j< i}\ell_j+\mu_{\ell_i}^{-1}\phi_i^{-1}(L'_i)\bigg], \\
&K_i^+ = \bigg[\sum_{j< i}\ell_j+\mu_{\ell_i}^{-1}\phi_i^{-1}(L'_i),\sum_{j\leq i}\ell_j\bigg].
\end{aligned}\]
The restriction $\gamma_\infty\rest_{[\sum_{j<i}\ell_j,\sum_{j\leq i}\ell_j]}$ goes from the initial state of the globular cell $c_{\nu_i}$ to its final state. We have therefore \[\gamma_\infty(K_i^m) \subset c_{\nu_i}.\] We deduce
\[
\gamma_\infty(K_i^m) \cap X^0 = \varnothing.
\]
We have
\[
\begin{aligned}
\gamma_\infty(K_i^-) \cap \widehat{c_{\nu_0}}[h] &= \bigg(\gamma_\infty(\{\sum_{j<i}\ell_j\})\cup \gamma_\infty\big(\big]\sum_{j<i}\ell_j,\sum_{j<i}\ell_j+\mu_{\ell_i}^{-1}\phi_i^{-1}(L_i)\big]\big)\bigg) \cap \widehat{c_{\nu_0}}[h]\\
&=\bigg(\{\widehat{g_{\nu_i}}(0)\}\cup \gamma_\infty\big(\big]\sum_{j<i}\ell_j,\sum_{j<i}\ell_j+\mu_{\ell_i}^{-1}\phi_i^{-1}(L_i)\big]\big)\bigg) \cap \widehat{c_{\nu_0}}[h]\\
&\subset \bigg(\{\widehat{g_{\nu_i}}(0)\}\cup \gamma_\infty\big(\big]\sum_{j<i}\ell_j,\sum_{j<i}\ell_j+\mu_{\ell_i}^{-1}\phi_i^{-1}(t_i)\big[\big)\bigg) \cap \widehat{c_{\nu_0}}[h]\\
&=\varnothing,
\end{aligned}
\]
the first equality by formal set identities, the second equality by definition of $\widehat{g_{\nu_i}}(0)$, the inclusion because $L_i<t_i$, and the last equality because $\widehat{g_{\nu_i}}(0)\in X^0$ and by definition of $t_i$. In the same way, we also have
\[
\begin{aligned}
\gamma_\infty(K_i^+) \cap \widehat{c_{\nu_0}}[h] &= \bigg(\gamma_\infty\big(\big[\sum_{j< i}\ell_j+\mu_{\ell_i}^{-1}\phi_i^{-1}(L'_i),\sum_{j\leq i}\ell_j\big[\big)\cup \gamma_\infty(\{\sum_{j\leq i}\ell_j\})\bigg) \cap \widehat{c_{\nu_0}}[h]\\
&=\bigg(\gamma_\infty\big(\big[\sum_{j< i}\ell_j+\mu_{\ell_i}^{-1}\phi_i^{-1}(L'_i),\sum_{j\leq i}\ell_j\big[\big)\cup \{\widehat{g_{\nu_i}}(1)\}\bigg) \cap \widehat{c_{\nu_0}}[h]\\
&\subset \bigg(\gamma_\infty\big(\big]\sum_{j< i}\ell_j+\mu_{\ell_i}^{-1}\phi_i^{-1}(t_i),\sum_{j\leq i}\ell_j\big[\big)\cup \{\widehat{g_{\nu_i}}(1)\}\bigg) \cap \widehat{c_{\nu_0}}[h]\\
&=\varnothing,
\end{aligned}
\]
the first equality by formal set identities, the second equality by definition of $\widehat{g_{\nu_i}}(1)$, the inclusion because $t_i<L'_i$, and the last equality because $\widehat{g_{\nu_i}}(1)\in X^0$ and by definition of $h_i$. Since $|X_\lambda|$ is weakly Hausdorff, the set $\widehat{c_{\nu_0}}[h]$ is a closed subset of $|X_\lambda|$. Moreover, $X^0$ is a closed subset of the space $|X_\lambda|$ as well by Proposition~\ref{p1}. Consequently, the set 
\[
\Omega = \bigcap_{i=1}^{i=n} \bigg( W\bigg(K_i^-,|X_\lambda|\backslash \widehat{c_{\nu_0}}[h]\bigg) \cap  W\bigg(K_i^m,|X_\lambda|\backslash X^0\bigg) \cap W\bigg(K_i^+,|X_\lambda|\backslash \widehat{c_{\nu_0}}[h]\bigg)\bigg),
\]
where \[W([a,b],U)=\{f\in \P^{\mathcal{G}}X_\lambda\mid f([a,b])\subset U\}\] is an open neighborhood of $\gamma_\infty$ of $\P^{\mathcal{G}}X_\lambda$ for the compact-open topology, and therefore for its $\Delta$-kelleyfication which adds open subsets. For all $\gamma\in \Omega$, one has 
\[\gamma(K_i^m)\cap X^0 =\varnothing\]
and, one has \[
\gamma(K_i^-)\cap \widehat{c_{\nu_0}}[h] = \gamma(K_i^+)\cap \widehat{c_{\nu_0}}[h] = \varnothing.
\]
It turns out that the segments of strictly positive length $K_i^-,K_i^m,K_i^+$ for $1\leq i \leq n$ are a finite partition of $[0,1]$ into nonoverlapping segments because we have by definition of the $K_i^-,K_i^m,K_i^+$ for $1\leq i \leq n$:
\[
[0,1] = \bigcup_{i=1}^{i=n} \bigg[\sum_{j<i}\ell_j,\sum_{j\leq i}\ell_j\bigg] = \bigcup_{i=1}^{i=n}\bigg(K_i^- \cup K_i^m \cup K_i^+\bigg).
\] 
Each $c_{\nu_0}$ appearing in $\carrier(\gamma)$ corresponds to a minimal execution path from $\widehat{g_{\nu_0}}(0)$ to $\widehat{g_{\nu_0}}(1)$ (note that these two states can be equal) of the decomposition of $\gamma$ obtained using Theorem~\ref{normal-form}. It necessarily intersects $\widehat{c_{\nu_0}}[h]$. Thus, the number of copies of $c_{\nu_0}$ in the carrier of $\gamma$ cannot exceed the number of $K_i^m$, i.e. the length of the carrier of $\gamma_\infty$. 
\epf

Theorem~\ref{bounded0} does not mean that the carriers of the execution paths of $\Omega$ are of length at most the length of the carrier of $\gamma_\infty$. Indeed, on the segments $K_0^-, K_0^+\cup K_1^-,K_1^+\cup K_2^-,\dots,K_n^+$, an execution path $\gamma$ of $\Omega$ can a priori intersect $X^0$ an arbitrarily large number of times. However, it cannot intersect $\widehat{c_{\nu_0}}[h]$. Therefore these segments do not add copies of $c_{\nu_0}$ in the carrier of $\gamma$.

\bth \label{bounded}
Let $(\gamma_k)_{k\geq 0}$ be a sequence of execution paths of $X_\lambda$ which converges in $\P^{\mathcal{G}}X_\lambda$. Let $c_{\nu_0}$ be a globular cell of $X_\lambda$. Let $i_k$ be the number of times that $c_{\nu_0}$ appears in $\carrier(\gamma_k)$. Then the sequence of integers $(i_k)_{k\geq 0}$ is bounded.
\eth

\bpf 
Write $\gamma_\infty$ for the limit of $(\gamma_k)_{k\geq 0}$ in $\P^{\mathcal{G}}X_\lambda$. By Theorem~\ref{bounded0}, there exists an open $\Omega$ containing $\gamma_\infty$ such that for all $\gamma\in \Omega$, the number of copies of $c_{\nu_0}$ in the carrier of $\gamma$ does not exceed the length of the carrier of $\gamma_\infty$. Since the sequence $(\gamma_k)_{k\geq 0}$ converges to $\gamma_\infty$, there exists $N>0$ such that for all $k\geq N$, $\gamma_k$ belongs to $\Omega$. The proof is complete. 
\epf

\bth \label{calcul_final_structure} Let $0\leq \nu<\lambda$. 
Then every execution path of $X_{\nu+1}$ can be written in a unique way as a finite Moore composition 
\[(f_1\gamma_1\mu_{\ell_1}) * \dots * (f_n\gamma_n\mu_{\ell_n})\] with $n\geq 1$ such that 
\begin{enumerate}[leftmargin=*]
	\item $\sum_i \ell_i = 1$.
	\item $f_i = f$ and $\gamma_i$ is an execution path of $X_\nu$ or $f_i=\widehat{g_\nu}$ and $\gamma_i=\delta_{z_i}\phi_i$ with $z_i\in \mathbf{D}^{n_\nu}\backslash\mathbf{S}^{n_\nu-1}$ and some $\phi\in \mathcal{G}(1,1)$.
	\item for all $1\leq i <n$, either $f_i\gamma_i$ or $f_{i+1}\gamma_{i+1}$ (or both) is (are) of the form $\widehat{g_\nu}\delta_{z}\phi$ for some $z\in \mathbf{D}^{n_\nu}\backslash\mathbf{S}^{n_\nu-1}$ and some $\phi\in \mathcal{G}(1,1)$: intuitively, there is no possible simplification using the Moore composition inside $X_\nu$. 
\end{enumerate}
\eth

\bpf We use the normal form of Theorem~\ref{normal-form} and we use Proposition~\ref{variable-length} to compose successive execution paths of $X_\nu$.
\epf

\section{Chains of globes}
\label{chain}

Let $Z_1,\dots,Z_p$ be $p$ nonempty topological spaces with $p\geq 1$. Consider the multipointed $d$-space 
\[
X=\globG(Z_1)* \dots *\globG(Z_p).
\]
with $p\geq 1$ where the $*$ means that the final state of a globe is identified with the initial state of the next one by reading from the left to the right. Let $\{\alpha_0,\alpha_1,\dots,\alpha_p\}$ be the set of states such that the canonical map $\globG(Z_i)\to X$ takes the initial state $0$ of $\globG(Z_i)$ to $\alpha_{i-1}$ and the final state $1$ of $\globG(Z_i)$ to $\alpha_{i}$.

As a consequence of the associativity of the semimonoidal structure on $\mathcal{G}$-spaces recalled in Theorem~\ref{closedsemimonoidal} and of \cite[Proposition~5.16]{Moore1}, we have

\bp \label{Ftenseur} 
Let $U_1,\dots,U_p$ be $p$ topological spaces with $p\geq 1$. Let $\ell_1,\dots,\ell_p>0$. There is the natural isomorphism of $\mathcal{G}$-spaces 
\[
\mathbb{F}^{\mathcal{G}^{op}}_{\ell_1}U_1\ot \dots \ot \mathbb{F}^{\mathcal{G}^{op}}_{\ell_p}U_p \iso \mathbb{F}^{\mathcal{G}^{op}}_{\ell_1+\dots+\ell_p}(U_1\p \dots\p U_p).
\]
\ep

The case $p=1$ of Proposition~\ref{comp-gl} is treated in Proposition~\ref{calcul-topology-glob} and already used in Proposition~\ref{calculM}. An additional argument is required for the case $p>1$. At first, we prove a lemma which is an addition to Proposition~\ref{morphG-metrizable}.

\begin{lem} \label{inverse}
	The set map $(-)^{-1}:\mathcal{G}(1,p) \to \mathcal{G}(p,1)$ which takes $f:[0,1]\iso^+[0,p]$ to its inverse $f^{-1}:[0,p]\iso^+ [0,1]$ is continuous. 
\end{lem}

\bpf
Since all $\Delta$-generated spaces are sequential, it suffices to prove that $(-)^{-1}:\mathcal{G}(1,p) \to \mathcal{G}(p,1)$ is sequentially continuous. Let $(f_n)_{n\geq 0}$ be a sequence of $\mathcal{G}(1,p)$ which converges to $f\in \mathcal{G}(1,p)$. Let $t\in [0,p]$. Then the sequence $(f_n^{-1}(t))_{n\geq 0}$ of $[0,1]$ has at least one limit point denoted by $L(t)$. By extracting a subsequence of the sequence $(f_n(f_n^{-1}(t)))_{n\geq 0}$, we obtain $f(L(t))=t$, which implies $L(t)=f^{-1}(t)$. Thus every subsequence of $(f_n^{-1}(t))_{n\geq 0}$ has a unique limit point $f^{-1}(t)$. Suppose that the sequence $(f_n^{-1}(t))_{n\geq 0}$ does not converge to $f^{-1}(t)$. Then there exists an open neighborhood $V$ of $f^{-1}(t)$ such that for all $n\geq 0$, $f_n^{-1}(t)\in V^c$ which is compact: contradiction. Therefore the sequence $(f_n^{-1})_{n\geq 0}$ pointwise converges to $f^{-1}$. By Proposition~\ref{morphG-metrizable}, we deduce that the sequence $(f_n^{-1})_{n\geq 0}$ converges to $f^{-1}$.
\epf

\bp  \label{comp-gl}
With the notations of this section. There is a homeomorphism 
\[
\P_{\alpha_0,\alpha_p}^{\mathcal{G}} X \iso  \mathcal{G}(1,p)\p Z_1\p \dots \p Z_p.
\]
\ep

\bpf
The Moore composition of paths induced a map of $\mathcal{G}$-spaces
\[
\P_{0,1} \moore^{\mathcal{G}}\globG(Z_1) \ot \dots \ot \P_{0,1} \moore^{\mathcal{G}}\globG(Z_p) \longrightarrow \P_{\alpha_0,\alpha_p} \moore^{\mathcal{G}}(X).
\]
By Proposition~\ref{calculM}, there is the isomorphism of $\mathcal{G}$-spaces 
\[
\P_{0,1} \moore^{\mathcal{G}}\globG(Z) \iso \mathbb{F}_1^{\mathcal{G}^{op}}Z
\]
for all topological spaces $Z$. We obtain a map of $\mathcal{G}$-spaces
\[
\mathbb{F}_1^{\mathcal{G}^{op}}Z_1 \ot \dots \ot \mathbb{F}_1^{\mathcal{G}^{op}}Z_p \longrightarrow \P_{\alpha_0,\alpha_p} \moore^{\mathcal{G}}(X).
\]
By Proposition~\ref{Ftenseur}, and since $\P_{\alpha_0,\alpha_p}^1 \moore^{\mathcal{G}}(X)=\P_{\alpha_0,\alpha_p}^{\mathcal{G}} X$ by definition of the functor $\moore^{\mathcal{G}}$, we obtain a continuous map 
\[
\begin{cases}
\Psi:&\mathcal{G}(1,p)\p Z_1\p \dots \p Z_p  \longrightarrow \P_{\alpha_0,\alpha_p}^{\mathcal{G}} X\\
&(\phi,z_1,\dots,z_p) \mapsto (\delta_{z_1}\phi_1)*\dots *(\delta_{z_p}\phi_p)
\end{cases}
\]
where $\phi_i\in \mathcal{G}(\ell_i,1)$ with $\sum_i \ell_i=1$ and $\phi=\phi_1\ot \dots \ot \phi_p$ being the decomposition of Proposition~\ref{decomposition-tenseur}. Since all executions paths of globes are one-to-one, the map $\Psi$ above is a continuous bijection. The continuous maps $Z_i\to \{0\}$ for $1\leq i \leq p$ induce by functoriality a map of multipointed $d$-spaces $X \to \vI^{\mathcal{G}}*\dots * \vI^{\mathcal{G}}$ ($p$ times) and then a continuous map \[
\begin{cases}
k:&\P_{\alpha_0,\alpha_p}^{\mathcal{G}}X\longrightarrow \P_{\alpha_0,\alpha_p}^{\mathcal{G}}(\vI^{\mathcal{G}}*\dots * \vI^{\mathcal{G}}) = \mathcal{G}(1,p)\\
&(\delta_{z_1}\phi_1)*\dots *(\delta_{z_p}\phi_p)\mapsto (\delta_{0}\phi_1)*\dots *(\delta_{0}\phi_p) = \phi_1\ot \dots \ot \phi_p.
\end{cases}\] 
Let $i\in \{1,\dots,p\}$. Then we have, with $\gamma=(\delta_{z_1}\phi_1)*\dots *(\delta_{z_p}\phi_p)$, 
\[
\gamma\bigg(k(\gamma)^{-1}(i-\frac{1}{2})\bigg)=\gamma(\phi^{-1}(i-\frac{1}{2}))=\delta_{z_i}\phi_i\phi^{-1}(i-\frac{1}{2}) = (z_i,\frac{1}{2}),
\]
the first equality by definition of $k:\P_{\alpha_0,\alpha_p}^{\mathcal{G}}X\to \mathcal{G}(1,p)$, the second equality since $i-1<i-\frac{1}{2}<i$ and by definition of $\gamma$, and the last equality by definition of the $\phi_i$'s. The set map 
\[
\begin{cases}
&\P_{\alpha_0,\alpha_p}^{\mathcal{G}}X \longrightarrow |\globG(Z_i)|\\
&\gamma \mapsto \gamma\bigg(k(\gamma)^{-1}(i-\frac{1}{2})\bigg)
\end{cases}
\]
is continuous since $k:\P_{\alpha_0,\alpha_p}^{\mathcal{G}}X\to \mathcal{G}(1,p)$ and $(-)^{-1}:\mathcal{G}(1,p)\to \mathcal{G}(p,1)$ are both continuous (see Lemma~\ref{inverse} for the latter map).  Consequently, the set map
\[
\begin{cases}
\overline{k}:&\P_{\alpha_0,\alpha_p}^{\mathcal{G}}X\longrightarrow Z_1\p  \dots \p Z_p\\
& (\delta_{z_1}\phi_1)*\dots *(\delta_{z_p}\phi_p) \mapsto (z_1,\dots,z_p)
\end{cases}
\]
is continuous. It implies that the set map 
\[
\Psi^{-1}=(k,\overline{k}):(\delta_{z_1}\phi_1)*\dots *(\delta_{z_p}\phi_p) \mapsto (\phi_1\ot \dots \ot \phi_p,z_1,\dots,z_p).
\]
is continuous as well and that $\Psi$ is a homeomorphism. 
\epf

Until the end of this section, we work like in Section~\ref{cellular-obj}  with a cellular multipointed $d$-space $X_\lambda$, with the attaching map of the globular cell $c_\nu$ for $\nu<\lambda$ denoted by \[\widehat{g_\nu}:\globG(\mathbf{D}^{n_\nu}) \to X_\lambda.\] Each carrier \[\underline{c}=[c_{\nu_1},\dots,c_{\nu_n}]\] gives rise to a map of multipointed $d$-spaces from a chain of globes to $X_\lambda$
\[
\widehat{g_{\underline{c}}}:\globG(\mathbf{D}^{n_{\nu_1}})*\dots * \globG(\mathbf{D}^{n_{\nu_n}}) \longrightarrow X_\lambda
\]
by ``concatenating'' the attaching maps of the globular cells $c_{\nu_1},\dots,c_{\nu_n}$. Let $\alpha_{i-1}$ ($\alpha_{i}$ resp.) be the initial state (the final state resp.) of $\globG(\mathbf{D}^{n_{\nu_i}})$ for $1\leq i\leq n$ in $\globG(\mathbf{D}^{n_{\nu_1}})* \dots * \globG(\mathbf{D}^{n_{\nu_n}})$. It induces a continuous map 
\[
\P^{\mathcal{G}}\widehat{g_{\underline{c}}}:X_{\underline{c}}:=\P^{\mathcal{G}}_{\alpha_0,\alpha_n}(\globG(\mathbf{D}^{n_{\nu_1}})*\dots * \globG(\mathbf{D}^{n_{\nu_n}})) \longrightarrow \P^{\mathcal{G}} X_\lambda.
\]

\bp \label{unique-param2}
Let $\gamma$ be an execution path of $X_\lambda$. Consider a nondecreasing set map $\phi:[0,1]\to [0,1]$ preserving the extremities such that $\gamma\phi=\gamma$. Then $\phi$ is the identity of $[0,1]$.
\ep

\bpf Note that it is not assumed that $\phi$ is continuous. Suppose that there exist $t<t'$ such that $\phi(t)=\phi(t')$. Then for $t''\in [t,t']$, $\gamma(t'')=\gamma(\phi(t''))=\gamma(\phi(t))$ because $\phi(t)\leq \phi(t'')\leq \phi(t')$, which contradicts the fact that $\gamma$ is locally injective by Proposition~\ref{all-regular}. Thus the set map $\phi$ is strictly increasing. Let $\carrier(\gamma)=[c_{\nu_1},\dots,c_{\nu_n}]$. Let $\gamma = (\gamma_1\mu_{\ell_1}) * \dots * (\gamma_n\mu_{\ell_n})$ with $\ell_1+\dots+\ell_n=1$ such that for all $1\leq i\leq n$, there exist $z_i\in \mathbf{D}^{n_{\nu_i}}\backslash \mathbf{S}^{n_{\nu_i}-1}$ and $\phi_i\in \mathcal{G}(1,1)$ such that for all $t\in ]0,1[$, $\gamma_i(t)=(z_i,\phi_i(t)) \in c_{\nu_i}$, $\gamma_i(0)=\widehat{g_{\nu_i}}(0)$ and $\gamma_i(1)=\widehat{g_{\nu_i}}(1)$. Then 
\[
\{t\in [0,1] \mid \gamma(t)\in X^0\} = \{0=t_0<t_1<\dots<t_n=1\}
\]
with $t_i=\sum_{1\leq j\leq i}\ell_j$ for $0\leq i\leq n$. We deduce that $0=\phi(t_0)<\phi(t_1)<\dots<\phi(t_n)=1$ because the set map $\phi$ is strictly increasing. Since $\gamma(\phi(t_i))=\gamma(t_i)\in X^0$ for $0\leq i\leq n$, one obtains $\phi(t_i)=t_i$ for $0\leq i\leq n$ and $\phi(]t_{i-1},t_{i}[)\subset ]t_{i-1},t_{i}[$ for all $1\leq i\leq n$. Then, observe that
\[
\forall 1\leq i\leq n, \forall t\in]t_{i-1},t_{i}[, (z_i,\phi_i(\phi(t))) = (z_i,\phi_i(t)).
\]
Since $\phi_i$ is bijective, it means that the restriction $\phi\rest_{]t_{i-1},t_{i}[}$ is the identity of $]t_{i-1},t_{i}[$ for all $1\leq i\leq n$. 
\epf

\begin{nota}
	Let $\phi$ be a set map from a segment $[a,b]$ to a segment $[c,d]$. Let \[
	\begin{aligned}
	&\phi(x^-) = \sup \{\phi(t)\mid t<x\},\\
	&\phi(x^+) = \inf \{\phi(t)\mid x<t\}.
	\end{aligned}
	\]
\end{nota}

\bth \label{unique-param3}
Let $\gamma_1$ and $\gamma_2$ be two execution paths of $X_\lambda$ such that there exist two nondecreasing set maps $\phi_1,\phi_2:[0,1]\to [0,1]$ preserving the extremities such that \[
\begin{aligned}
&\forall t\in [0,1],\gamma_1(\phi_1(t)) = \gamma_2(t)\\
&\forall t\in [0,1],\gamma_1(t) = \gamma_2(\phi_2(t)).
\end{aligned}
\]
Then $\phi_1,\phi_2\in \mathcal{G}(1,1)$ and $\phi_2=\phi_1^{-1}$. 
\eth

\bpf Note that it is not assumed that $\phi_1$ and $\phi_2$ are continuous. For all $t\in [0,1]$, we have  $\gamma_1(\phi_1(\phi_2(t))) = \gamma_2(\phi_2(t)) =\gamma_1(t)$. Using Proposition~\ref{unique-param2}, we deduce that $\phi_1\phi_2=\id_{[0,1]}$. In the same way, we have $\phi_2\phi_1=\id_{[0,1]}$. This proves that $\phi_1$ and $\phi_2$ are two bijective set maps preserving the extremities which are inverse to each other. Suppose e.g. that there exists $t\in [0,1]$ such that $\phi_1(t^-)<\phi_1(t)$. Then $\phi_1$ cannot be surjective: contradiction. By using similar arguments, we deduce that for all $t\in [0,1]$, $\phi_1(t^-)=\phi_1(t)=\phi_1(t^+)$ and $\phi_2(t^-)=\phi_2(t)=\phi_2(t^+)$. Consequently, the set maps $\phi_1$ and $\phi_2$ are continuous.
\epf

\bp
Let $\underline{c}$ be the carrier of some execution path of $X_\lambda$. Every execution path of the image of $\P^{\mathcal{G}}\widehat{g_{\underline{c}}}$ is of the form 
\[
(\widehat{g_{\nu_1}}\delta_{z^{1}} * \dots * \widehat{g_{\nu_n}}\delta_{z^{n}})\phi
\]
with $\phi\in \mathcal{G}(1,n)$ and $z^i\in \mathbf{D}^{n_{\nu_i}}$ for $1\leq i \leq n$.  
\ep

\bpf The first assertion is a consequence of the definition of $\widehat{g_{\underline{c}}}$ and of Proposition~\ref{comp-gl}. 
\epf

\begin{nota}
	Let $\underline{c}$ be the carrier of some execution path of $X_\lambda$. Using the identification provided by the homeomorphism of Proposition~\ref{comp-gl}, we can use the notation
	\[
	(\P^{\mathcal{G}}\widehat{g_{\underline{c}}})(\phi,z^1,\dots,z^n) = (\widehat{g_{\nu_1}}\delta_{z^{1}} * \dots * \widehat{g_{\nu_n}}\delta_{z^{n}})\phi. 
	\]
\end{nota}

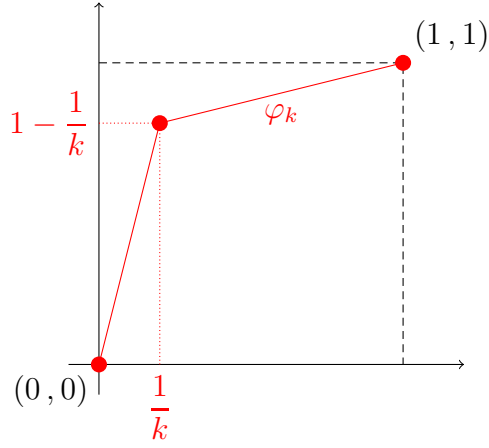
\begin{figure}
	\def\n{5}
	\begin{tikzpicture}[black,scale=4,pn/.style={circle,inner sep=0pt,minimum width=6pt,fill=red}]
	\draw[->] (-.1,0) -- (1.2,0);
	\draw[->] (0,-.1) -- (0,1.2);
	\draw[densely dashed] (0,1) -| (1,0);
	\draw[red] (0,0) node[pn] {}  node[black,below left] {$(0\,,0)$} -- (1/\n,1-1/\n) node[pn] {} -- (1,1) node[midway,below] {$\varphi_k$} node[pn] {} node[black,above right] {$(1\,,1)$};
	\draw[red,densely dotted] (0,1-1/\n) node[left] {$1-\dfrac{1}{k}$} -| (1/\n,0) node[below] {$\dfrac{1}{k}$};
	\end{tikzpicture}
	\caption{A sequence $(\phi_k)_{k\geq 1}$ of $\mathcal{G}(1,1)$ without limit point}
	\label{phi}
\end{figure}

Before proving the main theorem of this section, we need the following topological lemmas:

\begin{lem} \label{prod-seq0}
	Let $U_1,\dots,U_p$ be $p$ first-countable $\Delta$-Hausdorff $\Delta$-generated spaces with $p\geq 1$. Then the product $U_1\p \dots \p U_p$ in the category $\mathcal{T\!O\!P}$ of general topological spaces and continuous maps coincides with the product in $\top$. 
\end{lem}

\bpf
Consider $U_1\p \dots \p U_p$ equipped with the product topology in the category $\mathcal{T\!O\!P}$ of general topological spaces and continuous maps. This topology is first-countable as a finite product of first-countable topologies. Each space $U_i$ is locally path-connected, being $\Delta$-generated. Thus, the finite product $U_1\p \dots \p U_p$ equipped with the product topology in $\mathcal{T\!O\!P}$ is locally path-connected. We deduce that $U_1\p \dots \p U_p$ equipped with the product topology in $\mathcal{T\!O\!P}$ is $\Delta$-generated: the $\Delta$-kelleyfication functor is not required. Moreover since each $U_i$ is $\Delta$-Hausdorff, the product in $\mathcal{T\!O\!P}$ is $\Delta$-Hausdorff as well. It means that $U_1\p \dots \p U_p$ equipped with the product topology in $\mathcal{T\!O\!P}$ coincides with the product in $\top$.
\epf

\begin{lem} \label{prod-seq}
	Let $U_1,\dots,U_p$ be $p$ first-countable $\Delta$-Hausdorff $\Delta$-generated spaces with $p\geq 1$. Let $(u_n^i)_{n\geq 0}$ be a sequence of $U_i$ for $1\leq i \leq p$ which converges to $u_\infty^i\in U_i$. Then the sequence $((u_n^1,\dots,u_n^p))_{n\geq 0}$ converges to $(u_\infty^1,\dots,u_\infty^p)\in U_1\p \dots \p U_p$ for the product calculated in $\top$.
\end{lem}

Note that the converse is obvious: if the sequence $((u_n^1,\dots,u_n^p))_{n\geq 0}$ converges to $(u_\infty^1,\dots,u_\infty^p)\in U_1\p \dots \p U_p$, then the sequences $(u_n^i)_{n\geq 0}$ converge to $u_\infty^i\in U_i$ for all $1\leq i\leq p$ because of the existence of the projection maps $U_1\p \dots \p U_p\to U_i$ for all $1\leq i\leq p$. A sequence converges to some point in a $\Delta$-generated space if and only if the corresponding application from the one-point compactification $\overline{\mathbb{N}}=\mathbb{N}\cup\{\infty\}$ of the discrete space $\mathbb{N}$ to the $\Delta$-generated space is continuous. The point is that the one-point compactification of $\mathbb{N}$ is not $\Delta$-generated: its $\Delta$-kelleyfication is a discrete space. Therefore it does not seem possible to use the universal property of the finite product in $\top$ to prove Lemma~\ref{prod-seq}.

\bpf
Each convergent sequence gives rise to a continuous map $\overline{\mathbb{N}}\to U_i$ for $1\leq i \leq p$. We obtain a continuous map $\overline{\mathbb{N}}\to U_1\p \dots \p U_p$ by using the universal property of the finite product in $\mathcal{T\!O\!P}$ thanks to Lemma~\ref{prod-seq0} and the proof is complete. 
\epf

The sequence $(\phi_k)_{k\geq 1}$ of $\mathcal{G}(1,1)$ depicted in Figure~\ref{phi} has no limit point because the only possibility is the set map which takes $0$ to $0$ and the other points of $[0,1]$ to $1$: it does not belong to $\mathcal{G}(1,1)$. Thus, the topological space $\mathcal{G}(1,n)$, which is homeomorphic to $\mathcal{G}(1,1)$ for all $n\geq 1$, is not sequentially compact. However, Theorem~\ref{img-closed} holds anyway.

\bth \label{img-closed} Let $\underline{c}$ be the carrier of some execution path of $X_\lambda$. 
\begin{enumerate}[leftmargin=*]
	\item Consider a sequence $(\gamma_k)_{k\geq 0}$ of the image of $\P^{\mathcal{G}}\widehat{g_{\underline{c}}}$ which converges pointwise to $\gamma_\infty$ in $\P^{\mathcal{G}}X_\lambda$. Let \[\gamma_k=(\P^{\mathcal{G}}\widehat{g_{\underline{c}}})(\phi_k,z_k^1,\dots,z_k^n)\]
	with $\phi_k\in \mathcal{G}(1,n)$ and $z_k^i\in \mathbf{D}^{n_{\nu_i}}$ for $1\leq i \leq n$ and $k\geq 0$. Then there exist $\phi_\infty\in \mathcal{G}(1,n)$ and $z_\infty^i\in \mathbf{D}^{n_{\nu_i}}$ for $1\leq i \leq n$ such that 
	\[\gamma_\infty=(\P^{\mathcal{G}}\widehat{g_{\underline{c}}})(\phi_\infty,z_\infty^1,\dots,z_\infty^n)\]
	and such that $(\phi_\infty,z_\infty^1,\dots,z_\infty^n)$ is a limit point of the sequence $((\phi_k,z_k^1,\dots,z_k^n))_{k\geq 0}$.
	\item The image of $\P^{\mathcal{G}}\widehat{g_{\underline{c}}}$ is closed in $\P^{\mathcal{G}} X_\lambda$.
\end{enumerate}
\eth

\bpf (1) By a Cantor diagonalization argument, we can suppose that 
\begin{itemize}[leftmargin=*]
	\item The sequence $(z_k^i)_{k\geq 0}$ converges to $z_\infty^i\in \mathbf{D}^{n_{\nu_i}}$ for each $1\leq i \leq n$.
	\item The sequence $(\phi_k(r))_{k\geq 0}$ converges to a real number denoted by $\phi_\infty(r)\in [0,m]$ for each $r\in [0,1]\cap \mathbb{Q}$.
	\item The sequence $(\phi_k^{-1}(r))_{k\geq 0}$ converges to a real number denoted by $\phi_\infty^{-1}(r)\in [0,1]$ for each $r\in [0,n]\cap \mathbb{Q}$.
\end{itemize}
Since the sequence of execution paths $(\gamma_k)_{k\geq 0}$ converges pointwise to $\gamma_\infty$, we obtain 
\[
\gamma_\infty(r) = (\widehat{g_{\nu_1}}\delta_{z_\infty^{1}} * \dots * \widehat{g_{\nu_n}}\delta_{z_\infty^{n}})(\phi_\infty(r))
\]
for all $r\in [0,1]\cap \mathbb{Q}$ and 
\[
\gamma_\infty(\phi_\infty^{-1}(r)) = (\widehat{g_{\nu_1}}\delta_{z_\infty^{1}} * \dots * \widehat{g_{\nu_n}}\delta_{z_\infty^{n}})(r)
\]
for all $r\in [0,n]\cap \mathbb{Q}$. For $r_1<r_2 \in [0,1]\cap \mathbb{Q}$, $\phi_k(r_1)<\phi_k(r_2)$ for all $k\geq 0$. Therefore by passing to the limit, we obtain $\phi_\infty(r_1)\leq \phi_\infty(r_2)$. Note that $\phi_\infty(0)=0$ and $\phi_\infty(1)=n$ since $0,1\in \mathbb{Q}$. In the same way, we see that $\phi_\infty^{-1}:[0,n]\cap \mathbb{Q}\to [0,1]$ is nondecreasing and that $\phi_\infty^{-1}(0)=0$ and $\phi_\infty^{-1}(n)=1$. For $t\in ]0,1[$, let us extend the definition of $\phi_\infty$ as follows: 
\[
\phi_\infty(t) = \sup \{\phi_\infty(r)\mid r\in ]0,t]\cap \mathbb{Q}\}.
\]
The upper bound exists since $\{\phi_\infty(r)\mid r\in ]0,t]\cap \mathbb{Q}\}\subset [0,n]$. For each $t\in [0,1]\backslash \mathbb{Q}$, there exists a nondecreasing sequence $(r_k)_{k\geq 0}$ of rational numbers converging to $t$. Then
\[
\lim_{k\to\infty} \phi_\infty(r_k) = \phi_\infty(t).
\]
By continuity, we deduce that \[
\gamma_\infty(t) = (\widehat{g_{\nu_1}}\delta_{z_\infty^{1}} * \dots * \widehat{g_{\nu_n}}\delta_{z_\infty^{n}})(\phi_\infty(t))
\]
for all $t\in [0,1]$. It is easy to see that the set map $\phi_\infty:[0,1]\to [0,n]$ is nondecreasing and that it preserves extremities. For $t\in ]0,1[$, extend the definition of $\phi_\infty^{-1}$ as well as follows: 
\[
\phi_\infty^{-1}(t) = \sup \{\phi_\infty^{-1}(r)\mid r\in ]0,t]\cap \mathbb{Q}\}.
\]
The upper bound exists since $\{\phi_\infty^{-1}(r)\mid r\in ]0,t]\cap \mathbb{Q}\}\subset [0,1]$. For each $t\in [0,n]\backslash \mathbb{Q}$, there exists a nondecreasing sequence $(r_k)_{k\geq 0}$ of rational numbers converging to $t$. Then
\[
\lim_{k\to\infty} \phi_\infty^{-1}(r_k) = \phi_\infty^{-1}(t).
\]
By continuity, we deduce that \[
\gamma_\infty(\phi_\infty^{-1}(t)) = (\widehat{g_{\nu_1}}\delta_{z_\infty^{1}} * \dots * \widehat{g_{\nu_n}}\delta_{z_\infty^{n}})(t)
\]
for all $t\in [0,n]$. It is easy to see that the set map $\phi_\infty^{-1}:[0,n]\to [0,1]$ is nondecreasing and that it preserves extremities. We obtain for all $t\in [0,1]$
\[
\begin{aligned}
&\gamma_\infty(t) = (\widehat{g_{\nu_1}}\delta_{z_\infty^{1}} * \dots * \widehat{g_{\nu_n}}\delta_{z_\infty^{n}})(\mu_n^{-1}\mu_n\phi_\infty(t))\\
&\gamma_\infty(\phi_\infty^{-1}\mu_n^{-1}(t)) = (\widehat{g_{\nu_1}}\delta_{z_\infty^{1}} * \dots * \widehat{g_{\nu_n}}\delta_{z_\infty^{n}})(\mu_n^{-1}(t)).
\end{aligned}
\]
Using Theorem~\ref{unique-param3}, we obtain that $\mu_n\phi_\infty:[0,1]\to [0,1]$ and $\phi_\infty^{-1}\mu_n^{-1}:[0,1]\to [0,1]$ are homeomorphisms which are inverse to each other. We deduce that $\phi_\infty:[0,1]\to [0,n]$ and $\phi_\infty^{-1}:[0,n]\to [0,1]$ are homeomorphisms which are inverse to each other. Let $t\in [0,1]\backslash \mathbb{Q}$. Since the sequence $(\phi_k(t))_{k\geq 0}$ belongs to the sequential compact $[0,n]$, it has at least one limit point $\ell$. There exists a subsequence of $(\phi_k(t))_{k\geq 0}$ which converges to $\ell$. We obtain 
\[
\forall r\in [0,t]\cap \mathbb{Q},\forall r'\in [t,1]\cap \mathbb{Q}, \phi_\infty(r) \leq \ell\leq \phi_\infty(r').
\]
Since $\phi_\infty\in \mathcal{G}(1,n)$ and by density of $\mathbb{Q}$, we deduce that $\ell=\phi_\infty(t)$ necessarily. Now suppose that the sequence $(\phi_k(t))_{k\geq 0}$ does not converge to $\phi_\infty(t)$. Then there exists an open neighborhood $V$ of $\phi_\infty(t)$ in $[0,n]$ such that for all $k\geq 0$, $\phi_k(t)\notin V$. We deduce that the sequence $(\phi_k(t))_{k\geq 0}$ of $[0,n]$ has no limit point: contradiction. We have proved that the sequence $(\phi_k)_{k\geq 0}$ converges pointwise to $\phi_\infty$. Using Proposition~\ref{morphG-metrizable}, we deduce that $(\phi_k)_{k\geq 0}$ converges uniformly to $\phi_\infty$. We deduce that $(\phi_\infty,z_\infty^1,\dots,z_\infty^n)$ is a limit point of the sequence $((\phi_k,z_k^1,\dots,z_k^n))_{k\geq 0}$ in $\mathcal{G}(1,n)\p \mathbf{D}^{n_{\nu_1}}\p \dots \p \mathbf{D}^{n_{\nu_n}}$ by Proposition~\ref{morphG-metrizable} and Lemma~\ref{prod-seq}.

(2) Let $(\P^{\mathcal{G}}\widehat{g_{\underline{c}}}(\Gamma_n))_{n\geq 0}$ be a sequence of $(\P^{\mathcal{G}}\widehat{g_{\underline{c}}})(X_{\underline{c}})$ which converges in $\P^{\mathcal{G}} X_\lambda$. The limit $\gamma_\infty\in \P^{\mathcal{G}} X_\lambda$ of the sequence of execution paths $(\P^{\mathcal{G}}\widehat{g_{\underline{c}}}(\Gamma_n))_{n\geq 0}$ is also a pointwise limit. We can suppose by extracting a subsequence that the sequence $(\Gamma_n)_{n\geq 0}$ of $X_{\underline{c}}$ is convergent in $X_{\underline{c}}$. Thus, by continuity of $\P^{\mathcal{G}}\widehat{g_{\underline{c}}}$, we obtain $\gamma_\infty=(\P^{\mathcal{G}}\widehat{g_{\underline{c}}})(\Gamma_\infty)$ for some $\Gamma_\infty\in X_{\underline{c}}$. We deduce that $\P^{\mathcal{G}}\widehat{g_{\underline{c}}}(X_{\underline{c}})$ is sequentially closed in $\P^{\mathcal{G}} X_\lambda$. Since $\P^{\mathcal{G}} X_\lambda$ is sequential, being a $\Delta$-generated space, the proof is complete. 
\epf

\begin{cor} \label{Dini-cellular}
Suppose that $X_\lambda$ is a finite cellular multipointed  $d$-space without loops. Then the topology of $\P^{\mathcal{G}}X_\lambda$ is the topology of the pointwise convergence which is therefore $\Delta$-generated.
\end{cor}

We do not know whether the ``without loops'' hypothesis can be removed and whether finite can be replaced by locally finite.

\bpf Let $(\gamma_n)_{n\geq 0}$ be a sequence of execution paths of $X_\lambda$ which pointwise converges to $\gamma_\infty$. Since $X_\lambda$ is finite and without loop, the set 
\[
\mathcal{T} = \{\carrier(\gamma)\mid \gamma\in \P^{\mathcal{G}}X_\lambda\}
\]
is finite. We obtain a finite covering by (closed) subsets
\[
\P^{\mathcal{G}}X_\lambda = \bigcup_{\underline{c}\in \mathcal{T}}(\P \widehat{g_{\underline{c}}})(X_{\underline{c}})
\]
because each execution path has a carrier by Theorem~\ref{normal-form}. Suppose that $(\gamma_n)_{n\geq 0}$ does not converge to $\gamma_\infty$ in $\P^{\mathcal{G}}X_\lambda$. Then there exists an open neighborhood $V$ of $\gamma_\infty$ in $\P^{\mathcal{G}}X_\lambda$ such that $\forall n\geq 0, \gamma_n\notin V$. Since $\mathcal{T}$ is finite, one can suppose by extracting a subsequence that
\[
\exists \underline{c}\in \mathcal{T}, \forall n\geq 0, \gamma_n \in (\P \widehat{g_{\underline{c}}})(X_{\underline{c}}).
\]
By Theorem~\ref{img-closed}, the sequence $(\gamma_n)_{n\geq 0}$ has a limit point which belongs to the complement of $V$ which is closed. This limit point is necessarily the pointwise limit $\gamma_\infty$. We obtain $\gamma_\infty\notin V$: contradiction.
\epf

Corollary~\ref{Dini-cellular} can be viewed as a second Dini theorem for the space of execution paths of a finite cellular multipointed $d$-space without loops. Indeed, if $X_\lambda=\vI^{\mathcal{G}}$ (the directed segment), then $\P^{\mathcal{G}}X_\lambda=\mathcal{G}(1,1)$ and we recover the fact that the topology of $\mathcal{G}(1,1)$ coincides with the pointwise convergence topology by Proposition~\ref{morphG-metrizable}.

\section{The unit and the counit of the adjunction on q-cofibrant objects}
\label{unit}

Consider in this section the following situation: a pushout diagram of multipointed $d$-spaces 
\[
\xymatrix@C=4em@R=4em
{
	\globG(\mathbf{S}^{n-1}) \fd{} \fr{g} & A \ar@{->}[d]^-{f} \\
	\globG(\mathbf{D}^{n}) \fr{\widehat{g}} & \cocartesien X
}
\]
with $n\geq 0$ and $A$ cellular. Note that $A^0=X^0$. Let $D= \mathbb{F}^{\mathcal{G}^{op}}_1\mathbf{S}^{n-1}$ and $E=\mathbb{F}^{\mathcal{G}^{op}}_1\mathbf{D}^{n}$. Consider the Moore flow $\overline{X}$ defined by the pushout diagram of Figure~\ref{Xoverline} where the two equalities 
\[
\begin{aligned}
&\moore^{\mathcal{G}}(\globG(\mathbf{S}^{n-1})) = \globP(D)\\
&\moore^{\mathcal{G}}(\globG(\mathbf{D}^{n})) = \globP(E)
\end{aligned}
\]
come from Proposition~\ref{calculM} and where the map $\psi$ is induced by the universal property of the pushout. 

\begin{figure}
\[
\xymatrix@C=4em@R=4em
{
	\moore^{\mathcal{G}}(\globG(\mathbf{S}^{n-1})) = \globP(D) \fd{} \fr{\moore^{\mathcal{G}}(g)} & \moore^{\mathcal{G}}(A) \ar@/^20pt/@{->}[rdd]^-{\moore^{\mathcal{G}}(f)}\ar@{->}[d]^-{\overline{f}} \\
	\moore^{\mathcal{G}}(\globG(\mathbf{D}^{n})) = \globP(E) \fr{\overline{g}} \ar@/_20pt/@{->}[rrd]_-{\moore^{\mathcal{G}}(\widehat{g})} & \cocartesien \overline{X}\ar@{-->}[rd]|-{\psi}&  \\
	&& \moore^{\mathcal{G}}(X).
}
\]
	\caption{Definition of $\overline{X}$}
	\label{Xoverline}
\end{figure}

The $\mathcal{G}$-space of execution paths of the Moore flow $\overline{X}$ can be calculated by introducing a diagram of $\mathcal{G}$-spaces $\mathcal{D}^{f}$ over a Reedy category $\mathcal{P}^{g(0),g(1)}(A^0)$ whose definition is recalled in Appendix~\ref{reedy}. Let $T$ be the $\mathcal{G}$-space defined by the pushout diagram of $\topdgr_0$
\[
\xymatrix@C=4em@R=4em
{
	D  \fd{} \fr{\P\moore^{\mathcal{G}}(g)} & \P_{g(0),g(1)} \moore^{\mathcal{G}}(A) \ar@{->}[d]^-{\P\overline{f}} \\
	E  \fr{\P\overline{g}} & \cocartesien T.
}
\]

Consider the diagram of spaces $\D^{f}:\mathcal{P}^{g(0),g(1)}(A^0)\to \topdgr_0$ defined as follows:
\[
\D^f((u_0,\epsilon_1,u_1),(u_1,\epsilon_2,u_2),\dots ,(u_{n-1},\epsilon_n,u_n)) = Z_{u_0,u_1}\ot Z_{u_1,u_2} \ot \dots \ot Z_{u_{n-1},u_n}
\]
with 
\[
Z_{u_{i-1},u_i}=
\begin{cases}
\P_{u_{i-1},u_i}\moore^{\mathcal{G}}(A) & \hbox{if }\epsilon_i=0\\
T & \hbox{if }\epsilon_i=1
\end{cases}
\] 
In the case $\epsilon_i=1$, $(u_{i-1},u_i)=(g(0),g(1))$ by definition of $\mathcal{P}^{g(0),g(1)}(A^0)$. The inclusion maps $I_i's$ are induced by the map $\P\overline{f}:\P_{g(0),g(1)} \moore^{\mathcal{G}}(A) \to T$. The composition maps $c_i's$ are induced by the compositions of paths of the Moore flow $\moore^{\mathcal{G}}(A)$. 

\bth \cite[Theorem~9.7]{Moore1}
We obtain a well-defined diagram of $\mathcal{G}$-spaces \[\D^f:\mathcal{P}^{g(0),g(1)}(A^0)\to \topdgr_0.\] There is the isomorphism of $\mathcal{G}$-spaces $\liminj \D^f \iso \P \overline{X}$. 
\eth

By the universal property of the pushout, we obtain a canonical map of $\mathcal{G}$-spaces \[\P \psi:\liminj \D^f \longrightarrow \P\mathbb{M}^\mathcal{G}X.\] 
The goal of Theorem~\ref{pre-calculation-pathspace} and of Theorem~\ref{pre-calculation-pathspace2} is to prove that the canonical map of $\mathcal{G}$-spaces $\P\psi$ is an isomorphism of $\mathcal{G}$-spaces. The proof is twofold: at first, it is proved that it is an objectwise continous bijection, and then that it is an objectwise homeomorphism.

\bth \label{pre-calculation-pathspace}
Under the hypotheses and the notations of this section. The map of $\mathcal{G}$-spaces \[\P\psi:\liminj \D^f \longrightarrow \P\mathbb{M}^\mathcal{G}(X)\] is an objectwive bijection.
\eth

\bpf Throughout the proof, the reader must keep in mind that for any map of multipointed $d$-spaces \[q:X_1\longrightarrow X_2\] and for any execution path \[\gamma\in \P^\ell X_1 = \P^\ell \moore^{\mathcal{G}}(X_1)\] of length $\ell$ of the multipointed $d$-space $X_1$, or equivalently of the Moore flow $\moore^{\mathcal{G}}(X_1)$, there is by Definition~\ref{path-with-length} and Theorem~\ref{adj-multi-moore} the tautological equality
\[
\bigg(\P\moore^{\mathcal{G}}(q)\bigg)(\gamma)=|q|.\gamma,
\]
the right-hand term meaning the composite of the underlying continuous map $|q|:|X_1|\to |X_2|$ between the underlying spaces of $X_1$ and $X_2$ with the execution path $\gamma:[0,1]\to |X_1|$. It will be denoted $q\gamma$ or $q.\gamma$, as it was always done so far. In other terms, we will be using the abuse of notation 
\[
\P\moore^{\mathcal{G}}(q) =q
\]
for any map of multipointed $d$-spaces $q$. The reader must also keep in mind that if $\gamma\in \P^\ell X_1$ and $\gamma'\in\P^{\ell'}X_1$ are two composable execution paths of $X_1$ of length $\ell$ and $\ell'$ respectively, then the Moore composition of execution paths (cf. Proposition~\ref{Moore-comp}) \[\gamma*\gamma'\in \P^{\ell+\ell'} X_1\] is also by Theorem~\ref{adj-multi-moore} the composition of paths in $\moore^{\mathcal{G}}(X_1)$ for tautological reasons.

The proof of this theorem is divided in several parts.

$\bullet$ \underline{\textit{Objectwise calculation}}. 

It suffices to prove that the continuous map 
\[\P^1\psi:\liminj \D^f(1) \longrightarrow \P^1\mathbb{M}^\mathcal{G}(X)=\P^\mathcal{G}X.\]
is a bijection to complete the proof since all objects of the reparametrization category $\mathcal{G}$ are isomorphic and since colimits of $\mathcal{G}$-spaces are calculated objectwise. 

$\bullet$ \underline{\textit{The final topology}}.
If we can prove that the continuous \[\P^1\psi:\liminj \D^f(1) \longrightarrow \P^\mathcal{G}X\] is a bijection with the colimit $\liminj \D^f(1)$ equipped with the final topology, then the proof will be complete even in the category of $\Delta$-Hausdorff $\Delta$-generated topological spaces because of the following facts: 
\begin{itemize}[leftmargin=*]
	\item Let $i:A\to B$ be a continuous one-to-one map between $\Delta$-generated spaces such that $B$ is also $\Delta$-Hausdorff, then $A$ is $\Delta$-Hausdorff: let $f:[0,1]\to A$ be a continuous map; then $f$ being one-to-one, $f([0,1])=i^{-1}((i.f)([0,1]))$ is closed.
	\item The space $\P^{\mathcal{G}}X$ is, by definition, equipped with the $\Delta$-kelleyfication of the relative topology induced by the inclusion of set $\P^{\mathcal{G}}X\subset \ttop([0,1],|X|)$.
	\item If we work in the category of $\Delta$-Hausdorff $\Delta$-generated topological spaces, then the space $\ttop([0,1],|X|)$ will be $\Delta$-Hausdorff, hence the space $\P^{\mathcal{G}}X$ will be $\Delta$-Hausdorff, and therefore $\liminj \D^f$ equipped with the final topology will be $\Delta$-Hausdorff as well.
\end{itemize}

$\bullet$ \underline{\textit{Surjectivity of $\P^1\psi$}}. The map $\psi$ of Figure~\ref{Xoverline} is induced by the universal property of the pushout. It is bijective on states. The multipointed $d$-space $X$ is equipped with the $\Omega$-final structure because it is defined as a colimit in the category of multipointed $d$-spaces. By Theorem~\ref{final-structure-revisited}, every execution path of $X$ is therefore a Moore composition of the form 
\[
(f_1\gamma_1\mu_{\ell_1}) * \dots * (f_n\gamma_n\mu_{\ell_n})
\]
such that $f_i\in \{f,\widehat{g}\}$ for all $1\leq i \leq n$ with 
\[
\begin{cases}
\gamma_i\in {\P^{\mathcal{G}}\globG(\mathbf{D}^{n})} & \hbox{ if }f_i=\widehat{g}\\
\gamma_i\in {\P^{\mathcal{G}}A} & \hbox{ if }f_i=f
\end{cases}
\]
with $\sum_i \ell_i = 1$. Then for all $1\leq i \leq n$, $\gamma_i\mu_{\ell_i}\in \P^{\ell_i}\globG(\mathbf{D}^{n}) = \P^{\ell_i}\moore^{\mathcal{G}}(\globG(\mathbf{D}^{n}))$ or $\gamma_i\mu_{\ell_i}\in \P^{\ell_i}A=\P^{\ell_i}\moore^{\mathcal{G}}(A)$. It gives rise to the execution path 
\[
\P\overline{f_1}(\gamma_1\mu_{\ell_1}) * \dots * \P\overline{f_n}(\gamma_n\mu_{\ell_n})
\]
with 
\[
\begin{cases}
\overline{f_i}= \overline{g} & \hbox{ if }f_i=\widehat{g}\\
\overline{f_i}= \overline{f} & \hbox{ if }f_i=f
\end{cases}
\]
of length $1$ of the Moore flow $\overline{X}$. By the commutativity of Figure~\ref{Xoverline}, we obtain the equality 
\[
(f_1\gamma_1\mu_{\ell_1}) * \dots * (f_n\gamma_n\mu_{\ell_n}) = (\P^1 \psi)\bigg(\P\overline{f_1}(\gamma_1\mu_{\ell_1}) * \dots * \P\overline{f_n}(\gamma_n\mu_{\ell_n})\bigg).
\]
It means that the map of Moore flows $\psi:\overline{X} \to \moore^{\mathcal{G}}(X)$ induces a surjective continuous map from $\P^1\overline{X}$ to $\P^1\moore^{\mathcal{G}}(X)=\P^{\mathcal{G}}X$. In other terms, the map $\P^1\psi$ is a surjection. 

$\bullet$ \underline{\textit{The map $\widehat{\mathbb{C}}$}}. Consider the diagram of topological spaces \[\mathcal{E}^f:\mathcal{P}^{g({0}),g({1})}(A^0)\to \top\] defined as follows:
\[
\mathcal{E}^f((u_0,\epsilon_1,u_1),(u_1,\epsilon_2,u_2),\dots ,(u_{n-1},\epsilon_n,u_n)) = \displaystyle\bigsqcup_{\substack{(\ell_1,\dots,\ell_n)\\\ell_1+\dots+\ell_n=1}}  Z_{u_0,u_1}(\ell_1) \p \dots \p Z_{u_{n-1},u_n}(\ell_n)
\]
with 
\[
Z_{u_{i-1},u_i}(\ell_i)=
\begin{cases}
\P_{u_{i-1},u_i}^{\ell_i}\moore^{\mathcal{G}}(A) = \P_{u_{i-1},u_i}^{\ell_i}A & \hbox{if }\epsilon_i=0\\
T(\ell_i) & \hbox{if }\epsilon_i=1 \hbox{ ($\Rightarrow(u_{i-1},u_i)=(g({0}),g({1}))$).}
\end{cases}
\] 
The composition maps $c_i's$ are induced by the Moore composition of execution paths of $A$. The inclusion maps $I_i's$ are induced by the continuous maps $\P^\ell\overline{f}:\P^\ell_{g({0}),g({1})} \mathbb{M}^{\mathcal{G}}(A) \to T(\ell)$ for $\ell>0$. We obtain a well-defined diagram of topological spaces~\footnote{This point is left as an exercice; Verifying the commutativity relations is easy.} \[\mathcal{E}^f:\mathcal{P}^{g({0}),g({1})}(A^0)\to \top\] and, by Proposition~\ref{underlying-set-tensor}, there is an objectwise continuous surjective map
\[
k:\mathcal{E}^f \longrightarrow \D^f(1). 
\]
We deduce that $\liminj k$ is surjective. We want to prove that the composite map 
\[
\xymatrix@C=4em
{\widehat{\mathbb{C}}:\liminj \mathcal{E}^f  \ar@{->>}[r]^-{\liminj k}& (\liminj \D^f)(1) \fr{\P^1\psi}& \P^1\mathbb{M}^\mathcal{G}(X)=\P^\mathcal{G}X}
\]
is a continuous bijection. We already know that the map $\P^1\psi$ is surjective, and therefore that the map $\widehat{\mathbb{C}}:\liminj \mathcal{E}^f \to \P^\mathcal{G}X$ is surjective as well. To prove that $\widehat{\mathbb{C}}:\liminj \mathcal{E}^f \to \P^\mathcal{G}X$ is one-to-one, we must first introduce the notion of \textit{simplified} element. Let $\underline{x}$ be an element of some vertex of the diagram of spaces $\mathcal{E}^f $. We say that $\underline{x} \in \mathcal{E}^f (\underline{n})$ is \textit{simplified} if~\footnote{$d$ is the degree function of the Reedy category, see Appendix~\ref{reedy}.}
\[
d(\underline{n}) = \min \big\{d(\underline{m}) \mid \exists \underline{m}\in \Obj(\mathcal{P}^{g({0}),g({1})}(A^0)) \hbox{ and }\exists \underline{y}\in \mathcal{E}^f (\underline{m}),\underline{y}=\underline{x} \in \liminj \mathcal{E}^f 
\big\}.
\]
Let $\underline{x}$ be a simplified element belonging to some vertex $\mathcal{E}^f (\underline{n})$ of the diagram $\mathcal{E}^f $ with \[\underline{n} = ((u_0,\epsilon_1,u_1),(u_1,\epsilon_2,u_2),\dots ,(u_{n-1},\epsilon_n,u_n)).\]

$\bullet$ \underline{\textit{Case 1}}. It is impossible to have $\epsilon_i = \epsilon_{i+1}=0$ for some $1\leq i<n$. Indeed, otherwise $\underline{x}$ would be of the form \[(\dots,\gamma_i\mu_{\ell_i},\gamma_{i+1}\mu_{\ell_{i+1}},\dots)\] where $\gamma_i$ and $\gamma_{i+1}$ would be two execution paths of $A$. Using the equality \[c_i\big((\dots,\gamma_i\mu_{\ell_i},\gamma_{i+1}\mu_{\ell_{i+1}},\dots)\big) = (\dots,\gamma_i\mu_{\ell_i}*\gamma_{i+1}\mu_{\ell_{i+1}},\dots),\] the tuple $\underline{x}$ can then be identified in the colimit with the tuple \[\bigg(\dots,\underbrace{\big(\gamma_i\mu_{\ell_i}*\gamma_{i+1}\mu_{\ell_{i+1}}\big)\mu^{-1}_{\ell_i+\ell_{i+1}}}_{\in \P^{\mathcal{G}}A \hbox{ by Proposition~\ref{variable-length}}}\mu^{}_{\ell_i+\ell_{i+1}},\dots\bigg) \in \mathcal{E}^f (\underline{n'})
	\] 
with \[d(\underline{n'}) = n-1 + \sum_i \epsilon_i < d(\underline{n}).\] It is a contradiction because $\underline{x}$ is simplified by hypothesis.

$\bullet$ \underline{\textit{Case 2}}. Suppose that $\epsilon_i=1$ for some $1\leq i\leq n$ and that $\underline{x}$ is of the form 
\[(\dots,\P\overline{g}(\delta_{z_{i}}\phi_{i}\mu_{\ell_{i}}),\dots).\] If $z_{i}\in \mathbf{S}^{n-1}$, then using the equality 
\[
	I_{i}\big((\dots,g\delta_{z_{i}}\phi_{i}\mu_{\ell_{i}},\dots)\big) = (\dots,\P\overline{g}(\delta_{z_{i}}\phi_{i}\mu_{\ell_{i}}),\dots), 
\]
the tuple $\underline{x}$ can then be identified in the colimit with the tuple \[(\dots,g\delta_{z_{i}}\phi_{i}\mu_{\ell_{i}},\dots) \in \mathcal{E}^f(\underline{n'})\] 
with \[d(\underline{n'}) = n+ \big(\sum_i \epsilon_i\big) - 1< d(\underline{n}) .\] It is a contradiction because $\underline{x}$ is simplified by hypothesis. We deduce that in this case, $z_{i}\in \mathbf{D}^n\backslash\mathbf{S}^{n-1}$.

$\bullet$ \underline{\textit{Partial conclusion}}. Consequently, for all simplified elements $\underline{x}=(x_1,\dots,x_n)$ of $\mathcal{E}^f$, we have
\[
\widehat{\mathbb{C}}(\underline{x}) = (f_1 x_1)  * \dots * (f_n x_n)
\]
with for all $1\leq i \leq n$, 
\[
\begin{cases}
f_i=f &\hbox{ and } x_i \in \P^{\ell_i}A\\
f_i=\P\psi&\hbox{ and } x_i=\P\overline{g}(\delta_{z_{i}}\phi_{i}\mu_{\ell_{i}}) \hbox{ with }z_i\in \mathbf{D}^n\backslash\mathbf{S}^{n-1}
\end{cases}
\]
and there are no two consecutive terms of the first form (i.e. $f_i=f_{i+1}=f$ for some $i$). It means that it is the finite Moore composition of $\widehat{\mathbb{C}}(\underline{x})$ of Theorem~\ref{calcul_final_structure}. 

$\bullet$ \underline{\textit{Injectivity of $\widehat{\mathbb{C}}$}}. Let $\underline{x}$ and $\underline{y}$ be two elements of $\liminj \mathcal{E}^f$ such that $\widehat{\mathbb{C}}(\underline{x}) = \widehat{\mathbb{C}}(\underline{y})$. We can suppose that both $\underline{x}$ and $\underline{y}$ are simplified. Let $\underline{x}=(x_1,\dots,x_m)$ and $\underline{y}=(y_1,\dots,y_n)$. Then \[(f_1x_1)* \dots * (f_mx_m) = (g_1 y_1)* \dots * (g_n y_n).\] Since both members of the equality are the finite Moore composition of Theorem~\ref{calcul_final_structure}, we deduce that $m=n$ and that for all $1\leq i \leq m$, we have $f_i x_i = g_i y_i$. For a given $i\in [1,m]$, there are two mutually exclusive possibilities: 
\begin{enumerate}[leftmargin=*]
	\item $f_i=g_i=f$ and $x_i$ and $y_i$ are two paths of length $\ell_i$ of $A$. Since $f$ is one-to-one because $\mathbf{S}^{n-1}$ is a subset of $\mathbf{D}^n$, we deduce that $x_i=y_i$.
	\item $f_i=g_i=\P\psi$, $x_i=\overline{g}\delta_{z_{i}}\phi_i\mu_{\ell_i}$ and $y_i=\overline{g}\delta_{t_{i}}\psi_i\mu_{\ell_i}$, with $z_i,t_i\in \mathbf{D}^n\backslash\mathbf{S}^{n-1}$ and $\phi_i,\psi_i\in \mathcal{G}(1,1)$. We also have $\P\psi (x_i) = \widehat{g} \delta_{z_{i}}\phi_i\mu_{\ell_i}$ and $\P\psi (y_i) = \widehat{g} \delta_{t_{i}}\psi_i\mu_{\ell_i}$. The restriction of $\widehat{g}$ to $\globG(\mathbf{D}^n)\backslash\globG(\mathbf{S}^{n-1})$ being one-to-one, we deduce that $z_i=t_i$, $\phi_i=\psi_i$ and therefore once again that $x_i =  y_i$.
\end{enumerate}
We conclude that $\underline{x} = \underline{y}$ and that the map $\widehat{\mathbb{C}}:\liminj \mathcal{E}^f \to \P^{\mathcal{G}}X$ is one-to-one. 

$\bullet$ \underline{\textit{Informal summary}}. The arrows of the small category $\mathcal{P}^{g({0}),g({1})}(A^0)$ and the relations satisfied by them prove that each element of the colimit $\liminj \mathcal{E}^f$ has a simplified rewriting and this simplified rewriting coincides with the normal form of Theorem~\ref{calcul_final_structure}. The latter theorem relies on the fact that all execution paths of globes are one-to-one, and more generally that all execution paths of cellular multipointed $d$-spaces are locally injective.

$\bullet$ \underline{\textit{Injectivity of $\P^1\psi$}}. At this point of the proof, we have a composite continuous map 
\[
\xymatrix@C=4em
{
	\liminj \mathcal{E}^f  \ar@{->}@/_40pt/[rr]_-{\hbox{\tiny continuous bijection}}^-{\widehat{\mathbb{C}}}\ar@{->>}[r]^-{\liminj k} & \liminj \mathcal{D}^f(1) \ar@{->>}[r]^-{\P^1\psi} & \P^{\mathcal{G}}X.
}
\]
Let $a,b\in \liminj \mathcal{D}^f(1)$ such that $\P^1\psi(a)=\P^1\psi(b)$. Let $\overline{a},\overline{b}\in \liminj \mathcal{E}^f$ such that $(\liminj k)(\overline{a})=a$ and $(\liminj k)(\overline{b})=b$. Then $\overline{a}=\overline{b}$ and therefore $a=b$. We have proved that the continuous map $\P^1\psi:\liminj \mathcal{D}^f(1) \to \P^{\mathcal{G}}X$ is one-to-one. 
\epf

\bth \label{pre-calculation-pathspace2}
Under the hypotheses and the notations of this section. The map of $\mathcal{G}$-spaces \[\P\psi:\liminj \D^f \longrightarrow \P\mathbb{M}^\mathcal{G}(X)\] is an isomorphism of $\mathcal{G}$-spaces.
\eth

\bpf 
We already know by Theorem~\ref{pre-calculation-pathspace} that the map of $\mathcal{G}$-spaces \[\P\psi:\liminj \D^f \longrightarrow \P\mathbb{M}^\mathcal{G}(X)\] is an objectwise continuous bijection. We want to prove that it is an objectwise homeomorphism. Since all objects of the reparametrization category $\mathcal{G}$ are isomorphic, it suffices to prove that \[\P^1\psi:\liminj \D^f(1) \longrightarrow \P^\mathcal{G}X\] is a homeomorphism. Consider a set map $\xi:[0,1]\to \liminj \D^f(1)$ such that the composite map 
\[
\overline{\xi}:[0,1]\stackrel{\xi}\longrightarrow  \liminj \D^f(1) \stackrel{\P^1\psi}\longrightarrow \P^\mathcal{G}X
\]
is continuous. By Corollary~\ref{DeltaHomeo2}, it suffices to prove that the set map \[\xi:[0,1]\longrightarrow \liminj \D^f(1)\] is continuous as well.

$\bullet$ \underline{\textit{First reduction}}. The composite continuous map $\overline{\xi}$ gives rise by adjunction to a continuous map 
\[
\widehat{\xi}:[0,1]\p [0,1] \longrightarrow |X|.
\]
Since $[0,1]\p [0,1]$ is compact and since $|X|$ is weakly Hausdorff by Proposition~\ref{p1}, the subset $\widehat{\xi}([0,1]\p [0,1])$ is a compact subset of $|X|$. By Proposition~\ref{cpt-intersect-finite}, $\widehat{\xi}([0,1]\p [0,1])$ intersects a finite number of globular cells of the cellular multipointed $d$-space $X$. Therefore we can suppose that the multipointed $d$-space $X$ is finite by Proposition~\ref{pushout-qcof-plus}. Write \[\{c_j\mid j\in J\}\] for its finite set of globular cells.

$\bullet$ \underline{\textit{Second reduction}}. It suffices to prove that there exists a finite covering $\{F_1,\dots,F_n\}$ of $[0,1]$ by closed subsets such that each restriction $\overline{\xi}\rest_{F_i}$ factors through the colimit $\liminj \D^f(1)$. Let $\mathcal{T}$ be the set defined as follows: 
\[
\mathcal{T} = \bigg\{\carrier\big(\overline{\xi}(u)\big) \mid u\in [0,1]\bigg\}.
\]
Suppose that $\mathcal{T}$ is infinite. Since $J$ is finite, there exist $j_0\in J$ and a sequence $(\overline{\xi}(u_n))_{n\geq 0}$ of execution paths of $X$ such that the numbers $i_n$ of times that ${c_{j_0}}$ appears in the carrier of $\overline{\xi}(u_n)$ for $n\geq 0$ give rise to a strictly increasing sequence of integers $(i_n)_{n\geq 0}$. Since $[0,1]$ is sequentially compact, the sequence $(u_n)_{n\geq 0}$ has a convergent subsequence. By continuity, the sequence $(\overline{\xi}(u_n))_{n\geq 0}$ has therefore a convergent subsequence in $\P^{\mathcal{G}}X$. This contradicts Theorem~\ref{bounded}. Consequently, the set $\mathcal{T}$ is finite. For each carrier $\underline{c}\in \mathcal{T}$, let 
\[
U_{\underline{c}} = \bigg\{u\in [0,1]\mid\carrier(\overline{\xi}(u)) = \underline{c}\bigg\}.
\]
Consider the closure $\widehat{U_{\underline{c}}}$ of $U_{\underline{c}}$ in $[0,1]$. We obtain the finite covering of $[0,1]$ by closed subsets 
\[
[0,1] =\bigcup_{\underline{c}\in \mathcal{T}} \widehat{U_{\underline{c}}}.
\]
We replace $[0,1]$ by  $\widehat{U_{\underline{c}}}$ which is compact, metrizable and therefore sequential. The carrier \[\underline{c}=[c_{j_1},\dots,c_{j_m}]\] is fixed until the very end of the proof.

$\bullet$ \underline{\textit{Third reduction}}. The attaching maps 
\[
\widehat{g_{j_k}}:\globG(\mathbf{D}^{n_{j_k}}) \longrightarrow X
\]
for $1\leq k \leq m$ of the cells $c_{j_1},\dots,c_{j_m}$ yield a map of multipointed $d$-spaces 
\[
\widehat{g_{\underline{c}}}:\globG(\mathbf{D}^{n_{j_1}})* \dots * \globG(\mathbf{D}^{n_{j_m}})\longrightarrow X.
\] 
Let $\alpha_{i-1}$ ($\alpha_{i}$ resp.) be the initial state (the final state resp.) of $\globG(\mathbf{D}^{n_{j_i}})$ for $1\leq i\leq m$ in $\globG(\mathbf{D}^{n_{j_1}})* \dots * \globG(\mathbf{D}^{n_{j_m}})$. We obtain a map of $\mathcal{G}$-spaces
\[
\mathbb{F}^{\mathcal{G}^{op}}_1(\mathbf{D}^{n_{j_1}})\ot \dots \ot \mathbb{F}^{\mathcal{G}^{op}}_1(\mathbf{D}^{n_{j_m}}) \longrightarrow \D^f(\underline{m})
\]
for some $\underline{m}$ belonging to $\mathcal{P}^{g({0}),g({1})}(A^0)$ such that 
\[
\D^f(\underline{m}) = Z_{\widehat{g_{\underline{c}}}(\alpha_0),\widehat{g_{\underline{c}}}(\alpha_1)} \ot \dots \ot Z_{\widehat{g_{\underline{c}}}(\alpha_{m-1}),\widehat{g_{\underline{c}}}(\alpha_m)}.
\]
Using Proposition~\ref{Ftenseur}, we obtain a continuous map 
\[
y_{\underline{c}}:\mathcal{G}(1,m) \p \mathbf{D}^{n_{j_1}} \p \dots \p \mathbf{D}^{n_{j_m}} \longrightarrow Z_{\underline{c}}\subset\D^f(\underline{m})(1)
\]
where $Z_{\underline{c}}$ is, by definition, the image of $y_{\underline{c}}$. At this point, we have obtained that the continuous map \[\overline{\xi}\rest_{{U_{\underline{c}}}}:{U_{\underline{c}}}\longrightarrow \P^{\mathcal{G}}X\] factors as a composite of maps 
\[
\overline{\xi}\rest_{{U_{\underline{c}}}}:U_{\underline{c}} \longrightarrow Z_{\underline{c}}\subset\D^f(\underline{m})(1) \stackrel{p_{\underline{m}}}\longrightarrow \liminj \D^f(1) \longrightarrow \P^{\mathcal{G}}X.
\]
Consider a sequence $(u_n)_{n\geq 0}$ of $U_{\underline{c}}$ converging to $u_\infty\in \widehat{U_{\underline{c}}}$. Then for each $n\geq 0$, $\overline{\xi}(u_n)$ belongs to the image of $\P^{\mathcal{G}}\widehat{g_{\underline{c}}}$ which is a closed subset of the sequential space $\P^{\mathcal{G}}X$ by Theorem~\ref{img-closed}. Thus $\overline{\xi}(u_\infty)$ belongs to the image of $\P^{\mathcal{G}}\widehat{g_{\underline{c}}}$ as well. Since \[\P^{\mathcal{G}}\widehat{g_{\underline{c}}} = \P^1\psi.p_{\underline{m}}.y_{\underline{c}},\] we have obtained that the continuous map \[\overline{\xi}\rest_{\widehat{U_{\underline{c}}}}:\widehat{U_{\underline{c}}}\longrightarrow \P^{\mathcal{G}}X\] factors as a composite of maps 
\[
\overline{\xi}\rest_{\widehat{U_{\underline{c}}}}:\widehat{U_{\underline{c}}} \longrightarrow Z_{\underline{c}}\subset\D^f(\underline{m})(1) \stackrel{p_{\underline{m}}}\longrightarrow \liminj \D^f(1) \stackrel{\P^1\psi}\longrightarrow \P^{\mathcal{G}}X.
\]
They are all of them continuous except maybe the left-hand one from $\widehat{U_{\underline{c}}}$ to $Z_{\underline{c}}$ (cf. the remark after this proof). Since $\widehat{U_{\underline{c}}}$ is sequential, it remains to prove that the map 
\[
{\xi}\rest_{\widehat{U_{\underline{c}}}}:\widehat{U_{\underline{c}}} \longrightarrow Z_{\underline{c}}\subset\D^f(\underline{m})(1) \stackrel{p_{\underline{m}}}\longrightarrow \liminj \D^f(1)
\]
is sequentially continuous to complete the proof. 

$\bullet$ \underline{\textit{Sequential continuity}}. Consider a sequence $(u_n)_{n\geq 0}$ of $\widehat{U_{\underline{c}}}$ which converges to $u_\infty\in \widehat{U_{\underline{c}}}$.  Write \[{\xi}(u_n) = p_{\underline{m}}\big(y_{\underline{c}}(\phi_n,z^1_n,\dots,z^m_n)\big)\] for all $n\geq 0$. We obtain \[\overline{\xi}(u_n) = (\P^{\mathcal{G}}\widehat{g_{\underline{c}}})(\phi_n,z^1_n,\dots,z^m_n)\] for all $n\geq 0$. By Theorem~\ref{img-closed}, the sequence $((\phi_n,z^1_n,\dots,z^m_n))_{n\geq 0}$ has a limit point $(\phi_\infty,z^1_\infty,\dots,z^m_\infty)$. We deduce that the sequence $({\xi}(u_n))_{n\geq 0}$ has a limit point because both $y_{\underline{c}}$ and $p_{\underline{m}}$ are continuous. It is necessarily equal to ${\xi}(u_\infty)$ because the map \[\P^1\psi:\liminj \D^f(1) \to \P^\mathcal{G}X\] is continuous bijective by Theorem~\ref{pre-calculation-pathspace} and because \[\overline{\xi}=\P^1\psi.\xi.\] The same argument shows that every subsequence of $({\xi}(u_n))_{n\geq 0}$ has a limit point which is necessarily ${\xi}(u_\infty)$. Suppose that the sequence $({\xi}(u_n))_{n\geq 0}$ does not converge to ${\xi}(u_\infty)$. Then there exists an open neighborhood $V$ of ${\xi}(u_\infty)$ such that for all $n\geq 0$, ${\xi}(u_n)\notin V$. Since $V^c$ is closed in $\liminj \D^f(1)$,  it means that ${\xi}(u_\infty)$ cannot be a limit point of the sequence $({\xi}(u_n))_{n\geq 0}$. Contradiction. It implies that the sequence $({\xi}(u_n))_{n\geq 0}$ converges to ${\xi}(u_\infty)$.
\epf

Before expounding the consequences of Theorem~\ref{pre-calculation-pathspace} and of Theorem~\ref{pre-calculation-pathspace2}, let us add an additional remark about the proof of Theorem~\ref{pre-calculation-pathspace2}. It can be proved that the inverse image $p_{\underline{m}}^{-1}(\gamma)$ for each $\gamma\in \liminj \D^f(1)$ is always finite. When the multipointed $d$-space $X$ does not contain any loop, i.e. when $\P_{\alpha,\alpha}^{\mathcal{G}}X$ is empty for all $\alpha\in X^0$, the map $p_{\underline{m}}$ is even one-to-one and it is then possible to prove that the set map $\widehat{U_{\underline{c}}}\to Z_{\underline{c}}$ is always continuous. On the contrary, when $X$ contains loops, the set map $\widehat{U_{\underline{c}}}\to Z_{\underline{c}}$ is not necessarily continuous mainly because an inverse image $p_{\underline{m}}^{-1}(\gamma)$ may contain several points.

\begin{cor}  \label{calculation-pathspace}
Suppose that $A$ is a cellular multipointed $d$-space. Consider a pushout diagram of multipointed $d$-spaces 
\[
\xymatrix@C=4em@R=4em
{
	\globG(\mathbf{S}^{n-1}) \fd{} \fr{} & A \ar@{->}[d]^-{} \\
	\globG(\mathbf{D}^{n}) \fr{} & \cocartesien X
}
\]
with $n\geq 0$. Then there is the pushout diagram of Moore flows 
\[
\xymatrix@C=4em@R=4em
{
	\moore^{\mathcal{G}}(\globG(\mathbf{S}^{n-1}))=\globP(\mathbb{F}^{\mathcal{G}^{op}}_1\mathbf{S}^{n-1}) \fd{} \fr{} & \moore^{\mathcal{G}}(A) \ar@{->}[d]^-{} \\
	\moore^{\mathcal{G}}(\globG(\mathbf{D}^{n}))=\globP(\mathbb{F}^{\mathcal{G}^{op}}_1\mathbf{D}^{n}) \fr{} & \cocartesien \moore^{\mathcal{G}}(X).
}
\]
\end{cor}

\begin{cor} \label{pre-c2}
	Let $X$ be a q-cofibrant multipointed $d$-space. Then $\moore^{\mathcal{G}}(X)$ is a q-cofibrant Moore flow.
\end{cor}

\bpf 
For every q-cofibrant Moore flow $X$, the canonical map $\varnothing \to X$ is a retract of a transfinite composition of the q-cofibrations $C:\varnothing\to \{0\}$, $R:\{0,1\}\to \{0\}$ and of the q-cofibrations $\globP(\mathbb{F}^{\mathcal{G}^{op}}_\ell\mathbf{S}^{n-1})\subset \globP(\mathbb{F}^{\mathcal{G}^{op}}_\ell\mathbf{D}^{n})$ for $\ell>0$ and $n\geq 0$. The cofibration $R:\{0,1\}\to \{0\}$ is not necessary to reach all q-cofibrant objects. Therefore, this theorem is a consequence of Theorem~\ref{cof-accessible} and of Corollary~\ref{calculation-pathspace}.
\epf

\bth \label{c2}
Let $X$ be a q-cofibrant Moore flow. Then the unit of the adjunction $
X \to \moore^{\mathcal{G}}(\lmoore^{\mathcal{G}}(X))$
is an isomorphism. 
\eth

\bpf
By Proposition~\ref{calculM}, the theorem holds when $X$ is a globe. It also clearly holds for $X=\{0\}$. For every q-cofibrant Moore flow $X$, the canonical map $\varnothing \to X$ is a retract of a transfinite composition of the q-cofibrations $C:\varnothing\to \{0\}$, $R:\{0,1\}\to \{0\}$ and of the q-cofibrations $\globP(\mathbb{F}^{\mathcal{G}^{op}}_\ell\mathbf{S}^{n-1})\subset \globP(\mathbb{F}^{\mathcal{G}^{op}}_\ell\mathbf{D}^{n})$ for $\ell>0$ and $n\geq 0$. The cofibration $R:\{0,1\}\to \{0\}$ is not necessary to reach all q-cofibrant objects. Therefore, this theorem is also a consequence of Theorem~\ref{cof-accessible} and of Corollary~\ref{calculation-pathspace}.
\epf

\begin{cor} \label{general-fact}
	Let $X$ be a q-cofibrant Moore flow. Then there is the homeomorphism 
	\[
	\P^1 X \iso \P^{\mathcal{G}}(\lmoore^{\mathcal{G}}(X)).
	\]
\end{cor}

\bpf Apply the functor $\P^1(-)$ to the isomorphism $X\iso \moore^{\mathcal{G}}(\lmoore^{\mathcal{G}}(X))$.
\epf

\bth \label{counit-iso}
Let $\lambda$ be a limit ordinal. Let \[X:\lambda \longrightarrow \ptop{\mathcal{G}}\] be a colimit preserving functor such that
\begin{itemize}[leftmargin=*]
	\item The multipointed $d$-space $X$ is a set, in other terms $X=(X^0,X^0,\varnothing)$.
	\item For all $\nu<\lambda$, there is a pushout diagram of multipointed $d$-spaces 
	\[
	\xymatrix@C=4em@R=4em
	{
		\globG(\mathbf{S}^{n_\nu-1}) \fd{} \fr{g_\nu} & X_\nu \ar@{->}[d]^-{} \\
		\globG(\mathbf{D}^{n_\nu}) \fr{\widehat{g_\nu}} & \cocartesien X_{\nu+1}
	}
	\]
	with $n_\nu \geq 0$. 
\end{itemize}
Let $X_\lambda = \liminj_{\nu<\lambda} X_\nu$. For all $\nu \leq \lambda$, the counit map 
\[
\kappa_\nu:\lmoore^{\mathcal{G}} (\moore^{\mathcal{G}}(X_\nu)) \longrightarrow X_\nu
\]
is an isomorphism.
\eth

\bpf
The map $\kappa_0$ is an isomorphism because $X_0$ is a set. By Theorem~\ref{cof-accessible}, and since $\lmoore^{\mathcal{G}}$ is a left adjoint, it suffices to prove that if $\kappa_\nu$ is an isomorphism, then $\kappa_{\nu+1}$ is an isomorphism. Assume that $\kappa_\nu$ is an isomorphism. By Corollary~\ref{calculation-pathspace}, there is the pushout diagram of Moore flows 
\[
\xymatrix@C=4em@R=4em
{
	\moore^{\mathcal{G}}(\globG(\mathbf{S}^{n_\nu-1})) = \globP(\mathbb{F}^{\mathcal{G}^{op}}_1\mathbf{S}^{n_\nu-1}) \fd{} \fr{g_\nu} & \moore^{\mathcal{G}}(X_\nu) \ar@{->}[d]^-{} \\
	\moore^{\mathcal{G}}(\globG(\mathbf{D}^{n_\nu})) = \globP(\mathbb{F}^{\mathcal{G}^{op}}_1\mathbf{D}^{n_\nu}) \fr{\widehat{g_\nu}} & \cocartesien \moore^{\mathcal{G}}(X_{\nu+1}).
}
\]
Apply again the left adjoint $\lmoore^{\mathcal{G}}$ to this diagram, we obtain by using the induction hypothesis that $\kappa_{\nu+1}$ is an isomorphism. 
\epf

\begin{cor} \label{c1}
	For every q-cofibrant multipointed $d$-space $X$, the counit of the adjunction $\lmoore^{\mathcal{G}}(\moore^{\mathcal{G}}(X))\to X$ is an isomorphism of multipointed $d$-spaces. 
\end{cor}

\bpf
It is due to the fact that every q-cofibrant multipointed $d$-space $X$ is a retract of a cellular multipointed $d$-space (note that the cofibration $R:\{0,1\}\to \{0\}$ is not required to reach all cellular multipointed $d$-spaces) and that a retract of an isomorphism is an isomorphism.
\epf

\section{From multipointed \mins{d}-spaces to flows}
\label{conclusion}

The goals of this final section are to complete the proof of the Quillen equivalence between multipointed $d$-spaces and Moore flows in Theorem~\ref{adj-zigzag}, which together with the results of \cite{Moore1} establish that multipointed $d$-spaces and flows have Quillen equivalent q-model structures, and to give a new and conceptual proof of \cite[Theorem~7.5]{mdtop} in Theorem~\ref{lepb4} independent of \cite{model2}. We also give a new presentation of the underlying homotopy type of  flow in Proposition~\ref{underlyinghomotopytype}.

\bth \label{adj-zigzag}
	The adjunction $\lmoore^{\mathcal{G}}\dashv \moore^{\mathcal{G}}:\dtopG \leftrightarrows \ptop{\mathcal{G}}$ induces a Quillen equivalence between the q-model structure of Moore flows and the q-model structure of multipointed $d$-spaces.
\eth

\bpf Since the q-fibrations of Moore flows are the maps of Moore flows inducing an objectwise q-fibration on the $\mathcal{G}$-spaces of execution paths, the functor $\moore^{\mathcal{G}}$ takes q-fibrations of multipointed $d$-spaces to q-fibrations of Moore flows. Since $\moore^{\mathcal{G}}$ preserves the set of states and since trivial q-fibrations of Moore flows are maps inducing a bijection on states and an an objectwise trivial q-fibration on the $\mathcal{G}$-spaces of execution paths, the functor $\moore^{\mathcal{G}}$ takes trivial q-fibrations of multipointed $d$-spaces to trivial q-fibrations of Moore flows. Therefore, the functor $\moore^{\mathcal{G}}:\ptop{\mathcal{G}} \to \dtopG$ is a right Quillen adjoint. 

By Theorem~\ref{c2}, the map $X \to \moore^{\mathcal{G}}(\lmoore^{\mathcal{G}}(X))$ is a weak equivalence of Moore flows for every q-cofibrant Moore flow $X$. Let $Y$ be a (q-fibrant) multipointed $d$-space. Then the composite map of multipointed $d$-spaces
\[
\lmoore^{\mathcal{G}}(\moore^{\mathcal{G}}(Y^{cof})) \stackrel{\iso}\longrightarrow Y^{cof} \stackrel{\simeq} \longrightarrow Y
\]
where $Y^{cof}$ is a q-cofibrant replacement of $Y$, is a weak equivalence of multipointed $d$-spaces because: 1) the left-hand map is an isomorphism by Corollary~\ref{c1}; 2) the right-hand map is a weak equivalence by definition of a cofibrant replacement. 
\epf

Let us give now some reminders about flows and the categorization functor $cat$ from multipointed $d$-spaces to flows.

\bd \cite[Definition~4.11]{model3} \label{def-flow}
A {\rm flow} is a small semicategory enriched over the closed monoidal category $(\top,\p)$. The corresponding category is denoted by $\dtop$. 
\ed

Let us expand the definition above. A \textit{flow} $X$ consists of a topological space $\P X$ of execution paths, a discrete space $X^0$ of states, two continuous maps $s$ and $t$ from $\P X$ to $X^0$ called the source and target map respectively, and a continuous and associative map \[*:\{(x,y)\in \P X\p \P X; t(x)=s(y)\}\longrightarrow \P X\] such that $s(x*y)=s(x)$ and $t(x*y)=t(y)$.  A morphism of flows $f:X\longrightarrow Y$ consists of a set map $f^0:X^0\longrightarrow Y^0$ together with a continuous map $\P f:\P X\longrightarrow \P Y$ such that 
\[\begin{aligned}
&f^0(s(x))=s(\P f(x)),\\
&f^0(t(x))=t(\P f(x)),\\
&\P f(x*y)=\P f(x)*\P f(y).
\end{aligned}\] Let \[\P_{\alpha,\beta}X = \{x\in \P X\mid s(x)=\alpha \hbox{ and } t(x)=\beta\}.\]

\begin{nota}
	The map $\P f:\P X\longrightarrow \P Y$ can be denoted by $f:\P X\to \P Y$ is there is no ambiguity. The set map $f^0:X^0\longrightarrow Y^0$ can be denoted by $f:X^0\longrightarrow Y^0$ is there is no ambiguity.
\end{nota}

The category $\dtop$ is locally presentable. Every set can be viewed as a flow with an empty path space. The obvious functor $\set \subset \dtop$ is limit-preserving and colimit-preserving. One another example of flow is important for the sequel: 

\begin{exa} \label{ex-glob}
	For a topological space $Z$, let $\glob(Z)$ be the flow defined by 
	\[
	\begin{aligned}
	&\glob(Z)^0=\{0,1\},\\
	&\P \glob(Z)= \P_{0,1} \glob(Z)=Z,\\
	&s=0,\  t=1.
	\end{aligned}
	\]
	This flow has no composition law.
\end{exa}

\begin{nota}
	\[\begin{aligned}
	&C:\varnothing \to \{0\},\\
	&R:\{0,1\} \to \{0\},\\
	&\vI = \glob(\{0\}).
	\end{aligned}\]
\end{nota}

The \textit{q-model structure} of flows is the unique combinatorial model structure such that \[\{\glob(\mathbf{S}^{n-1})\subset \glob(\mathbf{D}^{n}) \mid n\geq 0\} \cup\{C,R\}\] is the set of generating cofibrations and such that \[\{\glob(\mathbf{D}^{n}\p\{0\})\subset \glob(\mathbf{D}^{n+1}) \mid n\geq 0\}\] is the set of generating trivial cofibrations (e.g. \cite[Theorem~7.6]{QHMmodel}) where the maps $\mathbf{D}^{n}\subset \mathbf{D}^{n+1}$ are induced by the mappings $(x_1,\dots,x_n) \mapsto (x_1,\dots,x_n,0)$. The weak equivalences are the maps of flows $f:X\to Y$  inducing a bijection $f^0:X^0\iso Y^0$ and a weak homotopy equivalence $\P f:\P X \to \P Y$ and the fibrations are the maps of flows $f:X\to Y$  inducing a q-fibration $\P f:\P X \to \P Y$ of topological spaces.

Let $X$ be a multipointed $d$-space. Consider for every $(\alpha,\beta)\in X^0 \p X^0$ the coequalizer of spaces \[\P_{\alpha,\beta}X = \liminj\left( \P^{\mathcal{G}}_{\alpha,\beta}X\p \mathcal{G}(1,1) \rightrightarrows \P^{\mathcal{G}}_{\alpha,\beta}X\right)\] where the two maps are $(c,\phi)\mapsto c$ and $(c,\phi)\mapsto c.\phi$. Let $[-]_{\alpha,\beta}:\P^{\mathcal{G}}_{\alpha,\beta}X \rightarrow \P_{\alpha,\beta}X$ be the canonical map.  

\bth \cite[Theorem~7.2]{mdtop} \label{def-cat} 
Let $X$ be a multipointed $d$-space. Then there exists a flow $cat(X)$ with $cat(X)^0=X^0$, $\P_{\alpha,\beta}cat(X)= \P_{\alpha,\beta}X$ and the composition law $*:\P_{\alpha,\beta}X \p \P_{\beta,\gamma}X \rightarrow \P_{\alpha,\gamma}X$ is for every triple $(\alpha,\beta,\gamma)\in X^0\p X^0\p X^0$ the unique map making the following diagram commutative:
\[
\xymatrix@C=4em@R=4em{
	\P_{\alpha,\beta}^{\mathcal{G}}X \p \P_{\beta,\gamma}^{\mathcal{G}}X
	\fr{*_{N}}\fd{[-]_{\alpha,\beta}\p [-]_{\beta,\gamma}} &
	\P_{\alpha,\gamma}^{\mathcal{G}}X \fd{[-]_{\alpha,\gamma}} \\
	\P_{\alpha,\beta}X \p \P_{\beta,\gamma}X \fr{} &
	\P_{\alpha,\gamma}X}
\] 
where $*_N$ is the normalized composition (cf. Definition~\ref{composition_map}).  The mapping $X \mapsto cat(X)$ induces a functor from $\ptop{\mathcal{G}}$ to $\dtop$.  \eth

\bd \label{cat-func} The functor $cat:\ptop{\mathcal{G}}\to \dtop$ is called the {\rm categorization functor}. \ed

The motivation for the constructions of this paper and of \cite{Moore1} comes from the following theorem which is added for completeness. 

\bth \label{noadjoint} The categorization functor $cat:\ptop{\mathcal{G}}\to\dtop$ is neither a left adjoint nor a right adjoint. In particular, it cannot be a left or a right Quillen equivalence. \eth

\bpf
This functor is not a left adjoint by \cite[Proposition~7.3]{mdtop}. Suppose that it is a right adjoint. Let $\mathcal{L}:\dtop\to \ptop{\mathcal{G}}$ be the left adjoint. Pick a nonempty topological space $Z$. The adjunction yields the bijection of sets \[\ptop{\mathcal{G}}(\mathcal{L}(\glob(Z)),\vI^{\mathcal{G}}) \iso \dtop(\glob(Z),\vI).\] Since $Z$ is nonempty, a map of flows from $\glob(Z)$ to $\vI$ is determined by the choice of a map from $Z$ to $\{0\}$. We deduce that there is exactly one map $f$  of multipointed $d$-spaces from $\mathcal{L}(\glob(Z))$ to $\vI^{\mathcal{G}}$. Suppose that $\mathcal{L}(\glob(Z))$ contains at least one execution path $\phi:[0,1]\to |\mathcal{L}(\glob(Z))|$. Then $f.\phi$ is an execution path of $\vI^{\mathcal{G}}$. Every map $g\in \ptop{\mathcal{G}}(\vI^{\mathcal{G}},\vI^{\mathcal{G}}) \iso \{[0,1]\iso^+ [0,1]\}$ gives rise to and execution path $g.f.\phi$ of $\vI^{\mathcal{G}}$. since $g.f\in \ptop{\mathcal{G}}(\mathcal{L}(\glob(Z)),\vI^{\mathcal{G}})$, we deduce that $g.f=f$. Contradiction. We deduce that the multipointed $d$-space $\mathcal{L}(\glob(Z))$ does not contain any execution path. Therefore this multipointed $d$-space is of the form $(U_Z,U_Z^0,\varnothing)$. We obtain the bijection  $\mtop((U_Z,U_Z^0),([0,1],\{0,1\}))\iso \{f\}$. Suppose that $U_Z$ is nonempty. Then for all $g\in \mtop(([0,1],\{0,1\}),([0,1],\{0,1\}))$, we have $g.f = f$. The only possibilities are that $f=0$ or $f=1$. Since $f$ is the unique element, we deduce that $U_Z=\varnothing$. There are also the natural bijections of sets  \[\ptop{\mathcal{G}}(\mathcal{L}(\{0\}),X) \iso \dtop(\{0\},cat(X)) \iso cat(X)^0 \iso X^0 \iso \ptop{\mathcal{G}}(\{0\},X).\] By the Yoneda lemma, we obtain $\mathcal{L}(\{0\})=\{0\}$. 

To summarize, if $\mathcal{L}:\dtop\to \ptop{\mathcal{G}}$ is a left adjoint to the functor $cat:\ptop{\mathcal{G}}\to \dtop$, then one has $\mathcal{L}(\{0\})=\{0\}$ and for all nonempty topological spaces $Z$, there is the equalities $\mathcal{L}(\glob(Z))=\varnothing$. By \cite[Theorem~6.1]{model3}, any flow is a colimit of globes and points. Since $\mathcal{L}$ is colimit-preserving, we deduce that for all flows $Y$, the multipointed $d$-space $\mathcal{L}(Y)$ is a set. We go back to the natural bijection given by the adjunction: \[\ptop{\mathcal{G}}(\mathcal{L}(Y),X) \iso \dtop(Y,cat(X)).\] Since $\mathcal{L}(Y)$ is a set, we obtain the natural bijection $\set(\mathcal{L}(Y),X^0) \iso \dtop(Y,cat(X))$. We obtain the natural bijection $\ptop{\mathcal{G}}(\mathcal{L}(Y),X^0) \iso \dtop(Y,cat(X))$ and by adjunction the natural bijection $\dtop(Y,X^0)\iso \dtop(Y,cat(X))$ since $cat(X^0)=X^0$. By Yoneda, we conclude that $cat(X)=X^0$ for all multipointed $d$-spaces $X$, which is a contradiction.
\epf

\bp \label{Ptenseur} \cite[Proposition~5.17]{Moore1}
Let $U$ and $U'$ be two topological spaces.  There is the natural isomorphism of $\mathcal{G}$-spaces
\[
\Delta_{\mathcal{G}^{op}} U \ot \Delta_{\mathcal{G}^{op}} U' \iso \Delta_{\mathcal{G}^{op}} (U \p U').
\]
\ep

Let $X$ be a flow. The Moore flow $\moore(X)$ is the enriched semicategory defined as follows: 
\begin{itemize}[leftmargin=*]
	\item The set of states is $X^0$.
	\item The $\mathcal{G}$-space $\P_{\alpha,\beta}\moore(X)$ is the $\mathcal{G}$-space $\Delta_{\mathcal{G}^{op}}(\P_{\alpha,\beta}X)$.
	\item The composition law is defined, using Proposition~\ref{Ptenseur} as the composite map 
	\[
	\xymatrix@C=4em
	{
		\Delta_{\mathcal{G}^{op}}(\P_{\alpha,\beta}X)\ot  \Delta_{\mathcal{G}^{op}}(\P_{\beta,\gamma}X) \iso \Delta_{\mathcal{G}^{op}}(\P_{\alpha,\beta}X \p \P_{\beta,\gamma}X)   \fr{\Delta_{\mathcal{G}^{op}}(*)}& \Delta_{\mathcal{G}^{op}}(\P_{\alpha,\gamma}) X.}
	\]
\end{itemize} 

The construction above yields a well-defined functor \[\moore:\dtop\to \dtopG.\] Consider a $\mathcal{G}$-flow $Y$. For all $\alpha,\beta\in Y^0$, let $Y_{\alpha,\beta}=\liminj \P_{\alpha,\beta}Y$. Let $(\alpha,\beta,\gamma)$ be a triple of states of $Y$. The composition law of the $\mathcal{G}$-flow $Y$ induces a continuous map \[Y_{\alpha,\beta} \p Y_{\beta,\gamma} \iso \liminj(\P_{\alpha,\beta}Y \ot \P_{\beta,\gamma}Y) \longrightarrow \liminj \P_{\alpha,\gamma}Y \iso Y_{\alpha,\gamma}\] which is associative. We obtain the 

\bp \label{def-N} \cite[Proposition~10.6 and Proposition~10.7]{Moore1}
For any $\mathcal{G}$-flow $Y$, the data 
\begin{itemize}[leftmargin=*]
	\item The set of states is $Y^0$
	\item For all $\alpha,\beta\in Y^0$, let $Y_{\alpha,\beta}=\liminj \P_{\alpha,\beta}Y$
	\item For all $\alpha,\beta,\gamma\in Y^0$, the composition law $Y_{\alpha,\beta}\p Y_{\beta,\gamma}\to Y_{\alpha,\gamma}$
\end{itemize}
assemble to a flow denoted by $\lmoore(Y)$. It yields a well-defined functor \[\lmoore:\dtopG \to \dtop.\]
There is an adjunction $\lmoore \dashv\moore$.
\ep

\bth \label{decomposing} There is the isomorphism of functors \[cat\iso\lmoore.\moore^{\mathcal{G}}.\] \eth

\bpf
First, let us notice that the functors $cat:\ptop{\mathcal{G}}\to \dtop$ (Theorem~\ref{def-cat}), $\moore^{\mathcal{G}}:\ptop{\mathcal{G}}\to \dtopG$ (Theorem~\ref{adj-multi-moore}) and $\lmoore:\dtopG\to \dtop$ (Proposition~\ref{def-N}) preserve the set of states by definition of these functors. Therefore, for every multipointed $d$-space $X$, the flows $cat(X)$ and $\lmoore.\moore^{\mathcal{G}}(X)$ have the same set of states $X^0$. Let $\mathcal{G}^1$ be the full subcategory of $\mathcal{G}$ generated by $1$: the set of objects is the singleton $\{1\}$ and $\mathcal{G}^1(1,1)=\mathcal{G}(1,1)$.  For $(\alpha,\beta)\in X^0\p X^0$, the inclusion functor $\iota:\mathcal{G}^1\subset \mathcal{G}$ induces a continuous map \[\liminj_{\mathcal{G}^1} \bigg(\left(\P_{\alpha,\beta}\moore^{\mathcal{G}}X\right).\iota\bigg) \to \liminj_{\mathcal{G}} \P_{\alpha,\beta}\moore^{\mathcal{G}}X.\] It turns out that there is the natural homeomorphisms \[\liminj_{\mathcal{G}^1} \bigg(\left(\P_{\alpha,\beta}\moore^{\mathcal{G}}X\right).\iota\bigg) \iso \liminj_{\mathcal{G}^1}\P_{\alpha,\beta}^1\moore^{\mathcal{G}}X \iso \liminj_{\mathcal{G}^1}\P_{\alpha,\beta}^{\mathcal{G}}X \iso \P_{\alpha,\beta}cat(X),\] 
the first one by definition of $\iota$, the second one by definition of $\moore^{\mathcal{G}}$ and the last one by definition of $cat$. We obtain a natural map of flows $cat(X)\to (\lmoore.\moore^{\mathcal{G}})(X)$ which is bijective on states. Let $\ell>0$ be an object of $\mathcal{G}$. Then the comma category $(\ell\!\downarrow\!\iota)$ is characterized as follows: 
\begin{itemize}[leftmargin=*]
	\item The set of objects is $\mathcal{G}(\ell,1)$ which is always nonempty for every $\ell>0$.
	\item The set of maps from an arrow $u:\ell\to 1$ to an arrow $v:\ell\to 1$ is an element of $\Mor(\mathcal{G})(u,v)$.
\end{itemize}
The comma category $(\ell\!\downarrow\!\iota)$ is connected since in any diagram of $\mathcal{G}$ of the form 
\[
\xymatrix@C=4em@R=4em
{[0,\ell] \ar@{=}[d]\fr{u} & [0,1] \ar@{-->}[d]^-{k}\\
	[0,\ell] \fr{v} & [0,1],}
\]
there exists a map $k\in \mathcal{G}([0,1],[0,1])$ making the square commute: take $k=v.u^{-1}$. By \cite[Theorem~1 p.~213]{MR1712872}, we deduce that the natural map of flows $cat(X)\to (\lmoore.\moore^{\mathcal{G}})(X)$ induces a homeomorphism between the spaces of paths. 
\epf

\begin{cor}  \label{calculation-pathspace3}
	Suppose that $A$ is a cellular multipointed $d$-space. Consider a pushout diagram of multipointed $d$-spaces 
	\[
	\xymatrix@C=4em@R=4em
	{
		\globG(\mathbf{S}^{n-1}) \fd{} \fr{} & A \ar@{->}[d]^-{} \\
		\globG(\mathbf{D}^{n}) \fr{} & \cocartesien X
	}
	\]
	with $n\geq 0$. Then there is the pushout diagram of flows 
	\[
	\xymatrix@C=4em@R=4em
	{
		\glob(\mathbf{S}^{n-1}) \fd{} \fr{} & cat(A) \ar@{->}[d]^-{} \\
		\glob(\mathbf{D}^{n}) \fr{} & \cocartesien cat(X).
	}
	\]
\end{cor}

\bpf It is a consequence of Corollary~\ref{calculation-pathspace}, Theorem~\ref{decomposing} and of the fact that $\lmoore:\dtopG\to\dtop$ is a left adjoint.
\epf

\bd  \label{inv-up-to-homotopy}
We consider the composite functors 
\[
\xymatrix@C=2.5em@R=2em
{
(\mathbf{L}cat):\ptop{\mathcal{G}}\fr{(-)^{cof}}& \ptop{\mathcal{G}} \fr{cat}& \dtop \\
(\mathbf{L}cat)^{-1}:\dtop \fr{\moore} & \dtopG \fr{(-)^{cof}}& \dtopG \fr{\lmoore^{\mathcal{G}}}& \ptop{\mathcal{G}}
}
\]
where $(-)^{cof}$ is a q-cofibrant replacement functor.
\ed

We can now write down the new proof of \cite[Theorem~7.5]{mdtop}.

\bth \label{lepb4}
The categorization functor from multipointed $d$-spaces to flows \[cat:\ptop{\mathcal{G}}\longrightarrow \dtop\] takes q-cofibrant multipointed $d$-spaces to q-cofibrant flows. Its total left derived functor in the sense of \cite{HomotopicalCategory} induces an equivalence of categories between the homotopy categories of the q-model structures.
\eth

\bpf The functor $cat\iso\lmoore.\moore^{\mathcal{G}}$ takes q-cofibrant multipointed $d$-spaces to q-cofibrant flows by Corollary~\ref{pre-c2} and because $\lmoore$ is a left Quillen adjoint. The rest of the proof is divided in four parts.

$\bullet$ \underline{\textit{$X\simeq Y\Rightarrow(\mathbf{L}cat)(X)\simeq(\mathbf{L}cat)(Y)$}}. 
Let $X\simeq Y$ be two weakly equivalent multipointed $d$-spaces in the q-model structure. Then there is the weak equivalence $X^{cof}\simeq Y^{cof}$. Since all multipointed $d$-spaces are q-fibrant, the right Quillen functor $\moore^{\mathcal{G}}$ takes weak equivalences of multipointed $d$-spaces to weak equivalences of Moore flows. Since $\lmoore$ is a left Quillen functor and since $\moore^{\mathcal{G}}$ preserves q-cofibrancy by Corollary~\ref{pre-c2}, we deduce using Theorem~\ref{decomposing} the sequence of isomorphisms and weak equivalences 
\[
(\mathbf{L}cat)(X) \iso \lmoore\moore^{\mathcal{G}}(X^{cof}) \simeq \lmoore\moore^{\mathcal{G}}(Y^{cof}) \iso (\mathbf{L}cat)(Y).
\] 

$\bullet$ \underline{\textit{$X\simeq Y\Rightarrow(\mathbf{L}cat)^{-1}(X)\simeq(\mathbf{L}cat)^{-1}(Y)$}}.
Let $X\simeq Y$ be two weakly equivalent flows in the q-model structure. Since $\moore$ is a right Quillen functor and since all flows are q-fibrant, we obtain the weak equivalence of Moore flows $\moore(X)\simeq \moore(Y)$. By definition of $(\mathbf{L}cat)^{-1}$ and since $\lmoore^{\mathcal{G}}$ is a left Quillen adjoint, we  deduce the sequence of isomorphisms and weak equivalences 
\[
(\mathbf{L}cat)^{-1}(X) \iso \lmoore^{\mathcal{G}} (\moore(X))^{cof} \simeq \lmoore^{\mathcal{G}} (\moore(Y))^{cof} \iso (\mathbf{L}cat)^{-1}(Y).
\]
The functors $(\mathbf{L}cat)$ and $(\mathbf{L}cat)^{-1}$ therefore induce functors between the homotopy categories.

$\bullet$ \underline{\textit{$(\mathbf{L}cat)^{-1}(\mathbf{L}cat)(X)\simeq X$}}. Let $X$ be a multipointed $d$-space. Then we have the sequence of isomorphisms and of weak equivalences
\[\begin{aligned}
(\mathbf{L}cat)^{-1}(\mathbf{L}cat)(X) & \iso \lmoore^{\mathcal{G}}\bigg(\moore\lmoore\overbracket{\moore^{\mathcal{G}}(X^{cof})}^{\substack{\hbox{\tiny q-cofibrant}\\ \hbox{\tiny by Corollary~\ref{pre-c2}}}}\bigg)^{cof}\\
&\simeq \lmoore^{\mathcal{G}}\bigg(\moore^{\mathcal{G}}(X^{cof})\bigg)^{cof}\\
&\simeq \lmoore^{\mathcal{G}}\moore^{\mathcal{G}}(X^{cof})\\
&\iso X^{cof}\\
&\simeq X,
\end{aligned}\]
the first isomorphism by definition of $(\mathbf{L}cat)$ and $(\mathbf{L}cat)^{-1}$ and by Theorem~\ref{decomposing}, the first weak equivalence since the adjunction $\lmoore\dashv \moore$ is a Quillen equivalence by \cite[Theorem~10.9]{Moore1} and since $\lmoore^{\mathcal{G}}$ is a left Quillen adjoint, the second weak equivalence by Corollary~\ref{pre-c2} and again since $\lmoore^{\mathcal{G}}$ is a left Quillen adjoint, the second isomorphism by Corollary~\ref{c1}, and the last weak equivalence by definition of a q-cofibrant replacement. 

$\bullet$ \underline{\textit{$(\mathbf{L}cat)(\mathbf{L}cat)^{-1}(Y)\simeq Y$}}. Let $Y$ be a flow. Then we have the sequence of isomorphisms and of weak equivalences
\[\begin{aligned}
(\mathbf{L}cat)(\mathbf{L}cat)^{-1}(Y) & \iso \big(\lmoore\moore^{\mathcal{G}}\big)\bigg(\lmoore^{\mathcal{G}}  (\mathbb{M}Y)^{cof}\bigg)^{cof}\\
&\simeq \big(\lmoore\moore^{\mathcal{G}}\big)\bigg(\lmoore^{\mathcal{G}}  (\mathbb{M}Y)^{cof}\bigg)\\
&\iso \lmoore(\mathbb{M}Y)^{cof}\\
&\simeq Y,
\end{aligned}\]
the first isomorphism by definition of $(\mathbf{L}cat)$ and $(\mathbf{L}cat)^{-1}$ and by Theorem~\ref{decomposing}, the first weak equivalence because $\moore^{\mathcal{G}}$ is a right Quillen adjoint, because $\lmoore^{\mathcal{G}}  (\mathbb{M}Y)^{cof}$ is q-cofibrant, because $\moore^{\mathcal{G}}$ preserves q-cofibrancy by Corollary~\ref{pre-c2} and finally because $\lmoore$ is a left Quillen adjoint, the second isomorphism by Theorem~\ref{c2}, and finally the last weak equivalence since the adjunction $\lmoore\dashv \moore$ is a Quillen equivalence by \cite[Theorem~10.9]{Moore1}. The proof is complete.
\epf

The underlying homotopy type of a flow is, morally speaking, the underlying space of states of a flow after removing the execution paths. It is defined only up to homotopy, not up to homeomorphism. We conclude the section and the paper by recovering it in a very intuitive way by using $(\mathbf{L}cat)^{-1}$. 

It is proved in \cite[Theorem~6.1]{4eme} that for every cellular flow $X$, there exists a cellular multipointed $d$-space $X^{top}$ such that there is an isomorphism $cat(X^{top})\iso X^{cof}$. 

\bd \cite[Section~6]{4eme}, The {\rm underlying homotopy type} of a flow $X$ is the topological space 
\[
||X|| := |X^{top}|
\]
where $|X^{top}|$ is the underlying topological space of the cellular multipointed $d$-space $X^{top}$. This yields a well defined functor \[||-||:\dtop\longrightarrow \ho(\top)\] from the category of flows to the homotopy category of the q-model structure of $\top$.
\ed

\bp \label{underlyinghomotopytype} For any flow $X$, there is the homotopy equivalence of topological spaces \[||X||\simeq |(\mathbf{L}cat)^{-1}(X)|.\]
\ep

\bpf
One has \[cat(\lmoore^{\mathcal{G}}(\moore(X)^{cof})) =(\lmoore\moore^{\mathcal{G}})(\lmoore^{\mathcal{G}}(\moore(X)^{cof})) \iso \lmoore(\moore(X)^{cof}),\] 
the equality by Theorem~\ref{decomposing} and the isomorphism by Theorem~\ref{c2}. Using the Quillen equivalence $\lmoore\dashv \moore$ of \cite[Theorem~10.9]{Moore1}, we obtain the weak equivalences of flows 
\[
cat(\lmoore^{\mathcal{G}}(\moore(X)^{cof})) \simeq X \simeq cat(X^{top}).
\]
Thanks to Theorem~\ref{lepb4}, we obtain the weak equivalence of q-cofibrant multipointed $d$-spaces 
\[
\lmoore^{\mathcal{G}}(\moore(X)^{cof})\simeq X^{top}.
\]
We deduce the homotopy equivalence between the underlying q-cofibrant spaces 
\[
|\lmoore^{\mathcal{G}}(\moore(X)^{cof})|\simeq |X^{top}|
\]
because the underlying topological space functor $|-|$ is a left Quillen functor by \cite[Proposition~8.1]{mdtop}.
\epf

The composite functor 
\[
\xymatrix@C=4em
{
\dtopG \fr{\lmoore^{\mathcal{G}}} & \ptop{\mathcal{G}} \fr{|-|} & \top
}
\]
is the composite of two left Quillen functors. Therefore the mapping \[X\mapsto |\lmoore^{\mathcal{G}}(X^{cof})|\] induces a functor from $\dtopG$ to $\ho(\top)$. For all $\mathcal{G}$-flows $X$, there is the isomorphism $X^{cof}\iso \moore^{\mathcal{G}}(\lmoore^{\mathcal{G}}(X^{cof}))$ by Theorem~\ref{c2}. Consequently, the functor $|\lmoore^{\mathcal{G}}((-)^{cof})|$ can be regarded as the underlying homotopy type functor for $\mathcal{G}$-flows.

\appendix

\section{The Reedy category \mins{\mathcal{P}^{u,v}(S)}: reminder}
\label{reedy}

Let $S$ be a nonempty set. Let $\mathcal{P}^{u,v}(S)$ be the small category defined by generators and relations as follows (see \cite[Section~3]{leftproperflow}): 
\begin{itemize}[leftmargin=*]
	\item $u,v\in S$ ($u$ and $v$ may be equal).
	\item The objects are the tuples of the form 
	\[\underline{m}=((u_0,\epsilon_1,u_1),(u_1,\epsilon_2,u_2),\dots ,(u_{n-1},\epsilon_n,u_n))\]
	with $n\geq 1$, $u_0,\dots,u_n \in S$, $\epsilon_1,\dots,\epsilon_n \in \{0,1\}$ and \[\forall i\hbox{ such that } 1\leq i\leq n, \epsilon_i = 1\Rightarrow (u_{i-1},u_i)=(u,v).\] 
	\item There is an arrow \[c_{n+1}:(\underline{m},(x,0,y),(y,0,z),\underline{n}) \to (\underline{m},(x,0,z),\underline{n})\]
	for every tuple $\underline{m}=((u_0,\epsilon_1,u_1),(u_1,\epsilon_2,u_2),\dots ,(u_{n-1},\epsilon_n,u_n))$ with $n\geq 1$ and every tuple $\underline{n}=((u'_0,\epsilon'_1,u'_1),(u'_1,\epsilon'_2,u'_2),\dots ,(u'_{n'-1},\epsilon'_{n'},u'_{n'}))$ with $n'\geq 1$. It is called a \textit{composition map}. 
	\item There is an arrow \[I_{n+1}:(\underline{m},(u,0,v),\underline{n}) \to (\underline{m},(u,1,v),\underline{n})\] for every tuple $\underline{m}=((u_0,\epsilon_1,u_1),(u_1,\epsilon_2,u_2),\dots ,(u_{n-1},\epsilon_n,u_n))$ with $n\geq 1$ and every tuple $\underline{n}=((u'_0,\epsilon'_1,u'_1),(u'_1,\epsilon'_2,u'_2),\dots ,(u'_{n'-1},\epsilon'_{n'},u'_{n'}))$ with $n'\geq 1$.
	It is called an \textit{inclusion map}. 
	\item There are the relations (group A) $c_i.c_j = c_{j-1}.c_i$ if $i<j$ (which means since $c_i$ and $c_j$ may correspond to several maps that if $c_i$ and $c_j$ are composable, then there exist $c_{j-1}$ and $c_i$ composable satisfying the equality). 
	\item There are the relations (group B) $I_i.I_j = I_j.I_i$ if $i\neq j$. By definition of these maps, $I_i$ is never composable with itself. 
	\item There are the relations (group C) \[c_i.I_j = \begin{cases}
	I_{j-1}.c_i&\hbox{if } j\geq i+2\\
	I_j.c_i&\hbox{if } j\leq i-1.
	\end{cases}\]
	By definition of these maps, $c_i$ and $I_i$ are never composable as well as $c_i$ and $I_{i+1}$. 
\end{itemize}

By \cite[Proposition~3.7]{leftproperflow}, there exists a structure of Reedy category on $\mathcal{P}^{u,v}(S)$ with the $\mathbb{N}$-valued degree map defined by \[d((u_0,\epsilon_1,u_1),(u_1,\epsilon_2,u_2),\dots ,(u_{n-1},\epsilon_n,u_n)) = n + \sum_i \epsilon_i.\]
The maps raising the degree are the inclusion maps. The maps decreasing the degree are the composition maps.

\section{An explicit construction of the left adjoint \mins{\lmoore^{\mathcal{G}}}}
\label{lmoore-explicit}

The proof of Theorem~\ref{adj-multi-moore} uses a well-known characterization of right adjoint functors between locally presentable categories. It is possible to describe more explicitly the functor $\lmoore^{\mathcal{G}}:\dtopG \to \ptop{\mathcal{G}}$ as follows. 

\begin{nota}
	The composite of two natural transformations $\mu:F\Rightarrow G$ and $\nu:G\Rightarrow H$ is denoted by $\nu \odot \mu$ to make the distinction with the composition of maps and of functors.
\end{nota}

The category $\mathcal{D}(\dtopG)$ of all small diagrams of Moore flows over all small categories is defined as follows. An object is a functor $F:{I}\to \dtopG$ from a small category ${I}$ to $\dtopG$. A morphism from $F:{I}_1\to \dtopG$ to $G:{I}_2\to \dtopG$ is a pair $(f:{I}_1\to {I}_2,\mu:F \Rightarrow G.f)$ where $f$ is a functor and $\mu$ is a natural transformation. If $(g,\nu)$ is a map from $G:{I}_2\to \dtopG$ to $H:I_3\to \dtopG$, then the composite $(g,\nu).(f,\mu)$ is defined by $(g.f,(\nu.f)\odot\mu)$. The identity of $F:{I}_1\to \dtopG$ is the pair $(\id_{{I}_1},\id_F)$. If $(h,\xi):(H:I_3\to \dtopG) \to (I:I_4\to \dtopG)$ is another map of $\mathcal{D}(\dtopG)$, then we have
\[\begin{aligned}
\left((h,\xi).(g,\nu)\right).(f,\mu)&=  (h.g,\xi.g \odot\nu).(f,\mu) \\
& = (h.g.f,\xi.g.f\odot\nu.f\odot\mu)\\
& = (h,\xi). (g.f,\nu.f\odot \mu)\\
& = (h,\xi) .\left((g,\nu).(f,\mu)\right).
\end{aligned}\]
Thus the composition law is associative and the category $\mathcal{D}(\dtopG)$ is well-defined. It is well-known that the colimit of small diagrams induces a functor $\liminj : \mathcal{D}(\dtopG) \to \dtopG$ (see e.g. \cite[Proposition~A.2]{leftproperflow}).

\bth \label{bigD}
There exists a functor $\mathbb{D}:\dtopG \to \mathcal{D}(\dtopG)$ such that the composite functor 
\[
\xymatrix@C=4em
{
	\dtopG \fr{\mathbb{D}} & \mathcal{D}(\dtopG) \fr{\liminj} & \dtopG
}
\]
is the identity functor and such that every vertex of $D(X)$ for any Moore flow $X$ is one of the following kind: 
\begin{enumerate}
	\item the Moore flow $\{0\}$,
	\item the globe $\globP(D)$ of some $\mathcal{G}$-space $D$,
	\item the Moore flow $\globP(D)*\globP(E)$ for two $\mathcal{G}$-spaces $D$ and $E$ where the final state of $\globP(D)$ is identified with the initial state of $\globP(E)$ (it is the ``composition'' of the two globes, hence the notation).
\end{enumerate}
Moreover, each $\mathcal{G}$-space $D$ and $E$ used by the diagram is of the form $\P_{\alpha,\gamma}X$ or $\P_{\alpha,\beta}X\ot \P_{\beta,\gamma}X$. 
\eth

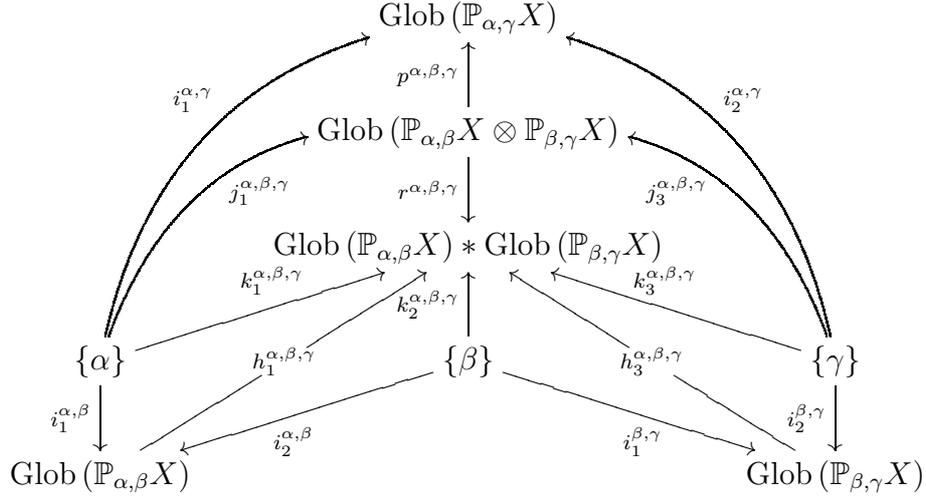
\begin{figure}
	\[
	\xymatrix{
		& \glob\left(\P_{\alpha,\gamma}X\right) & \\
		& \glob\left(\P_{\alpha,\beta}X\ot \P_{\beta,\gamma}X\right)\fu{p^{\alpha,\beta,\gamma}} \fd{r^{\alpha,\beta,\gamma}}&\\
		& \glob\left(\P_{\alpha,\beta}X\right)*\glob\left(\P_{\beta,\gamma}X\right)& \\
		{\{\alpha\}}\fd{i_1^{\alpha,\beta}}\ar@{->}[ru]^{k_1^{\alpha,\beta,\gamma}}\ar@/^35pt/[ruu]_{j_1^{\alpha,\beta,\gamma}}\ar@/^35pt/[ruuu]^{i_1^{\alpha,\gamma}}
		& \fu{k_2^{\alpha,\beta,\gamma}}
		{\{\beta\}}\ar@{->}[ld]^{i_2^{\alpha,\beta}}\ar@{->}[rd]_{i_1^{\beta,\gamma}}
		&
		\ar@{->}[lu]_{k_3^{\alpha,\beta,\gamma}}{\{\gamma\}\fd{i_2^{\beta,\gamma}}\ar@/_35pt/[luu]^{j_3^{\alpha,\beta,\gamma}}\ar@/_35pt/[luuu]_{i_2^{\alpha,\gamma}}}\\
		\glob\left(\P_{\alpha,\beta}X\right)\ar@{->}[ruu]|{h_1^{\alpha,\beta,\gamma}}
		& &
		\glob\left(\P_{\beta,\gamma}X\right)\ar@{->}[luu]|{h_3^{\alpha,\beta,\gamma}}
	}
	\]
	\caption{The Moore flow $X$ as a colimit of globes and points (the definition of the maps are easily understandable, cf. the proof of \cite[Theorem~6.1]{model3} for further explanations)}
	\label{UFO}
\end{figure}

\bpf
This theorem is proved in \cite[Theorem~6.1]{model3} for the category of flows which are semicategories enriched over the closed (semi)monoidal bicomplete category $(\top,\p)$ (see Definition~\ref{def-flow}). The diagram is depicted in Figure~\ref{UFO}. We refer to the proof of \cite[Theorem~6.1]{model3} for the definitions of the maps. Since a Moore flow is a semicategory enriched over the closed semimonoidal bicomplete category $([\mathcal{G}^{op},\top]_0,\ot)$, the proof is complete. 
\epf

\begin{nota}
	Denote by $\mathbb{B}(\dtopG)$~\footnote{$\mathbb{B}$ like ``brick'': the globes and the point are the elementary bricks to build flows.} the full subcategory of $\dtopG$ generated by $\{0\}$ and $\globP(D)$ where $D$ runs over the class of all $\mathcal{G}$-spaces.
\end{nota}

\bth \label{const-left-adj}
Let $\K$ be a bicomplete category. A functor \[F:\K \longrightarrow \dtopG\] has a left adjoint \[F_!:\dtopG \longrightarrow \K\] if and only if there exists a functor $m:\mathbb{B}(\dtopG) \to \K$ such that there are the natural bijections 
\[
\begin{aligned}
&\K(m(\{0\}),Y)\iso F(Y)^0\\
&\K(m(D),Y)\iso \dtopG(\globP(D),F(Y))
\end{aligned}
\]
for all objects $Y$ of $\K$ and all $\mathcal{G}$-spaces $D$.
\eth

\bpf
The ``only if'' direction comes from the fact that there is a natural bijection \[F(X)^0 \iso \dtopG(\{0\},F(X)).\] Let $D$ and $E$ be two $\mathcal{G}$-spaces. Let $m(\globP(D)*\globP(E))$ be the object of $\K$ defined by the pushout diagram of $\K$
\[
\xymatrix@C=4em@R=4em
{
	m(\{0\})  \fr{m(0 \mapsto 1)} \fd{m(0 \mapsto 0)} & m(D) \fd{}\\
	m(E) \fr{} & \cocartesien m(\globP(D)*\globP(E)).
}
\]
By taking the image by the functor $\K(-,Y):\K\to \set$, we obtain the pullback diagram of sets 
\[
\xymatrix@C=4em@R=4em
{
	\K(m(\globP(D)*\globP(E)),Y) \cartesien \fr{} \fd{} & \K(m(D),Y) \fd{}\\
	\K(m(E),Y) \fr{} &  \K(m(\{0\}),Y)
}
\]
for all objects $Y$ of $\K$. We therefore obtain the natural bijection of sets 
\[
\K(m(\globP(D)*\globP(E)),Y) \iso \dtopG(\globP(D)*\globP(E),F(Y))
\]
for all objects $Y$ of $\K$. Let
\[
F_!(X) := \liminj  m(\mathbb{D}(X)).
\]
This defines a functor from $\dtopG$ to $\K$. For all objects $Y$ of $\K$, there is the sequence of natural bijections (note that in the calculation below, the colimits are taken over a same small category which depends only on $X$)
\[\begin{aligned}
\K(F_!(X),Y) 
&\iso \K(\liminj  m(\mathbb{D}(X)),Y)\\
&\iso \limproj \K(m(\mathbb{D}(X)),Y) \\
&\iso \limproj \dtopG(\mathbb{D}(X),F(Y)) \\
&\iso \dtopG(\liminj \mathbb{D}(X),F(Y)) \\
&\iso \dtopG(X,F(Y)),
\end{aligned}\]
the first one by definition of $F_!$, the second one and the fourth one by the universal property of the (co)limit, the third one by hypothesis and by the calculation above, and finally the last one by Theorem~\ref{bigD}.
\epf

After Theorem~\ref{const-left-adj}, it suffices now to find a multipointed $d$-space denoted by $\mathbb{M}_!^{\mathcal{G}}(\{0\})$ such that there is a natural bijection with respect to $X$ \[\ptop{\mathcal{G}}\big(\mathbb{M}_!^{\mathcal{G}}(\{0\}),X\big) \iso \moore^{\mathcal{G}}(X)^0\]  and a multipointed $d$-space denoted by $\mathbb{M}_!^{\mathcal{G}}(\globP(D))$ natural with respect to the $\mathcal{G}$-space $D$ such that there is a natural bijection with respect to $D$ and $X$
\[\ptop{\mathcal{G}}\big(\mathbb{M}_!^{\mathcal{G}}(\globP(D)),X\big) \iso \dtopG(\globP(D),\moore^{\mathcal{G}}(X)).\] 

We have the natural bijections 
\[
\ptop{\mathcal{G}}\big(\{0\},X\big) \iso X^0 \iso  \moore^{\mathcal{G}}(X)^0,
\]
and therefore \[\mathbb{M}_!^{\mathcal{G}}(\{0\})=\{0\}.\] 

We have the sequence of natural bijections
\[
\begin{aligned}
\ptop{\mathcal{G}}\left(\int^{\ell}\glob^{\mathcal{G}}_\ell(D(\ell)),X\right) &\iso \int_{\ell} \ptop{\mathcal{G}}(\glob^{\mathcal{G}}_\ell(D(\ell)),X) \\
&\iso \int_{\ell} \bigsqcup_{(\alpha,\beta)\in X^0\p X^0} \top (D(\ell),\P_{\alpha,\beta}^\ell\moore^{\mathcal{G}}(X)) \\
&\iso \bigsqcup_{(\alpha,\beta)\in X^0\p X^0} \int_{\ell} \top (D(\ell),\P_{\alpha,\beta}^\ell\moore^{\mathcal{G}}(X)) \\
&\iso \bigsqcup_{(\alpha,\beta)\in X^0\p X^0}\top^{\mathcal{G}^{op}}(D,\P_{\alpha,\beta}\moore^{\mathcal{G}}(X)) \\
&\iso \bigsqcup_{(\alpha,\beta)\in X^0\p X^0}\topdgr_0(D,\P_{\alpha,\beta}\moore^{\mathcal{G}}(X))\\
&\iso \dtopG(\globP(D),\moore^{\mathcal{G}}(X)),
\end{aligned}
\]
the first bijection by definition of a (co)limit, the second bijection by Proposition~\ref{precalcul}, the third bijection because the underlying diagram of this end is connected, the fourth bijection by \cite[page 219 (2)]{MR1712872}, the fifth bijection since $\topdgr_0$ is a full subcategory of $\top^{\mathcal{G}^{op}}$, and finally the last bijection by Proposition~\ref{map-from-glob}. We obtain 
\[
\mathbb{M}_!^{\mathcal{G}}(\globP(D)) = \int^{\ell}\glob^{\mathcal{G}}_\ell(D(\ell)).
\]

\section{The setting of \mins{k}-spaces}
\label{kspaces}

In this appendix, the category of $k$-spaces is denoted by $\top_k$ and the category of $\Delta$-generated spaces by $\top_\Delta$. The proofs are just sketched.

We must at first prove the existence of the projective q-model structure of $[\mathcal{G}^{op},\top_k]_0$: \cite{dgrtop} is written in the locally presentable setting. We do not know whether the arguments of \cite{PiaDgr} are valid here since they are written in the category of Hausdorff $k$-spaces. Anyway, it is possible to give a much simpler argument. The inclusion $\top_\Delta\subset \top_k$ has a right adjoint $k_\Delta:\top_k\to \top_\Delta$, which gives rise to a Quillen equivalence. In fact, the q-model structure of $\top_k$ is right-induced by $k_\Delta:\top_k\to \top_\Delta$ in the sense of \cite{HKRS17,GKR18}. From the functor $k_\Delta$, we obtain a right adjoint $U:[\mathcal{G}^{op},\top_k]_0 \to [\mathcal{G}^{op},\top_\Delta]_0$. For an object $X$ of $[\mathcal{G}^{op},\top_k]_0$, let $\cocyl(X)=\ell \mapsto \ttop_k([0,1],X(\ell))$ where $\ttop_k$ is the internal hom of $\top_k$. Since the composite functor $k_\Delta.\ttop_k$ is the internal hom of $\top_\Delta$, the Quillen path object argument can be used to obtain that $U:[\mathcal{G}^{op},\top_k]_0 \to [\mathcal{G}^{op},\top_\Delta]_0$ right-induces the projective q-model structure of $[\mathcal{G}^{op},\top_k]_0$. This technique still works in the cofibrantly generated non-combinatorial setting: it dates back to \cite{MR36:6480} (see also \cite[Theorem~11.3.2]{ref_model2}).

The q-model category of multipointed $d$-spaces is constructed in \cite[Theorem~6.14]{QHMmodel} by right-inducing it and by using the Quillen path object argument again. The q-model category of $\mathcal{G}$-flows is constructed in \cite[Theorem~8.8]{Moore1} by mimicking the method used in \cite[Theorem~3.11]{leftdetflow} which works for any convenient category of topological spaces. Indeed, it uses Isaev's work \cite{Isaev} about model categories of fibrant objects which does not require to work in a locally presentable setting. 

Theorem~\ref{pre-right-adj} is not valid anymore. See \cite[Theorem~5.10]{leftproperflow} for a treatment of the similar problem for flows. The reason is that a $k$-space is not necessarily homeomorphic to the disjoint sum of its path-connected components (e.g. the Cantor space). It is used in Theorem~\ref{adj-multi-moore} together with the local presentability of the category of $\Delta$-generated spaces to prove the existence of the left adjoint $\lmoore^{\mathcal{G}}: \dtopG \to \ptop{\mathcal{G}}$. In the framework of $k$-spaces, the existence of $\lmoore^{\mathcal{G}}$ can be proved using Appendix~\ref{lmoore-explicit}. We have therefore a Quillen adjunction \[\lmoore^{\mathcal{G}}\dashv \moore^{\mathcal{G}}:\dtopG(\top_k) \leftrightarrows \ptop{\mathcal{G}}(\top_k)\]
between the q-model structures, where the notation $(\top_k)$ is to specify the category of topological spaces which is used. 

A $k$-space, unlike a $\Delta$-generated space, is not necessarily sequential. It is not clear how to adapt the proof of Theorem~\ref{img-closed} by replacing sequences by nets since there is a Cantor diagonalization argument which does not seem to be adaptable at least with uncountable nets. It is not clear either how to modify accordingly the last part of the proof of Theorem~\ref{pre-calculation-pathspace2} about the sequential continuity because $\widehat{U_{\underline{c}}}$ is an arbitrary compact space now, and not necessarily a closed subset of $[0,1]$ anymore. To obtain Theorem~\ref{synth} for $k$-spaces, another method must be used. For all $k$-spaces $Z$, the canonical map $k_\Delta(Z)\to Z$ is a weak homotopy equivalence which induces a bijection $\top_k([0,1],k_\Delta(Z)) \iso \top_\Delta([0,1],Z)$ because $[0,1]\in \top_\Delta\subset \top_k$. From these observations, we obtain two left Quillen equivalences $\dtopG(\top_\Delta)\to \dtopG(\top_k)$ and $\ptop{\mathcal{G}}(\top_\Delta)\to \ptop{\mathcal{G}}(\top_k)$. We obtain a diagram of left Quillen adjunctions 
\[
\xymatrix@C=3em@R=3em
{
	\dtopG(\top_\Delta) \fr{\subset} \fd{} & \dtopG(\top_k) \fd{}\\
	\ptop{\mathcal{G}}(\top_\Delta) \fr{\subset}  & \ptop{\mathcal{G}}(\top_k)
}
\]
which is commutative because the $\mathcal{G}$-flows $\{0\}$ and $\glob(D)$ of $\dtopG(\top_\Delta)$ are taken to the same multipointed $d$-space of $\ptop{\mathcal{G}}(\top_k)$ by Appendix~\ref{lmoore-explicit}. Using the two-out-of-three property, we obtain that the Quillen adjunction \[\lmoore^{\mathcal{G}}\dashv \moore^{\mathcal{G}}:\dtopG(\top_k) \leftrightarrows \ptop{\mathcal{G}}(\top_k)\] is a Quillen equivalence. Let $X$ be a q-cofibrant object of $\ptop{\mathcal{G}}(\top_k)$. Then $X\in \ptop{\mathcal{G}}(\top_\Delta)$. By Corollary~\ref{c1}, there is the isomorphism $\lmoore^{\mathcal{G}}(\moore^{\mathcal{G}}(X))\iso X$. Let $X$ be a q-cofibrant object of $\dtopG(\top_k)$. Then $X\in \dtopG(\top_\Delta)$. By Theorem~\ref{c2}, there is the isomorphism $X \iso \moore^{\mathcal{G}}(\lmoore^{\mathcal{G}}(X))$. 

At this point, it is legitimate to ask whether the main results of the companion paper \cite{Moore1} are valid for $k$-spaces. The answer is that they are. The main tool of \cite{Moore1} is the Quillen equivalence $\liminj:[\mathcal{G}^{op},\top_\Delta]_0 \leftrightarrows \top_\Delta:\Delta_{\mathcal{G}^{op}}$ proved in \cite[Theorem~7.6]{dgrtop}. There is the commutative diagram of left Quillen adjoints 
\[
\xymatrix@C=3em@R=3em
{
	[\mathcal{G}^{op},\top_\Delta]_0 \fr{} \fd{} & \top_\Delta \fd{}\\
	[\mathcal{G}^{op},\top_k]_0 \fr{}  & \top_k.
}
\]
All left Quillen adjoints except maybe the bottom horizontal one are left Quillen equivalences. Therefore the bottom horizontal one is a left Quillen equivalence as well.

As a conclusion, most of the results, but not all, of this paper and of the companion paper \cite{Moore1} are still valid for $k$-spaces. However, there is no known proofs avoiding to use $\Delta$-generated spaces.


\end{document}